\documentclass[12pt]{amsart}

\usepackage[text={400pt,650pt},centering]{geometry}

\usepackage{color}
\usepackage{esint,amssymb}
\usepackage{graphicx}
\usepackage{MnSymbol}
\usepackage{mathtools}
\usepackage[colorlinks=true, pdfstartview=FitV, linkcolor=blue, citecolor=blue, urlcolor=blue,pagebackref=false]{hyperref}
\usepackage{microtype,dsfont}

\newtheorem{theorem}{Theorem}
\newtheorem{proposition}[theorem]{Proposition}
\newtheorem{lemma}[theorem]{Lemma}
\newtheorem{corollary}[theorem]{Corollary}
\theoremstyle{definition}
\newtheorem{remark}[theorem]{Remark}
\newtheorem{example}{Example}
\newtheorem*{definition}{Definition}

\numberwithin{equation}{section}
\numberwithin{theorem}{section}

\newcommand{\Z}{\mathbb{Z}}
\newcommand{\N}{\mathbb{N}}
\newcommand{\R}{\mathbb{R}}

\newcommand{\E}{\mathbb{E}}
\renewcommand{\P}{\mathbb{P}}
\newcommand{\F}{\mathcal{F}}
\newcommand{\G}{\mathcal{G}}
\renewcommand{\H}{\mathcal{H}}
\newcommand{\K}{\mathcal{K}}

\newcommand{\Rd}{\mathbb{R}^d}

\newcommand{\ep}{\varepsilon}

\renewcommand{\subset}{\subseteq}
\newcommand{\indc}{\mathds{1}}

\DeclareMathOperator{\dist}{dist}

\DeclareMathOperator*{\esssup}{ess\,sup}
\DeclareMathOperator{\diam}{diam}
\DeclareMathOperator{\tr}{tr}

\DeclareMathOperator{\USC}{USC}
\DeclareMathOperator{\LSC}{LSC}
\DeclareMathOperator{\BUC}{BUC}

\renewcommand{\bar}{\overline}
\renewcommand{\tilde}{\widetilde}

\newcommand{\data}{\mathrm{data}}

\begin{document}

\title[Homogenization of geometric motions]{Stochastic homogenization of quasilinear Hamilton-Jacobi equations and geometric motions}

\begin{abstract}
We study random homogenization of second-order, degenerate and quasilinear Hamilton-Jacobi equations which are positively homogeneous in the gradient. Included are the equations of forced mean curvature motion and others describing geometric motions of level sets 
as well as a large class of viscous, non-convex Hamilton-Jacobi equations. The main results include the first proof of qualitative stochastic homogenization for such equations.  We also present quantitative error estimates which give an algebraic rate of homogenization.
\end{abstract}

\author[S. Armstrong]{Scott Armstrong}
\address[S. Armstrong]{Ceremade (UMR CNRS 7534) \\ Universit\'e Paris-Dauphine \\ Place du Mar\'echal De Lattre De Tassigny \\ 75775 Paris CEDEX 16, France}
\email{armstrong@ceremade.dauphine.fr}

\author[P. Cardaliaguet]{Pierre Cardaliaguet}
\address[P. Cardaliaguet]{Ceremade (UMR CNRS 7534) \\ Universit\'e Paris-Dauphine \\ Place du Mar\'echal De Lattre De Tassigny \\ 75775 Paris CEDEX 16, France}
\email{cardaliaguet@ceremade.dauphine.fr}

\keywords{stochastic homogenization, mean curvature equation, Hamilton-Jacobi equation, error estimate}
\subjclass[2010]{35B27, 35F21, 60K35}
\date{\today}

\maketitle

\section{Introduction}

\subsection{Motivation and informal summary of results}
In this paper, we study time-dependent, quasilinear, viscous Hamilton-Jacobi equations taking the form
\begin{equation}
\label{e.pde}
\partial_t u^\ep - \ep \tr\left( A\left(\frac{Du^\ep}{\left|Du^\ep\right|},\frac x\ep\right)D^2u^\ep \right) + H\left(Du^\ep,\frac x\ep \right) = 0 \quad \mbox{in} \ \Rd\times (0,\infty).
\end{equation}
We give a brief summary of the main assumptions, which are given precisely in~Section~\ref{s.statements}, below. 
The diffusion matrix $A(e,x)$ is assumed to be nonnegative definite for each $(e,x) \in \partial B_1\times\Rd$; in particular, the diffusive term may vanish or be the Laplacian. The Hamiltonian $H(\xi,x)$ is assumed to be positively homogeneous of order $p\in[1,\infty)$ in $\xi$, but is not necessarily convex in $\xi$. Both $A(e,\cdot)$ and $H(\xi,\cdot)$ are assumed to be stationary random fields sampled by a probability measure~$\P$ which satisfies a finite range of dependence.

\smallskip

The interest is in describing the behavior of solutions of~\eqref{e.pde} for $0<\ep\ll 1$. The main result is a characterization of the limit of $u^\ep(x,t)$, subject to suitable initial conditions, as $\ep\to 0$. We show that $u^\ep$ converges locally uniformly, with $\P$--probability one, to the solution $u$ of a deterministic equation of the form
\begin{equation} 
\label{e.homogenized}
\partial_t u+ \overline H(Du) = 0 \quad \mbox{in} \ \Rd \times (0,\infty),
\end{equation}
with an effective Hamiltonian~$\overline H:\Rd\to\R$ which has sublevel sets that are star-shaped with respect to the origin. This is the first stochastic homogenization result for a viscous equation with a Hamiltonian which may be nonconvex in~$Du^\ep$, or for a diffusion matrix which may have dependence on $Du^\ep$. 

A particular case of~\eqref{e.pde} satisfying our assumptions is the \emph{equation of forced mean curvature motion}
\begin{equation}
\label{e.MCM}
\partial_t u^\ep - \ep \tr\left( \left( I_d - \frac{Du^\ep\otimes Du^\ep}{\left|Du^\ep\right|^2} \right)D^2u^\ep \right) + a\left(\frac x\ep\right) \left| Du^\ep\right| = 0 \quad \mbox{in} \ \Rd\times (0,\infty),
\end{equation}
where the forcing field~$a$ is positive, Lipschitz, bounded and satisfies
\begin{equation}
\label{e.MCMcondition}
 \inf_{x\in\Rd} \left( a^2(x) - (d-1) \left|Da(x) \right| \right) >0  \quad \mbox{$\P$--a.s.} 
\end{equation}
The level sets of solutions of~\eqref{e.MCM} follow a generalized evolution with normal velocity $\ep \kappa + a\left( \tfrac x\ep \right)$, where $\kappa$ represents the mean curvature of the surface. 
In this context, the limiting homogenized equation takes the form
\begin{equation}
\label{e.MCMbar}
\partial_t u + \overline a\left( \frac{Du}{|Du|} \right)|Du| = 0,
\end{equation}
where $\overline a:\partial B_1 \to \R$ is a positive function which describes the velocity of level sets of $u$ and thus the effective velocity of the original flow. The homogenization of~\eqref{e.MCM} in the random setting has been an open problem for some time; its importance was recently highlighted in the recent review article~\cite{C} (see page~773). The condition~\eqref{e.MCMcondition} was introduced in~\cite{LS0} in the context of periodic homogenization of~\eqref{e.MCM} and its role is to ensure the Lipschitz regularity of solutions (roughly speaking, it is the condition under which the Bernstein method for estimating $\|Du^\ep\|_{L^\infty}$ is applicable). It was recently shown in~\cite{CM} to be necessary for homogenization to hold, in general, even in that much simpler context: in other words, without Lipschitz regularity, homogenization may fail. 

\smallskip

Another special case of~\eqref{e.pde} includes the general class of viscous Hamilton-Jacobi equations of the form 
\begin{equation}
\label{e.viscous}
\partial_t u^\ep - \ep \Delta u^\ep + H\left( Du^\ep,\frac x\ep \right) = 0 \quad \mbox{in} \ \Rd \times (0,\infty), 
\end{equation}
where $H$ satisfies, for some $p>1$ and $0<c_0\leq C_0$,
\begin{equation*}
H(t e ,x) = t^p H(e,x) \quad \mbox{and} \quad c_0 \leq H(e,x) \leq C_0, \quad \forall t\geq 0, \ e\in\partial B_1, \ x\in\Rd.
\end{equation*}
Our results therefore give the first large class of non-convex Hamilton-Jacobi equations for which homogenization holds in $d>1$. Even for~\eqref{e.viscous} with the second-order term removed the conclusions are new: the only previous results even for first-order nonconvex equations to our knowledge are found in the recent papers~\cite{ATY1,ATY2}. The latter articles demonstrated homogenization in one space dimension and treated special cases in higher dimensions for Hamiltonians taking the form $H(p,x) = \tilde H(p)+ W(x)$, an entirely different structure from the one here. 

\smallskip

Every proof of qualitative homogenization for a Hamilton-Jacobi equation in the random setting has been based in some way on an application of the subadditive ergodic theorem. This requires the identification of a subadditive quantity whose limiting behavior controls that of the solutions to the equation. Such subadditive structures have only been found, with the exception of the results in~\cite{ATY1,ATY2}, in the case of convex or quasi-convex Hamiltonians and equations with linear diffusion terms. For this reason, a general qualitative theory of stochastic homogenization for equations with nonconvex Hamiltonians or quasilinear viscous terms has proved elusive. Indeed, even identifying a single example for which it can be proven that homogenization holds for a viscous equation with nonconvex or quasilinear structure has remained  open until now. 

\smallskip

In this paper, we propose a new strategy for obtaining qualitative homogenization results, based on the simple idea that the lack of a subadditive structure can be overcome by a quantitative approach. Rather than using soft arguments based on ergodic theorems, we assume a much stronger mixing assumption for the coefficients (a finite range of dependence condition) and attempt to prove more: homogenization with an explicit error estimate. The quantitative theory of stochastic homogenization for Hamilton-Jacobi equations originated in~\cite{ACS} for first-order equations and in~\cite{AC} for semi-linear viscous equations. The strategy here is to build on the techniques introduced in~\cite{ACS,AC} to handle more general equations. While the arguments in those papers seem to still rely on subadditivity, we demonstrate here that the ideas in fact do not require it.

\smallskip

Very recently (and several months after this paper was written), Ziliotto~\cite{Zil} produced an example of a first-order, coercive (and nonconvex) Hamilton-Jacobi equation with stationary-ergodic coefficients for which homogenization fails. Thus, in addition to negatively resolving the question of whether a general homogenization result holds for nonconvex Hamilton-Jacobi equations, Ziliotto demonstrated conclusively the impossibility of obtaining homogenization results by soft or qualitative arguments. This progress provides further motivation for developing quantitative approaches to this and similar problems. 

\smallskip

The assumption of finite range dependence is well-motivated physically and is analogous to the standard assumption of~i.i.d. in discrete probability models. It is not an assumption made for simplicity: we do not know how to relax it even to allow very quick decaying correlations of the coefficients. However, by stability arguments, we can obtain homogenization results for coefficients fields which are uniform limits of finite-range fields. This covers many typical examples, including for instance coefficients fields built by convolutions of smooth (but not compactly supported) functions against Poisson point clouds. 
 
 \smallskip
 
 On a technical level, as we will see in Section~\ref{s.fluctuations}, the finite range assumption gives us \emph{almost sure} bounds on the increments of a certain martingale, so that the fluctuations of this martingale can be strongly controlled by Azuma's inequality. An assumption which allows for long-range correlations would not have this property. On a more philosophical level, our strategy is quantitative and therefore requires a quantitative ergodic assumption-- but it turns out that the finite range of dependence condition happens to be the only assumption under which it is known how to prove quantitative results. Indeed, it is completely open to obtain quantitative results for the homogenization of even first-order Hamilton-Jacobi equations under any assumption that allows for long-range correlations. Even obtaining a convergence rate for the shape theorem in first-passage percolation is well-known to be open when the edges do not satisfy a finite range of dependence! 

\smallskip

In the case of the forced mean curvature equation~\eqref{e.MCM} with periodic coefficients, a lot of attention has been given to the existence of \emph{plane-like} solutions, which began with the work of~\cite{CD}. In our setting, these are solutions, for a given $\mu>0$ of the stationary problem 
\begin{equation*}
 -\tr\left( A(Du,x) D^2u\right) + H(Du,x) = \mu \quad \mbox{in}  \ \Rd
\end{equation*}
whose graphs stay within a bounded distance from an affine function. Plane-like solutions are intimately connected to homogenization; it is not hard to see that their existence for every given slope implies homogenization for general initial data (at least in the context in which one has Lipschitz solutions for sufficiently smooth initial data, which is assured by our hypotheses). In the random setting, we do not expect that plane-like solutions exist, in general. Indeed, their existence would imply a rate of homogenization of $O(\ep)$ for affine initial data, while in dimension $d=2$ we expect a convergence rate of $O(\ep^{\frac23})$, in line with the conjectured bound for first-passage percolation (a discrete version of a first-order Hamilton-Jacobi equation) and other growing surface models expected to possess a scaling limit related to the KPZ equation.

\smallskip

Our strategy for homogenization is nevertheless based on the construction of a weaker variant of a plane-like solution in half spaces. We consider the problem posed in the half-space $\H_e^+ = \{ x\cdot e>0\}$, for a given parameter~$\mu>0$ and unit vector~$e$:
\begin{equation*} \label{}
\left\{ 
\begin{aligned}
&  -\tr\left( A(Dm_\mu,x) D^2m_\mu\right) + H(Dm_\mu,x) = \mu & \mbox{in} & \ \H_e^+, \\
& m_\mu = 0 & \mbox{on} & \ \partial \H_e^+
\end{aligned}
\right.
\end{equation*}
We call this the \emph{planar metric problem}, as the value of $m_\mu(x)$ can be thought of as a ``distance"  from the point $x$ to the plane $\partial \H_e^+$. Rather than prove that $m_\mu$ stays a bounded distance from an affine function, it turns out to be sufficient for homogenization to show roughly that, for some exponent $\alpha>0$ and a deterministic constant~$\overline m_\mu(e)$, we have
\begin{equation*} \label{}
\sup_{ x\in B_r\cap \H_e^+} \frac1r \left| m_\mu (x) - \overline m_\mu(e) (x\cdot e) \right| = O\left( r^{-\alpha} \right) \quad \mbox{as} \ r \to \infty \quad \mbox{with high probability.}
\end{equation*}
In other words, we may soften the requirement that $m_\mu$ be a bounded distance from a plane by allowing the permitted distance to depend on the distance from the boundary plane~$\partial \H_e^+$ and by only checking points in a bounded set. The heart of the paper is the proof of this estimate, which is stated precisely in Proposition~\ref{p.planarquantitative} and proved in Sections~\ref{s.fluctuations} and~\ref{s.means}. The argument is naturally split into two steps: first we show by a concentration argument, exploiting the finite range dependence of the coefficients, that the random fluctuations of $m_\mu(x)$ are at most of order $O((x\cdot e)^{\frac12})$. Then we argue that the means $\E\left[ m_\mu(x) \right]$, which by stationarity can be written as a deterministic function of the form $f(x\cdot e)$, must therefore be close to a plane, that is, $f$ is close to a linear function. For the last step we think of the distance to the boundary plane as ``time" and consider the ``semigroup" generated by the ``flow" and use a maximum principle argument. The fact that $\partial \H_e^+$ is unbounded and the fluctuations estimate is not uniform raises another difficulty which is overcome by a new ``approximate finite speed of propagation" property presented in Section~\ref{s.propagation}.  

\subsection{Precise assumptions and statement of homogenization}
\label{s.statements}
We begin with the structural conditions on the coefficients in~\eqref{e.pde} before giving the probabilistic formulation of the problem. Throughout that paper, we fix an exponent $p\geq 1$, a dimension $d\in\N_*$, a positive integer $n\in\N_*$, constants $0<c_0 \leq C_0 <\infty$ and a modulus $\rho \in C([1,\infty))$ satisfying 
\begin{equation*} \label{}
\lim_{R\to \infty} \rho(R) = \infty. 
\end{equation*}
We also fix parameters $\theta,\kappa >0$. It is convenient to set
\begin{equation*} \label{}
\data:= \left(d,p,n,c_0,C_0,\rho,\theta,\kappa \right). 
\end{equation*}

\smallskip

We consider diffusion matrices $A: \partial B_1 \times \Rd \to \R^{d\times d}$ which take the form 
\begin{equation}
\label{e.Asigma}
A = \frac12 \sigma\sigma^T, \qquad \mbox{where $\sigma\in C^{1}(\partial B_1 \times\Rd; \, \R^{d\times n})$}
\end{equation}
is a matrix-valued function (with $\sigma^T$ denoting its transpose) which satisfies 
\begin{equation}
\label{e.sigma}
\left| \sigma(e,x) \right| + \left| D_x \sigma(e,x) \right| + \left| D_\xi \sigma(e,x) \right| \leq C_0\qquad \mbox{in} \ \partial B_1\times \R^d.
\end{equation}
For notational purposes, it is convenient to extend $\sigma(\cdot,x)$ (and $A(\cdot,x)$) to $\Rd\setminus\{0\}$ by defining $\sigma(\xi,x):= \sigma(\xi/|\xi|,x)$, which makes them $0$-homogeneous functions of their first argument. 
The Hamiltonian 
\begin{equation}
\label{e.H}
H \in C^{1}( \Rd \times\Rd)
\end{equation}
 is assumed to satisfy, for every $t>0$ and $\xi,x\in\Rd$, 
\begin{equation}
\label{e.Hposhom}
H(t\xi,x) = t^pH(\xi,x)  \qquad \mbox{and} \qquad 
c_0\left|\xi\right|^p \leq H(\xi,x) \leq C_0\left|\xi\right|^p
\end{equation}
and
\begin{equation} 
\label{e.HLip}
\left| D_x H( \xi,x) \right| + \left| \xi \right| \left| D_\xi H(\xi,x) \right| \leq C_0 \left| \xi \right|^p\qquad \mbox{in} \ \left(\R^d\setminus\{0\}\right)\times \R^d.
\end{equation}
The main structure condition on the coefficients is what we call the \emph{Lions-Souganidis (LS) coercivity condition}, since it was introduced (albeit in a slightly different form) in~\cite{LS0}: we suppose that, with the modulus $\rho$ fixed as above, 
\begin{equation}\label{e.LScoercivity}
\inf\left\{ \mathcal C_{\sigma,H}(\xi,x) \,:\,  \xi,x\in\Rd,  \ \left|\xi\right| \geq R  \right\} \geq \rho(R),
\end{equation}
where the quantity $\mathcal C_{\sigma,H}$, which measures the coercivity of the equation, is defined by 
\begin{multline*} \label{}
\mathcal C_{\sigma,H}(\xi,x) := \inf_{\eta \in B_\kappa(\xi)} \bigg( \theta(1-2\theta)  \left( H(\eta,x)\right)^2 
-(1+\kappa)^3\left|\sigma(\eta,x)\right|^2\left|D_x\sigma(\eta,x)\right|^2|\xi|^2  \\
-\theta (1+\kappa)^2|\sigma(\eta,x)|^2 \left| \xi \right| \left( \left|D_xH(\eta,x)\right| + \kappa\left|D_\xi H(\eta,x)\right|\right) 
 \bigg).
\end{multline*}
At first glance,~\eqref{e.LScoercivity} appears to be quite a technical assumption. Let us mention that it is redundant in the first-order case ($A=0$) or in the case that the Hamiltonian grows faster than linearly ($p>1$). The full condition is necessary to allow for $p=1$ in the viscous setting, which includes the forced mean-curvature model. We show below in Section~\ref{s.assumptions} that~\eqref{e.LScoercivity} is satisfied by each of our motivating examples. It is a generalization of the condition~\eqref{e.MCMcondition} for the forced mean curvature equation, which, as mentioned above, has been shown to be necessary for homogenization, even in the periodic case~\cite{CM}. The main role of~\eqref{e.LScoercivity} is to ensure that Lipschitz estimates hold for solutions of the equations we consider (this is proved in Appendix~\ref{s.metricproblem}). The reason it is so technical is because it is the condition needed to ensure that Bernstein's method for obtaining gradient bounds is applicable.

\smallskip

We work with the probability space $\Omega$, defined to be the set of all such coefficient fields:
\begin{equation*}
\Omega:= \left\{ 
(\sigma,H) \,:\, 
\mbox{$\sigma$ and $H$ satisfy~\eqref{e.Asigma},~\eqref{e.sigma},~\eqref{e.H},~\eqref{e.Hposhom},~\eqref{e.HLip} and~\eqref{e.LScoercivity}}
 \right\}.
\end{equation*}
We endow $\Omega$ with a family $\{ \F(U) \}$ of $\sigma$--algebras, indexed by the family of Borel subsets $U$ of $\Rd$, and defined by
\begin{multline}
\label{e.FU}
\F(U):= \mbox{the $\sigma$--algebra generated by the family of maps $\Omega \to \R^{d\times n} \times \R$,} \\ 
(\sigma,H) \mapsto (\sigma(e,x),H(\xi,x)), \quad\mbox{where} \ e\in\partial B_1, \ \xi\in\Rd \ \mbox{and} \ x\in U.
\end{multline}
The largest of these we denote by $\F:=\F(\Rd)$. The interpretation of~$\F(U)$ is that it is the $\sigma$--algebra containing ``all of the information which can be obtained by observing the coefficients restricted to $U$." 

\smallskip

Throughout the paper, we consider a fixed probability measure~$\P$ on $(\Omega,\F)$ which satisfies the following two conditions:
\begin{enumerate}

\item[(P1)] $\P$ is stationary: for every $y\in\Rd$, we have that $$\P = \P \circ T_y,$$ where $T_y:\Omega\to \Omega$ acts on $\Omega$ by translation in $y$, i.e., 
\begin{equation*}
T_y(\sigma,H):= ((e,z)\mapsto \sigma(e,z+y), (\xi,z)\mapsto H(\xi,z+y)).
\end{equation*}

\item[(P2)] $\P$ satisfies a unit range of dependence: for all Borel subsets $U,V\subseteq \Rd$ such that $\dist(U,V) \geq 1$, we have
\begin{equation*}
\mbox{$\F(U)$ and $\F(V)$ are $\P$--independent.}
\end{equation*}

\end{enumerate}

\smallskip

Throughout the paper, all differential inequalities are to be understood in the viscosity sense. Since the quasilinear diffusions we consider have singularities at $\xi=0$, and for the readers' convenience, we recall the appropriate definitions in Section~\ref{s.notation} below.

\smallskip

The main result of the paper is Theorem~\ref{t.main}, stated at the beginning of Section~\ref{s.outline}. Here we present a consequence of it which is simpler, more qualitative and thus easier to read. It asserts that the initial value problem for~\eqref{e.pde} homogenizes almost-surely at an algebraic rate, at least for Lipschitz continuous solutions (the Lipschitz assumption can be removed, see Remark~\ref{r.lipschitz} below). 

\begin{theorem}
\label{t.homogenization}
Let~$\P$ be a probability measure on~$(\Omega,\F)$ which satisfies~(P1) and~(P2). Then there exists a universal exponent $\alpha>0$ and a function $\overline H\in C^{0,\frac 27-}_{\mathrm{loc}}(\Rd)$ satisfying, for every $\xi\in\Rd$,
\begin{equation*}
c_0 \left| \xi \right|^p \leq \overline H(\xi) \leq C_0 \left| \xi \right|^p
\end{equation*}
such that, for every~$T\ge1$ and~$u,u^\ep \in W^{1,\infty} (\Rd \times [0,T])$ satisfying
\begin{equation*}
\left\{ 
\begin{aligned}
& \partial_t u^\ep -\ep \tr\left( A\left(Du^\ep,\frac x\ep \right)D^2u^\ep \right) + H\left(Du^\ep,\frac x\ep \right) = 0 &  \mbox{in} & \ \Rd\times (0,T], \\
& \partial_t u +  \overline H(Du) = 0 & \mbox{in} & \ \Rd\times (0,T], \\
& u^\ep(\cdot,0) = u(\cdot,0) & \mbox{on} & \ \Rd,
\end{aligned}
\right.
\end{equation*}
we have 
\begin{equation*}
 \P \left[ \sup_{R\geq1} \limsup_{\ep \to 0} \, \sup_{(x,t) \in B_R\times [0,T]}  \ep^{-\alpha} \left| u^\ep(x,t)  - u(x,t) \right| =  0 \right]  = 1.
\end{equation*}
\end{theorem}

The proof of Theorem~\ref{t.homogenization} is given at the end of Section~\ref{s.outline}. 

\smallskip

As we will see from the argument, the exponent~$\alpha$ in Theorem~\ref{t.homogenization} can be taken to be any positive number smaller than~$\frac{1}{90}$. Needless to say, this is not optimal, and we made no attempt to optimize our proof to obtain the best exponent. It is less obvious that the limiting H\"older  exponent of $\tfrac 27$ for the regularity of $\overline H$ can be improved, at least in the general quasilinear setting (it is easy to show that $\overline H$ is Lipschitz in the semilinear case). Note that in the statement of the theorem we used the notation $C^{0,\beta-}_{\mathrm{loc}}(\Rd) := \cap_{0<\gamma<\beta} C^{0,\gamma}_{\mathrm{loc}}(\Rd)$.

\subsection{Examples}
\label{s.assumptions}

In this subsection, we check that the motivating examples~\eqref{e.MCM} and~\eqref{e.viscous} satisfy the (LS) condition~\eqref{e.LScoercivity}.

\begin{example}[Superlinear case] 
The assumption~\eqref{e.LScoercivity} is redundant in the case~$p>1$, that is, we may remove the (LS) condition in the case of a superlinear Hamiltonian. Therefore our results apply to a large family of quasilinear, viscous and, in general, nonconvex Hamilton-Jacobi equations. 

\smallskip

To check this, it is enough to show that, if $\mathcal C_{\sigma,H}(\xi,x)\leq  K$ for some $K\geq 1$, then $|\xi|\leq R_K$ for a sufficiently large real number $R_K$. The condition $\mathcal C_{\sigma,H}(\xi,x) \leq K$ implies that there exists $\eta\in B_\kappa(\xi)$ such that
\begin{multline*}
K \geq \theta(1-2\theta) c_0^2|\eta|^{2p} -\theta (1+\kappa)^2|\sigma|^2C_0 (1+|\eta|^p) |\xi|
 -(1+\kappa)^3|\sigma|^2|D_x\sigma|^2|\xi|^2 \\ -\theta \kappa (1+\kappa)^2 |\sigma|^2C_0 (1+|\eta|^{p-1})|\xi|, 
\end{multline*}
where $\sigma$ and $D_x\sigma$ are evaluated at $(\eta,x)$. 
If we choose $\theta=\frac14$, $\kappa=\frac12$, then, in view of the fact that $|\xi-\eta|\leq \kappa\leq 1$, the above inequality implies
\begin{multline*}
K \geq C^{-1} c_0^2|\xi|^{2p}-Cc_0 -C\|\sigma\|^2_{L^\infty(\Rd)} C_0 (1+|\xi|^{p+1})  -C\|\sigma\|^2_{L^\infty(\Rd)}\|D_x\sigma\|^2_{L^\infty(\Rd)}|\xi|^2 \\
-C \|\sigma\|^2_{L^\infty(\Rd)} C_0  (1+|\xi|^{p}),
\end{multline*}
where $C$ depends only on $p$. Since $p>1$, this yields that $|\xi|\leq R_K$, as desired.
\end{example}

\begin{example}[Forced mean curvature motion] 
Consider the case that
\begin{equation*}
\sigma(\xi,x) = \sigma(\xi) = \sqrt{2}\left( I_d- \frac{\xi}{\left| \xi \right|}\otimes \frac{\xi}{\left| \xi \right|} \right) \quad \mbox{and} \quad   H(\xi,x)= a(x)\left|\xi \right|,
\end{equation*}
where $a$ is a random field satisfying $c_0 \leq a(x) \leq C_0$. Then $$A(\xi) = \frac12\sigma(\xi) \sigma(\xi)^T =  I_d- \frac{\xi}{\left| \xi \right|}\otimes \frac{\xi}{\left| \xi \right|}$$ and so equation~\eqref{e.pde} is the forced mean curvature equation~\eqref{e.MCM}. 

\smallskip

We claim that, as pointed out in~\cite{LS0}, the (LS) condition~\eqref{e.LScoercivity} is satisfied provided that 
\begin{equation}\label{e.CondLSMCM}
\inf_{x\in \R^d} \left(a(x))^2 -(d-1) \left|D_xa(x)\right|\right) > 0. 
\end{equation}
To check this, we suppose $\mathcal C_{\sigma,H}(\xi,x) \leq K$. Then there exists $\eta\in B_\kappa(\xi)$ such that
$$
K \geq \theta(1-2\theta) a(x)^2|\eta|^2 -\theta (1+\kappa)^2(d-1) \left|D_xa(x)\right||\eta||\xi| 
-\theta \kappa (1+\kappa)^2 (d-1)^2 c_0 |\xi|.
$$
Using that $|\eta-\xi|\leq \kappa\leq 1$, we get 
$$
K \geq \theta(1-2\theta) \left(  a(x)^2 -\frac{(1+\kappa)^2}{(1-2\theta)}(d-1) \left|D_xa(x) \right|\right)|\xi|^2 -C |\xi|.
$$
In view of assumption \eqref{e.CondLSMCM}, we may choose $\kappa$, $\delta$ and $\theta$ so small that 
$$
a^2 -\frac{(1+\kappa)^2}{(1-2\theta)}(d-1) \left|D_xa \right|\geq \delta>0\quad \mbox{in} \ \R^d.
$$
This implies that $|\xi|\leq R_K$ for some $R_K>0$ depending only on $(K,d,\kappa,\delta,\theta)$. 
\end{example}

\begin{example}[Anisotropic forced mean curvature motion] Our assumptions allow for the previous example to be generalized to the anisotropic setting. We may consider the case that 
$$
A(\xi,x)=\left(I_d- \frac{\xi}{\left| \xi \right|}\otimes \frac{\xi}{\left| \xi \right|} \right)\tilde A(\hat \xi,x)\left(I_d- \frac{\xi}{\left| \xi \right|}\otimes \frac{\xi}{\left| \xi \right|}\right),  \qquad H(\xi,x)= \left|B\left( \frac{\xi}{\left| \xi\right|},x\right)\xi\right|,
$$
where $\tilde A$ takes the form $\tilde A= \frac12 \tilde \sigma\tilde \sigma^T$ and $B$ is a matrix-valued random field satisfying $c_0I_d \leq B^TB \leq C_0 I_d$. Note that in this context $\sigma= \left(I_d-\hat \xi\otimes \hat \xi\right)\tilde \sigma$. 
Then the following condition implies that the (LS) condition holds:
\begin{equation*}
\inf_{(e,x) \in \partial B_1\times \Rd} \left[\left|B(e,x)e\right| -| \sigma(e,x)|^2\left(| \sigma_x(e,x)|^2+|B(e,x)e||B_x(e,x)|\right)\right] >0.
\end{equation*}
We leave the confirmation of this claim to the reader.
\end{example}

\subsection{Brief review of the literature}
Periodic and almost-periodic homogenization results were proved in~\cite{LS0} for the type of quasilinear and geometric Hamilton-Jacobi equations considered here. The fundamental qualitative homogenization results for convex Hamilton-Jacobi equations in the random setting were first proved in~\cite{RT,S} in first-order case and later in the semilinear, viscous case in~\cite{KRV,LS}. New proofs and extensions of these results appeared later in~\cite{LS3,Sch,AS1,AS2}.

\smallskip

As explained above, for geometric equations, the existence of plane-like solutions in periodic media was proved in~\cite{CD}; see also~\cite{chambolle2009homogenization} for a BV approach and~\cite{torres2004plane} for a construction of plane-like solutions with periodic exclusions. Our assumption~\eqref{e.LScoercivity} implies that the forcing term does not change sign. This restriction has been lifted, under suitable restriction, for periodic media: see~\cite{CLS}, which contains an explicit computation in the one-dimensional case; in~\cite{dirr2008pulsating} pulsating waves are constructed under a smallness condition of the forcing term; forced mean curvature motion for graphs, with sign changing velocities, is studied in~\cite{barles2011homogenization} under a rather sharp condition on the forcing term; in that setting,~\cite{cesaroni2013long} explains the construction of generalized traveling waves and analyses the long-time behavior of the motion. When the velocity changes sign, pinning phenomena may occur: this amounts to finding a stationary, positive supersolution at non-vanishing applied load~\cite{dirr2011pinning}. Another very interesting problem is to study the properties of the homogenized motion (the so-called stable norm, or equivalently, the effective Hamiltonian), examining for example its regularity: for periodic coefficients, this question is studied in~\cite{chambolle2014plane}.

\subsection{Notation}
\label{s.notation}

The symbols $C$ and $c$ denote positive constants which may vary from line to line and, unless otherwise indicated, depend only on the data and on an upper bound for $|p|$ or $\mu$. For $s,t\in \R$, we write $s\wedge t : = \min\{ s,t\}$ and $s \vee t := \max\{s,t\}$.  We denote the $d$-dimensional Euclidean space by $\Rd$,  $\N$ is the set of natural numbers and $\N^* := \N \setminus \{ 0 \}$. For each $x,y \in \Rd$, $|x|$ denotes the Euclidean length of $x$ and $x\cdot y$ the scalar product. For $r>0$, we set $B_r(x): = \{ y\in \Rd : |x-y| < r\}$ and $B_r : = B_r(0)$. If $K$ is a subset of $\R^d$, we denote by $\overline K$ and  $\partial K$ its closure and its boundary and let $K+B_r$ be the set of points which are at a distance at most $r$ of $K$. The Hausdorff  distance between two subsets $U,V\subseteq \Rd$ is  $\dist_H(U,V)=\inf\{r\geq 0 : U\subseteq V+r B_1\ {\rm and }\ V\subseteq U+r B_1\}$.  The set of bounded and Lipschitz continuous maps on $\Rd$ is denoted $W^{1,\infty}(\R^d)$. If $E$ is a set, then $\indc_E$ is the indicator function of $E$. The denote the set of upper and lower semicontinuous functions on a domain $E\subseteq\Rd$ by $\USC(E)$ and $\LSC(E)$, respectively. The space of bounded and uniformly continuous functions is denoted $\BUC(E)$.

\smallskip

We usually {\it do not} display the dependence of the various quantities with respect to the random parameter $\omega=(\sigma, H)\in\Omega$, unless this is necessary. 
 
\smallskip 
 
Throughout the paper, differential equations and inequalities are to be understood in the viscosity sense (see \cite{CIL} for the general background). 
We next recall the appropriate notion of viscosity solution for quasilinear equations like~\eqref{e.pde} which may be singular at $\xi=0$.

\begin{definition}
\label{d.viscsol}
Given $U\subseteq \Rd\times (0,\infty)$ and a function $u\in C(U)$, we say that $u$ is a \emph{subsolution} (respectively, \emph{supersolution}) of the equation
\begin{equation*} \label{}
\partial_t u - \tr\left( A \left( Du,x \right) D^2u  \right) + H(Du,x) = 0 \quad \mbox{in} \ U
\end{equation*}
if, for every $\phi\in C^2(U)$ and $(x_0,t_0) \in U$ such that
\begin{equation*} \label{}
(x,t) \mapsto u(x,t) - \phi(x,t) \quad \mbox{has a local maximum (resp., minimum) at} \ (x_0,t_0),
\end{equation*}
we have \begin{equation*} \label{}
\partial_t \phi(x_0,t_0) - \tr^*\left( A \left( D\phi(x_0,t_0),x \right) D^2\phi(x_0,t_0)  \right) + H(D\phi(x_0,t_0),x) \leq 0
\end{equation*}
(resp.,
\begin{equation*} \label{}
\partial_t \phi(x_0,t_0) - \tr_*\left( A \left( D\phi(x_0,t_0),x \right) D^2\phi(x_0,t_0)  \right) + H(D\phi(x_0,t_0),x) \geq 0 \, \mbox{)}.
\end{equation*}
Here, for any symmetric matrix $X\in \R^{d\times d}$, $\tr^*(A(\xi,x)X)$ and $\tr_*(A(\xi,x)X)$ stand for the upper and lower semicontinuous envelopes of the map $(\xi,x) \to \tr (A(\xi,x)X)$, which agree with $ \tr (A(\xi,x)X)$ on the domain of $A$ and are defined at $\xi=0$ for any matrix $X$ by 
\begin{equation*} \label{}
\tr^*(A(0,x)X) := \limsup_{\xi\to 0, \ \xi\neq 0,\  x'\to x}\tr(A(\xi,x')X), 
$$
$$
 \tr_*(A(0,x)X) := \liminf_{\xi\to 0, \ \xi\neq 0, \ x'\to x} \tr(A(\xi,x') X).
\end{equation*}
The definitions of solution for other equations encountered in this paper (such as the metric problem and approximate corrector problem) are completely analogous. 
\end{definition}

\subsection{Outline of the paper}
In the next section, we state the main result, reduce it to auxiliary results which are the focus of the rest of the paper and show that Theorem~\ref{t.homogenization} is a corollary of it.  Sections~\ref{s.fluctuations} and~\ref{s.means} are the heart of the paper: there we give the proof of homogenization for the planar metric problem. We give an estimate on the stochastic fluctuations in Section~\ref{s.fluctuations} and then Section~\ref{s.means} contains the convergence of the statistical bias. The proof of Theorem~\ref{t.homogenization} is completed in the final two sections, where we provide a deterministic link between the planar metric problem and the approximate correctors (Section~\ref{s.correctors}) and between the approximate correctors and the full time-dependent, initial-value problem (Section~\ref{s.timedependent}). Appendix~\ref{s.metricproblem} contains some auxiliary results on well-posedness and global Lipschitz bounds for the metric, approximate corrector and full time-dependent problems. 

\section{Overview of the proof of Theorem~\ref{t.homogenization}}
\label{s.outline}

In this section, we state the main result of the paper and reduce its proof to four statements which are the focus of the rest of the paper. 

\subsection{Statement of the main result}
We begin by stating the main result of the paper. Given  $L,T\geq 1$ and $\ep>0$, we consider solutions $u,u^\ep\in W^{1,\infty}(\Rd\times [0,T])$ of
\begin{equation}
\label{e.ptfmpdes}
\left\{ 
\begin{aligned}
& \partial_t u^\ep -\ep \tr\left( A\left(Du^\ep,\frac x\ep \right)D^2u^\ep \right) + H\left(Du^\ep,\frac x\ep \right) = 0 &  \mbox{in} & \ \Rd\times (0,T], \\
& \partial_t u +  \overline H(Du) = 0 & \mbox{in} & \ \Rd\times (0,T], \\
& u^\ep(\cdot,0) = u(\cdot,0) & \mbox{on} & \ \Rd,
\end{aligned}
\right.
\end{equation}
such that $u^\ep$ and $u$ satisfy the following Lipschitz estimate: for every $x,y\in\Rd$ and $t,s\in [0,T]$,
\begin{equation} \label{e.LipLs}
\left| u^\ep(x,t) - u^\ep (y,s) \right| \vee \left| u(x,t) - u(y,s) \right| \leq L \left( |x-y| + |t-s| \right).
\end{equation}
(See Remark~\ref{r.lipschitz} for comments on removing the Lipschitz hypothesis~\eqref{e.LipLs}.)

\begin{theorem}
\label{t.main}
Consider a probability measure~$\P$ on~$(\Omega,\F)$ satisfying~(P1) and~(P2). Then there exists $\overline H:\Rd\to [0,\infty)$, depending only on $\P$, such that:
\begin{itemize}
\item For every $L\geq 1$ and $\alpha \in \left( 0,\tfrac27\right)$, there exists~$C=C(\data,L,\alpha)\geq 1$ such that, for every $\xi,\eta \in \overline B_L$, 
\begin{equation}
\label{e.Hbarprops}
c_0 \left| \xi \right|^p \leq \overline H(\xi) \leq C_0 \left| \xi \right|^p
\quad \mbox{and} \quad 
\left| \overline H(\xi) - \overline H(\eta) \right| \leq C\left| \xi - \eta \right|^{\alpha}.
\end{equation}

\item For every $L,R,T\geq1$, there exists $C(\data,L,R,T)\geq 1$ and $q(\data) < \infty$ such that, for every $k\in\N$ and $\lambda\in (0,1]$, we have 
\begin{multline}
\label{e.bigestimate}
\begin{multlined}
\P \bigg[ \mbox{there exist} \ \ep \in \left[ 2^{-(k+1)}, 2^{-k} \right) \ \mbox{and} \ u,u^\ep \in W^{1,\infty}(\Rd \times [0,T])  \\
\mbox{satisfying~\eqref{e.ptfmpdes},~\eqref{e.LipLs} and} \ \sup_{(x,t) \in B_R\times [0,T]} \left| u^\ep (x,t) - u(x,t) \right| \geq \lambda \bigg]
\end{multlined} \\
 \leq C 2^{kq} \exp\left( -\frac{2^{k/5} \lambda^{18}}{C} \right).
\end{multline}
\end{itemize}
\end{theorem}

In this section we give an overview of the proof of Theorem~\ref{t.main}. Similar to the strategy in~\cite{AT1,AC}, the homogenization of time-dependent problems is reduced to the convergence of a \emph{metric problem}. Unlike in the convex case, however, we need to study the \emph{planar} metric problem, which roughly measures the distance from a point to a plane (rather than between points, as in the usual metric problem considered in the convex case~\cite{AS2}).

\subsection{Ingredients in the proof of Theorem~\ref{t.main}}
In this subsection, we state the key auxiliary propositions which are used to prove Theorem~\ref{t.main}. Recall from the introduction that the planar metric problem is
\begin{equation} \label{e.planarmetric}
\left\{ 
\begin{aligned}
& -\tr\left( A(Dm_\mu,x) D^2m_\mu\right) + H(Dm_\mu,x) = \mu & \mbox{in} & \ \H_e^+, \\
& m_\mu = 0 & \mbox{on} & \ \partial \H_e^+,
\end{aligned}
\right.
\end{equation}
where $\mu>0$ and, for a unit direction $e\in\partial B_1$, we define $\H_e^+$ and $\H_e^-$ to be the half-spaces
\begin{equation}
\label{e.He}
\H_e^+: = \left\{ x\in\Rd \,:\, x\cdot e > 0\right\} \qquad \mbox{and} \qquad  \H_e^-: = \left\{ x\in\Rd \,:\, x\cdot e < 0\right\}. 
\end{equation}
In Appendix~\ref{s.metricproblem}, we show that this problem is well-posed and there exists a unique Lipschitz solution which we denote by $m_\mu(\cdot,\H_e^-)$. See Theorem~\ref{t.LipschitzAppenA1}. The main step in the proof of Theorem~\ref{t.homogenization} is the demonstration that~\eqref{e.planarmetric} homogenizes, that is, that for every $\mu>0$ and $e\in \partial B_1$, there exists a deterministic quantity $\overline m_\mu(e)>0$ such that 
\begin{equation} \label{e.MPlim}
\lim_{t\to \infty} \frac {m_\mu(te,\H_e^-)}{t} = \overline m_\mu(e), \quad \mbox{$\P$--a.s.}
\end{equation}
We actually prove more, namely the quantitative version of~\eqref{e.MPlim}  given in the following proposition.

\begin{proposition}
\label{p.planarquantitative}
Fix $L\geq 1$, $e \in \partial B_1$ and $\mu\in (0,L]$, there exist $\overline m_\mu(e)\geq 0$, $C(\data,L)\geq1$ and $q(\data)<\infty$ such that, for every $x\in \H_e^+$ and $\lambda>0$,
we have the estimate
\begin{equation}\label{e.Pm-barm}
\P \left[ \left| m_\mu(x,\H_e^-) - \overline m_\mu (e)(x\cdot e) \right| > \lambda \right]  \leq C \exp\left( \frac{-\lambda^2}{(1+x\cdot e)^{\frac85}} \right).
\end{equation}
Moreover, the map $(\mu,e) \mapsto \overline m_\mu(e)$ is continuous on $(0,\infty) \times \partial B_1$.
\end{proposition}

Proposition~\ref{p.planarquantitative} is a consequence of the results in Sections~\ref{s.fluctuations} and~\ref{s.means} and its proof comes near the end of Section~\ref{s.means}.

\smallskip

Using the result of Proposition~\ref{p.planarquantitative}, we can identify~$\overline H$.

\begin{definition}[{The effective Hamiltonian}]
We define the effective Hamiltonian $\overline H:\Rd \to \R$ by setting $\overline H(0) := 0$ and, for every $t>0$ and $e\in\partial B_1$, 
\begin{equation*}
\overline H(te):= \inf\left\{ \mu>0 \,:\, \overline m_\mu(e) > t \right\}.
\end{equation*}
\end{definition}

It is immediate that $\overline H:\Rd \to \R$ is continuous and coercive and, since $\overline m_\mu(e)$ is increasing in $\mu$, its sublevel sets are star-shaped with respect to the origin. We do not expect~$\overline H$ to be positively homogeneous, in general, unless $p=1$ or $\P \left[ A\equiv 0 \right]=1$, due to the interaction between the diffusion and the Hamiltonian. 

\smallskip

The quantitative error estimates for the time-dependent initial-value problem require an explicit H\"older estimate for $\overline H$, which is proved at the end of Section~\ref{s.correctors}.

\begin{proposition}
\label{p.Hbarprops}
For every $L\geq 1$ and $\alpha \in \left( 0,\tfrac 27 \right)$, there exists a constant $C=C(\data,L,\alpha)\geq 1$ such that, for every $\xi,\eta \in \overline B_L$, 
\begin{equation}
\label{e.Hbarpgrowth1}
c_0 \left| \xi \right|^p \leq \overline H(\xi) \leq C_0 \left| \xi \right|^p
\end{equation}
and
\begin{equation}
\label{e.Hbarholder1}
\left| \overline H(\xi) - \overline H(\eta) \right| \leq C\left| \xi - \eta \right|^{\alpha}.
\end{equation}
\end{proposition}

Once we have proved Proposition~\ref{p.planarquantitative}, the rest of the proof of Theorem~\ref{t.homogenization} is deterministic and consists of transfering the limit for the planar metric problem to a limit for solutions of the time-dependent problem by comparison arguments. This part of the proof of our main result is a fairly routine, if technical, adaptation of the perturbed test function method~\cite{E2}. However, it is more difficult in our setting than usual, due to the presence of the singular, quasilinear diffusion (as has been noticed previously~\cite{CM}) and due to the need for a quantitative statement. 

\smallskip

We state the deterministic comparison results in two parts. First, we give the link between the planar metric problem and the \emph{approximate correctors}, which are the solutions $v^\delta(\cdot,\xi)$, for each $\xi\in\Rd$ and $\delta > 0$, of the problem 
\begin{equation*} \label{}
\delta v^\delta (x,\xi) - \tr\left( A\left(\xi+Dv^\delta(x,\xi),x\right)D^2v^\delta (x,\xi) \right) + H\left(\xi+Dv^\delta(x,\xi),x\right) = 0 \quad \mbox{in} \ \Rd.
\end{equation*}
It turns out that~\eqref{e.MPlim} is essentially equivalent to the limit
\begin{equation} \label{e.dvdlim}
\lim_{\delta\to 0} - \delta v^\delta (0,\xi) = \overline H(\xi), \quad \mbox{$\P$--a.s.}
\end{equation}
A quantitative version of this fact is summarized in the following proposition, which is proved in Section~\ref{s.correctors}.

\begin{proposition}
\label{p.correctors}
Fix $\lambda,\delta \in (0,1]$, $L\geq 1$, $e\in \partial B_1$ and $\mu\in (0,L]$. Select~$t>0$ such that $\xi = te$ satisfies $\mu = \overline H(\xi)$. Then there exists $C(\data,L) \geq 1$ such that, for every $s\geq C/\lambda\delta$, 
\begin{equation*} \label{}
\sup_{x\in B_{s/2}} \left| m_\mu(x,\H_e^--se) - \overline m_\mu(e)\left( s+(x\cdot e)\right) \right| \leq \frac\lambda \delta 
\implies \left| \delta v^\delta(0,\xi) + \overline H(\xi) \right| \leq C \lambda^{\frac 15}.
\end{equation*}
\end{proposition}

\smallskip

It is essentially well-known that, in a fairly general framework, the limit~\eqref{e.dvdlim} implies that~\eqref{e.pde} homogenizes. The final ingredient for the proof of Theorem~\ref{t.homogenization} is a quantitative version of this statement. It is proved in Section~\ref{s.timedependent}. 

\begin{proposition}
\label{p.PTFM}
Fix~$0< \ep \leq \delta \leq \lambda \leq 1$ and $L,R,T\in[1,\infty)$. Suppose that $u^\ep,u \in W^{1,\infty}(\Rd \times [0,T])$ satisfy~\eqref{e.ptfmpdes} and~\eqref{e.LipLs}. Then
\begin{equation} \label{e.errorassumps}
\sup_{(x,t) \in B_R \times [0,T] } \left| u^\ep (x,t) - u(x,t) \right| \geq \lambda 
\end{equation}
implies that, for a constant $C=C(\data,L,R,T)\geq 1$, 
\begin{equation} \label{e.bigtdassums}
\sup_{(x,\xi) \in B_{C/\ep } \times \overline B_L} \left| \delta v^\delta (x,\xi) + \overline H(\xi) \right| \geq \frac\lambda {C} - C\left( \frac{\ep}{\delta\lambda} \right)^{\frac1{10}}.
\end{equation}
\end{proposition}

\subsection{Proof of Theorems~\ref{t.homogenization} and~\ref{t.main}}

We next give the proof of Theorem~\ref{t.main}, subject to the four results stated in the previous subsection (and some more standard auxiliary estimates proved later in the paper). 

\begin{proof}[{Proof of Theorem~\ref{t.main}}]
We prove only the second statement of the theorem, since the first is contained in Proposition~\ref{p.Hbarprops}. Fix $L,R,T\geq 1$. We denote by $C$ and $c$ positive constants with depend only on~$(\data,L,R,T)$ and may vary in each occurrence.

\smallskip

Also fix $k\in\N$, $2^{-k} \leq \delta \leq \lambda \leq 1$. Below we will select~$\delta$ in terms of $k$ and $\lambda$. The goal is to estimate the probability of the event that, for some~$\ep \in \left( 2^{-(k+1)}, 2^{-k} \right]$ and $u,u^\ep\in W^{1,\infty}(\Rd \times [0,T])$ satisfying~\eqref{e.ptfmpdes} and~\eqref{e.LipLs}, we have 
\begin{equation}
\label{e.badevent}
\sup_{(x,t) \in B_R\times [0,T]} \left| u^\ep (x,t) - u(x,t) \right| \geq \lambda. 
\end{equation}
According to Proposition~\ref{p.PTFM}, we see that~\eqref{e.badevent} implies 
\begin{equation}
\label{e.badnessdvd}
\sup_{(x,\xi) \in B_{C2^k} \times \overline B_L} \left| \delta v^\delta (x,\xi) + \overline H(\xi)  \right| \geq c \lambda,
\end{equation}
provided that
\begin{equation*}  
\left( \frac{C\ep}{\delta\lambda} \right)^{\frac1{10}} \leq c\lambda,
\end{equation*}
which is equivalent to 
\begin{equation} \label{e.qualifier}
  \lambda \geq C \left( \frac{2^{-k}}{\delta} \right)^{\frac1{11}}.
\end{equation}
We next apply Proposition~\ref{p.correctors}. We deduce that~\eqref{e.badnessdvd} implies, for $s:= C\delta^{-1} \lambda^{-5}$,
\begin{multline}
\label{e.badnessmmu}
\sup \Big\{  \big| m_\mu(y+x,\H_e^-+x-se) - \overline m_\mu(e)\left(s+ y\cdot e\right) \big| \\ 
:\, ( \mu, e, x, y) \in (0,CL^p]\times \partial B_1\times B_{C2^k}\times  B_{s/2} 
\Big\} 
\geq c \lambda^{5} \delta^{-1}.
\end{multline}
In order to estimate the probability of~\eqref{e.badnessmmu}, we need to snap to a finite grid so that we can apply union bounds to Proposition~\ref{p.planarquantitative}. This requires some (deterministic) continuity of $\overline m_\mu(e)$ and $m_\mu(x,\H^-_e)$ in all three parameters $x$, $e$ and $\mu$, which is given in~Lemmas~\ref{l.mbarcont} and~\ref{l.mconte}. We deduce from these lemmas and~\eqref{e.badnessmmu} the existence of a finite set 
\begin{equation*}
\Lambda \subseteq \left( 0, CL^p \right] \times \partial B_1 \times B_{C2^k} \times B_{s/2},
\end{equation*}
depending only on $(\data,k,\delta,\lambda)$, such that $\Lambda$ has at most $C 2^{kq}\left(\lambda \delta \right)^{-q}$ elements, for an exponent $q=q(\data)<\infty$, and\begin{equation}
\label{e.snaptogrid}
\sup_{\left( \mu, e,x,y \right) \in \Lambda} \left| m_\mu(y+x ,\H_e^-+x-se) - \overline m_\mu(e) (s+y\cdot e) \right| \geq c \lambda^{5} \delta^{-1} - C.
\end{equation}
Making $q(\data)<\infty$ larger, if necessary, and applying Proposition~\ref{p.planarquantitative}, we deduce that, for every $\left( \mu,e,x,y\right) \in \Lambda$,
\begin{multline*}
\P\left[\left| m_\mu(y+x,\H_e^-+x-se) - \overline m_\mu(e) (s+y\cdot e) \right| \geq c \lambda^{5} \delta^{-1} -C \right] 
\\
\leq Cs^q\exp\left( \frac{-c\lambda^{10}}{\delta^2s^{\frac85}} \right) 
 \leq C\lambda^{-q}\delta^{-q}\exp\left( \frac{-c\lambda^{18}}{\delta^{\frac25}} \right),
\end{multline*}
provided that $\lambda^5\delta^{-1} \geq C$ (so that $c \lambda^{5} \delta^{-1} - C \geq c \lambda^{5} \delta^{-1}$) and $s\geq C/(\lambda\delta)$, which in view of the definition of~$s$ is equivalent to $\lambda \leq c$. Moreover, the first restriction that~$\lambda^5\delta^{-1} \geq C$ can be removed since otherwise the last term on the right side is larger than~$1$. Therefore, a union bound gives, up to a redefinition of $q(\data)<\infty$, 
\begin{multline}
\label{e.breakingbad}
\P \big[  \mbox{$\exists u^\ep,u\in W^{1,\infty}(\Rd \times [0,T])$ satisfying~\eqref{e.ptfmpdes},~\eqref{e.LipLs} and \eqref{e.badevent}} \big] \\
 \leq  \P \left[  \mbox{\eqref{e.snaptogrid} holds} \right]   \leq C 2^{kq}\left(\lambda \delta \right)^{-q} \exp\left( \frac{-c\lambda^{18}}{\delta^{\frac25}} \right).
\end{multline}
It remains to specify $\delta$. Making no attempt to be optimal, we take $\delta:= 2^{-k/2}$.
It is clear that this choice satisfies~\eqref{e.qualifier} provided that 
\begin{equation*} \label{}
\lambda \geq C2^{-k/22}.
\end{equation*}
After another redefinition of~$q=q(\data)<\infty$, and we see that the right side of~\eqref{e.breakingbad} is at most  
\begin{equation*}
C 2^{kq} \exp\left( -c2^{k/5} \lambda^{18}  \right). 
\end{equation*}
This completes the proof of the theorem, since the restriction on $\lambda$ may be removed, since if it is false, the quantity in the previous line is larger than~$1$.
\end{proof}

We now show that Theorem~\ref{t.homogenization} is a consequence of Theorem~\ref{t.main}. 

\begin{proof}[Proof of Theorem~\ref{t.homogenization}]
Take $\alpha<\frac1{90}$, set $\lambda(k):= 2^{-\alpha k}$ and observe that, for constants $R,T\geq 1$, the estimate ~\eqref{e.bigestimate} yields, for every $m\in\N$ sufficiently large,
\begin{multline*}
\P \left[ \exists \ep \in (0,2^{-m}], \, \sup_{(x,t) \in B_R\times [0,T]}  \ep^{-\alpha} \left| u^\ep(x,t)  - u(x,t) \right| \geq 1 \right]\\ 
\leq C \sum_{k=m}^\infty 2^{kq} \exp\left( -2^{k/5- 18k\alpha}/C \right) \leq C \exp\left( - 2^{\left(\frac15 - 18\alpha \right) m} \right),
\end{multline*}
where $C$ depends on $\data$ and the Lipschitz constant of $u^\ep$ and $u$ on $\Rd \times[0,T]$.
Summing over $m$, applying the Borel-Cantelli lemma and then shrinking $\alpha$ slightly yields 
\begin{equation*}
\P \left[\limsup_{\ep \to 0} \, \sup_{(x,t) \in B_R\times [0,T]}  \ep^{-\alpha} \left| u^\ep(x,t)  - u(x,t) \right| =0 \right]  = 1.
\end{equation*}
Taking the intersection of these events for a sequence $R = R_j \to \infty$ then yields the conclusion of the theorem.
\end{proof}

\begin{remark}
\label{r.lipschitz}
Theorems~\ref{t.homogenization} and~\ref{t.main} are stated and proved under the condition that both~$u^\ep$ and~$u$ are Lipschitz continuous in space and time. For~$u$, this is not controversial, since it satisfies a first-order equation with a coercive Hamiltonian and will thus be locally Lipschitz on~$\Rd \times (0,T]$, at least if it is assumed to have at most affine growth initially. However, $u^\ep$ will not be Lipschitz, in general, unless it is bounded in~$C^{1,1}(\Rd)$ at~$t=0$ (that initial data belonging to $C^{1,1}(\Rd)$ suffices for a Lipschitz estimate is explained in Appendix~\ref{s.metricproblem}). Therefore it may appear that Theorem~\ref{t.homogenization} only implies even qualitative homogenization for initial-value problems with sufficiently regular initial data.  

\smallskip

However, this can be overcome easily, using the comparison principle and approximating any  bounded and uniformly continuous initial condition from above and below by $C^{1,1}$ functions. The monotonicity of the solutions as functions of the initial data guarantees that we may interchange the two limits (for the approximation and homogenization), yielding a quite general qualitative homogenization result.  

\smallskip

This interpolation trick also works at the level of the quantitative estimates, provided we assume H\"older continuous initial data, permitting us to deduce error estimates and an algebraic rate of convergence for more general initial-value problems. In order to check this, it is necessary to track the dependence on~$L$ (the upper bound for parameters such as~$\mu$,~$|\xi|$, etc) of all the constants~$C$ in each quantitative estimate in  the paper, in order to ensure that the dependence is polynomial in~$L$ (which it is). For the readability of the paper, we have chosen not to display such dependence. 
\end{remark}

\section{The fluctuations estimate}
\label{s.fluctuations}

In this section, we prove the following estimate on the fluctuations of the metric problem to any nonempty compact target set $S\subseteq \Rd$ which satisfies the following interior ball condition: 
\begin{equation}\label{e.IntBallPpt}
S = \bigcup_{\overline B_1(x) \subseteq S} \overline B_1(x).
\end{equation}
The metric problem is
\begin{equation} 
\label{e.MPforS}
\left\{ 
\begin{aligned}
&  -\tr\left(A(Dm,x)D^2m\right) +H(Dm, x) = \mu & \mbox{in} & \ \R^d\setminus S, \\
& m = 0 & \mbox{on} & \ \partial S.
\end{aligned}
\right.
\end{equation}
In Appendix A we show that~\eqref{e.MPforS} is well-posed and give some properties of its solution, which we denote by $m_\mu(\cdot,S) \in W^{1,\infty}_{\mathrm{loc}}\big(\overline{\Rd\setminus S}\big)$.

\smallskip

The main result of this section is the following estimate for the stochastic fluctuations of~$m_\mu(x,S)$ for a point $x\in\Rd\setminus S$.

\begin{proposition}
\label{p.fluctuations}
Let $L\geq 1$ and $S\subseteq\Rd$ be a compact set satisfying~\eqref{e.IntBallPpt}. Then there exists $C(\data,L)\geq1$ such that, for every $\mu\in (0,L]$, $x\in\Rd$ and $\lambda>0$, 
\begin{equation}
\label{e.fluctuations}
\P \big[ \left| m_\mu(x,S) - \E\left[ m_\mu(x,S) \right] \right| > \lambda \big] \leq C\exp\left( -\,\frac{\mu^3 \lambda^2}{C\left( 1+\dist(x,S)\right)}\right).
\end{equation}
\end{proposition}

Since the constant~$C$ in Proposition~\ref{p.fluctuations} does not depend on $S$, we obtain the same result for the planar metric problem by considering an increasing sequence~$\left\{ S_n \right\}_{n\geq 1}$ of compact sets whose union is~$\H_e^-$ and using the stability of viscosity solutions under local uniform convergence and the obvious monotonicity of~$m_\mu(x,S_n)$ in~$n$. 

\begin{corollary}
\label{c.fluctuationsH} 
Let $L\geq 1$. Then there exists $C(\data,L)\geq 1$ such that, for every $\mu\in (0,L]$, $e\in\partial B_1$, $x\in \H_e^+$ and $\lambda>0$,
\begin{equation}
\label{e.fluctuationsP}
\P \left[ \left| m_\mu(x,\H_e^-) - \E\left[ m_\mu(x,\H_e^-) \right] \right| > \lambda \right] \leq C\exp\left( -\,\frac{\mu^3 \lambda^2}{C\left( 1+e\cdot x\right)}\right).
\end{equation}
\end{corollary}

The proof of Proposition~\ref{p.fluctuations} is based on Azuma's inequality and is similar to the arguments used by the authors in the viscous convex case~\cite{AC}, which were partially based on those introduced in the first-order convex case~\cite{ACS,MN} and some previous ideas originating in first-passage percolation~\cite{K,Z}. 

\smallskip

We continue with some notation and conventions in force throughout the remainder of this section. We fix $L\geq 1$, $\mu\in (0,L]$ and a compact set $S\subseteq \Rd$ satisfying~\eqref{e.IntBallPpt}. Unless otherwise stated, we denote by $C$ and $c$ positive constants  which may vary in each occurrence and depend only on~$(\data,L)$. Some of our estimates below depend on a \emph{lower} bound for $\mu$, these are typically denoted by $C_\mu$ or $c_\mu$ with the dependence on $\mu$ made explicit. The constants $\ell_\mu$ and $L_\mu$ are as in the statement of Lemma~\ref{l.EstiMetric} in Appendix~\ref{s.metricproblem}, and we have that, for some constants $C(\data,L)\geq 1$ and $c(\data)>0$,
\begin{equation}
\label{e.Lmubounds}
c\mu \leq \ell_\mu \leq L_\mu \leq C.
\end{equation}
For technical reasons, it is convenient to consider solutions of the metric problem for coefficients $(\sigma, H)$ belonging to the closure $\overline \Omega$ of the set $\Omega$ with respect to the topology of local uniform convergence. The well-posedness and global Lipschitz estimates for coefficients belonging to $\overline \Omega$ is presented in Appendix~\ref{s.metricproblem}. If it necessary to display the dependence of $m_\mu(\cdot,S)$ on the coefficients, we write $m_\mu(\cdot,S,\omega)$ for $\omega=(\sigma,H) \in \overline \Omega$. 

\subsection{Localization in sublevel sets}

A key step in the proof of Proposition~\ref{p.fluctuations}, following~\cite{AC}, is to show that the solutions of the planar metric problem depend almost entirely on the coefficients restricted to their sublevel sets. This paves the way for a martingale argument to estimate the stochastic fluctuations. The result is summarized in Proposition~\ref{p.localizationPC}, below. The main new observation here is that the proof of~\cite[Lemma 3.3]{AC} does not require convexity of the Hamiltonian, rather a weak form of positive homogeneity. Nevertheless, we give a complete argument here for the reader's convenience and because the statement presented here is slightly different than the one in~\cite{AC}.

\begin{lemma} 
\label{l.seige}
Fix coefficients $\omega_1=(\sigma_1,H_1)\in\overline\Omega$ and $\omega_2=(\sigma_2,H_2)\in\overline \Omega$. Suppose that $t\geq 1$ and
\begin{equation*}
(\sigma_1,H_1)\equiv (\sigma_2,H_2)  \quad \mbox{in} \quad   \Rd \times \left\{ x\in\Rd\setminus S\,:\, m_\mu(x,S, \omega_1)\leq t\right\}.
\end{equation*}
Then, for every $x\in\Rd\setminus S$ such that $m_\mu(x,S,\omega_1) \leq t$, we have 
\begin{multline}\label{e.seige}
m_\mu(x,S,\omega_1 ) - m_\mu(x,S,\omega_2 ) \\
\leq \frac{4 C_0  L_\mu^3}{\mu l_\mu} \exp\left( \frac{4L_\mu}{l_\mu} \right) \exp\left( - \frac{\mu }{ C_0  L_\mu^2} \big( t - m_\mu\left(x,S,\omega_2 \right) \big)\right).
\end{multline}
\end{lemma}
\begin{proof}
For notational simplicity, we denote $m_i:= m_\mu(\cdot,S,\omega_i)$ for $i\in\{1,2\}$. The argument is a comparison between $m_1$ and $w:= \varphi(m_2)$, where $\varphi:\R_+ \to \R_+$ is given by
\begin{equation*}\label{}
\varphi(s) := s + k  \exp\left( \alpha (s-t+k) \right),
\end{equation*}
and the constants $k$ and $\alpha$ are defined by
\begin{equation*}
\left\{
\begin{aligned}
& k := \sup \left\{ m_1(x)- m_2(x )  \,: \, x\in \Rd\setminus S, \, m_1(x ) = t \right\}, \\
& \alpha := ( C_0  L_\mu^2)^{-1}\mu.
\end{aligned}
\right.
\end{equation*}
We may assume without loss of generality that $k>0$, since otherwise \eqref{e.seige} is immediate and there is nothing more to show. 

\smallskip

As in the proof of~\cite[Lemma 3.3]{AC}, we show by a direct computation that $w$ is a supersolution of the equation with coefficients $\omega_1$. We perform the computation as though $m_2$ is smooth; however what follows can be made rigorous in the viscosity sense in the usual manner, by performing the analogous computation on a smooth test function and using Proposition~\ref{p.compaAppenA2}.

\smallskip

We compute
\begin{align*}
Dw(x) & = \left( 1 + \alpha k \exp\left( \alpha \left( m_2(x) - t +k \right) \right)  \right) Dm_2(x),  \\
D^2w(x) & = \left( 1 + \alpha k \exp\left(\alpha \left( m_2(x) - t +k \right) \right)  \right)   D^2m_2(x) \\ & \qquad +\alpha^2 k\exp\left( \alpha \left( m_2(x) - t +k \right) \right) Dm_2(x) \otimes Dm_2(x).
\end{align*}
Using the homogeneity of $A_1$ and $H_1$ with respect to the gradient variable, the bound $|A_1|\leq  C_0 $ and the gradient estimate $|Dm_2|\leq L_\mu$, we find:
\begin{multline*}
-\tr\left( A_1(Dw, x) D^2w \right) + H_1(Dw,x )\\  \geq (1+\alpha k) \left( -\tr\left( A_1(Dm_2, x ) D^2m_2 \right) + (1+\alpha k)^{p-1} H_1(Dm_2,x )\right)
-  C_0  \alpha^2 k L_\mu^2.
\end{multline*}
As  $(\sigma_1,H_1)=(\sigma_2,H_2)$ in $\{m_1\leq t\}$ and recalling the definition of $\alpha$, we have therefore, for every $x\in  \{m_1\leq t\}$, 
$$
-\tr\left( A_1(Dw(x), x ) D^2w(x) \right) + H_1(Dw(x),x )  \geq (1+\alpha k) \mu -  C_0  \alpha^2 k L_\mu^2 = \mu.
$$
So $w$ satisfies
\begin{equation}\label{e.pert1}
-\tr\left( A_1(Dw,x ) D^2w \right) + H_1(Dw,x ) \geq \mu \quad \mbox{in} \  \left\{ m_1 \leq t \right\}
\end{equation}
and, by the definition of $k$:
\begin{equation}\label{e.pert2}
w\geq 0\quad \mbox{on} \ \overline  S \quad \mbox{and} \quad w \geq m_1 \quad \mbox{on} \ \partial \left\{ m_1  \leq t \right\}.
\end{equation}
By  comparison principle (Proposition \ref{p.compaAppenA2}), we obtain that
\begin{equation} \label{e.cmppert}
w \geq m_1 \quad \mbox{in} \ \left\{ m_1  \leq t\right\}.
\end{equation}
Rewriting the previous inequality in terms of~$m_2$ yields
\begin{equation}\label{e.gumpt}
m_1( x ) \leq  m_2(  x ) + k \exp\left( \frac{\mu}{ C_0  L_\mu^2}(m_2(x )-t+k) \right), \quad x\in\left\{ m_1 \leq t \right\}.
\end{equation}

\smallskip

The rest of the argument follows the proof of Lemma 3.3 in \cite{AC} and consists in estimating the constant $k$ by using Lipschitz estimates. 
For this we first note that, since $k= t-\min\left\{m_2(x )\,: \, m_1(x ) =  t\right\}$, there exists $x_0\in \R^d$ such that $m_1(x_0 ) = t$ and $m_2(x_0 ) = t -k$.  Then the Lipschitz estimate on $m_1$ in \eqref{e.Lipschpc} implies that  $\dist(x_0,S)\geq L_\mu^{-1} t$ while the lower bound on the growth of $m_2$ in~\eqref{control2} yields that $(l_\mu/L_\mu) t-2 \leq m_2(x_0)= t-k$. We get a first, rough bound on $k$: 
\begin{equation}\label{e.roughbdk}
k\leq t(1-l_\mu/L_\mu)+2.
\end{equation}
Next we fix $h\in [0,t-k]$ and note that, by the growth of $m_2$ in \eqref{e.movefron}, there exists $x_h\in \Rd$ such that $m_2(x_h )=t-k-h$ and $|x_h-x_0|\leq h/l_\mu+2$. Note that, by definition of $k$, the set $\{m_2(\cdot)\leq t-k\}$ is entirely contained in $\{m_1(\cdot)\leq t\}$, so that $x_h$ also belongs to $\{m_1(\cdot)\leq t\}$. Using \eqref{e.gumpt} and~\eqref{e.Lipschpc}, we find that
\begin{align}
\label{kjhbefvz}
t-L_\mu (l_\mu^{-1}h+2)
\leq m_1(x_h ) 
& \leq m_2(x_h )+ k \exp\left( \frac{\mu}{ C_0  L_\mu^2}(m_2( x_h )-t+k) \right) \\
& \leq t-k-h+ k \exp\left(- \frac{\mu}{ C_0  L_\mu^2}h \right). \notag
\end{align}
Fix $\ep:=\exp(-1)$ and set $h:=  \mu^{-1}( C_0  L_\mu^2)$. Observe that, in view of~\eqref{e.roughbdk}, we have $t-k-h\geq 0$ provided that $t\geq (2 C_0  L_\mu^3)/(\mu l_\mu)$. Then \eqref{kjhbefvz} gives
$$
k \leq \frac{1}{1-\ep} \left( -\frac{ C_0  L_\mu^2}{\mu}+ \frac{ C_0  L_\mu^3}{\mu l_\mu}+2L_\mu\right)\leq 4 \frac{ C_0  L_\mu^3}{\mu l_\mu}.
$$
Inserting this into~\eqref{e.gumpt} yields~\eqref{e.seige} for $t\geq (2 C_0  L_\mu^3)/(\mu l_\mu)$. We conclude by noting that~\eqref{e.seige}  always holds for $t\leq 2 C_0  L_\mu^3/ (\mu l_\mu)$.
\end{proof}

It is convenient to rewrite the statement of Lemma~\ref{l.seige} in terms of the sublevel sets of $m_\mu(\cdot,S)$. To this end, we set
\begin{equation*}
\bar a_\mu:= 3\bar a'_\mu+2l_\mu, \quad \mbox{where}\quad  \bar a_\mu':= \frac{ C_0  L_\mu^3}{\mu l_\mu}-\frac{ C_0  L_\mu^2}{\mu }\log\left(\frac{\mu l_\mu^3}{4 C_0  L_\mu^3}\right).
\end{equation*}
Notice that~$\bar a_\mu = C(\mu l_\mu)^{-1}$.

\begin{proposition}
\label{p.localizationPC}
Fix coefficients $\omega_1=(\sigma_1,H_1)\in\overline\Omega$ and $\omega_2=(\sigma_2,H_2)\in\overline \Omega$. Suppose that $t\geq \bar a_\mu$ and 
\begin{equation*}
(\sigma_1,H_1) \equiv (\sigma_2,H_2)  \quad \mbox{in} \ \Rd \times \left\{  m_\mu(\cdot,S, \omega_1)\leq t\right\}.
\end{equation*}
Then 
\begin{equation}\label{e.seige2}
\left|m_\mu(x,S,\omega_1)-m_\mu(x,S,\omega_2)\right|\leq l_\mu, \quad x\in  \left\{ m_\mu(\cdot,S, \omega_1)\leq t - \overline a_\mu \right\}
\end{equation}
and, for any $s\in [0,t- \overline a_\mu]$, 
\begin{equation}\label{e.seige3}
\dist_H\big( \left\{m_\mu(\cdot,S,\omega_1)\leq s\right\},\, \left\{m_\mu(\cdot,S,\omega_2)\leq s\right\}\big) \leq 3. 
\end{equation}
\end{proposition}

\begin{proof} 
We use the notation $m_i:= m_\mu(\cdot,S,\omega_i)$ for $i\in\{1,2\}$ as in the proof of the previous lemma. 
By the definition of $\bar a_\mu'$, above, we have
\begin{equation*}
l_\mu = \frac{4 C_0  L_\mu^3}{\mu l_\mu} \exp\left( \frac{4L_\mu}{l_\mu} \right) \exp\left( - \frac{\mu }{ C_0  L_\mu^2} \bar a_\mu'\right).
\end{equation*}
Fix $x\in\Rd\setminus S$ such that $m_1(x)\leq t-\bar a_\mu'$. If $m_1(x)>m_2(x)$, then Lemma \ref{l.seige} gives  
\begin{align*}
m_1(x ) - m_2(x )  & \leq \frac{4 C_0  L_\mu^3}{\mu l_\mu} \exp\left( \frac{4L_\mu}{l_\mu} \right) \exp\left( - \frac{\mu }{ C_0  L_\mu^2} ( t - m_2(x ) )\right)  \\
& \leq \frac{4 C_0  L_\mu^3}{\mu l_\mu} \exp\left( \frac{4L_\mu}{l_\mu} \right) \exp\left( - \frac{\mu }{ C_0  L_\mu^2} ( t - m_1(x ) )\right) \\
& \leq \frac{4 C_0  L_\mu^3}{\mu l_\mu} \exp\left( \frac{4L_\mu}{l_\mu} \right) \exp\left( - \frac{\mu }{ C_0  L_\mu^2} \bar a_\mu' )\right) \leq l_\mu.
\end{align*}
%
%
Let $s$ be the largest real number such that $\{m_2(\cdot)\leq s\}\subset \{m_1(\cdot)\leq t-\bar a_\mu'\}$. Then there exists $x\in\Rd\setminus S$ such that $m_2(x)=s$ and $m_1(x)=t-\bar a_\mu'$, so that 
$$
t-\bar a_\mu'= m_1(x)\leq m_2(x)+l_\mu= s+l_\mu.
$$
Therefore $\{m_2(\cdot)\leq t-\bar a_\mu'-l_\mu\}\subset \{m_1(\cdot)\leq t\}$ and hence 
\begin{equation*}
(A_1,H_1)\equiv (A_2,H_2) \quad \mbox{in} \ \left\{m_2(\cdot)\leq t-\bar a_\mu'-l_\mu\right\}.
\end{equation*}
Reversing the roles of $m_1$ and $m_2$ in the above argument, we obtain
\begin{equation*}
m_2\leq m_1+l_\mu \quad \mbox{in} \ \left\{m_2(\cdot)\leq t-2\bar a_\mu'-l_\mu\right\}.
\end{equation*}
Arguing as above, one also has that $\{m_1(\cdot)\leq t-3\bar a_\mu'-2l_\mu\}\subset \{m_2(\cdot)\leq t-2\bar a_\mu'-l_\mu\}$, which shows that
\begin{equation*}
\left|m_2-m_1\right| \leq l_\mu \quad \mbox{in} \ \left\{m_1(\cdot)\leq t-3\bar a_\mu'-2l_\mu\right\}.
\end{equation*}
Recalling that $\bar a_\mu=3\bar a_\mu'+2l_\mu$, we get~\eqref{e.seige2}. 

\smallskip

We now prove \eqref{e.seige3}. If $s\leq t-\bar a_\mu$, then  using \eqref{e.seige2} and~\eqref{e.movefron},
$$
\left\{m_1(\cdot)\leq s\right\}\subset \left\{m_2(\cdot)\leq s+ l_\mu\right\}\subseteq \left\{m_2(\cdot)\leq s\right\}+ \overline B_{1+2l_\mu} \subseteq \left\{m_2(\cdot)\leq s\right\}+ \overline B_{3}.
$$
Hence $\{m_1(\cdot)\leq s\}\subset \{m_2(\cdot)\leq s\}+ \overline B_{3}$. Reversing the roles of $m_1$ and $m_2$ in the above argument completes the proof of~\eqref{e.seige3}.   
\end{proof}

\subsection{Construction of the localized approximations \texorpdfstring{$m_\mu^U$}{}}

Given a compact set $U\subset \R^d$ such that $S\subset U$, we now define the localized approximations to $m_\mu(\cdot,S)$ which we denote by $m_\mu^U(\cdot,S)$. For this purpose we fix a family $\left\{\omega_n=(\sigma_n,H_n)\right\}$ which is dense in $\Omega$ with respect to the topology of local uniform convergence. 
We define, for every $k,n\in\N$, 
\begin{equation*}
B_{U,k}(\omega_n)=\left\{ \omega=(\sigma,H)\in \Omega \,:\,  \sup_{(\xi,x)\in B_k\times U} \left|(\sigma,H)(\xi,x)-(\sigma_n,H_n)(\xi,x)\right| \leq \frac1k\right\}
\end{equation*}
and, for each $x\in\Rd$ and $\omega\in\Omega$,
\begin{equation*}
m_\mu^U(x,S,\omega):= \inf_{k\in\N} \sup_{n\in \N}  m_\mu(x,S,\omega_n)\indc_{B_{U,k}(\omega_n)}(\omega).
\end{equation*}
Note that the events $B_{U,k}(\omega_n)$ are nonincreasing in $k$, and therefore the infimum in the definition of $m_\mu^U$ is also a limit as $k\to \infty$. As usual, we suppress the dependence of $m_\mu$ on~$\omega$ if there is no loss of clarity. 

\smallskip

We next verify some basic properties of $m_\mu^U$.

\begin{lemma} 
\label{l.basicpptmU} 
For each compact subset $U\subseteq\Rd$ satisfying $S\subseteq U \subseteq \Rd$ and every $x\in \Rd\setminus S$, the random variable 
$m_\mu^U(x,S)$ is $\mathcal F(U)$--measurable. For every $x,y\in\Rd\setminus S$,
\begin{equation*}
\left| m_\mu^U(x,S) - m_\mu^U(y,S)\right|  \leq L_\mu\left| x-y \right| 
\end{equation*}
and 
\begin{equation*}
m_\mu(\cdot,S)\leq m_\mu^U(\cdot,S) \quad \mbox{in} \ \R^d \setminus S.
\end{equation*}
Finally, for every $\omega=(\sigma,H)\in \Omega$ and $x\in \R^d$, there exists $\omega'=(\sigma',H')\in \overline \Omega$ such that $(\sigma,H)\equiv(\sigma',H')$ in $\R^d\times U$ and $m_\mu^U(x,S,\omega) = m_\mu(x,S,\omega')$.
\end{lemma}

\begin{proof} 
As the coefficient fields $(\sigma,H)\in \Omega$ are locally uniformly continuous, the event $B_{U,k}(\omega_n)$ belongs to $\F(U)$. So $m_\mu^U(x,S)$ is $\F(U)$--measurable for any $x$. Moreover, the map $x\to m_\mu^U(x,S)$ is $L_\mu-$Lipschitz continuous on $\R^d$ because so is $m_\mu(\cdot,S)$. 

Let us now check that $m_\mu\leq m_\mu^U$. Fix $\omega\in \Omega$. As the family $\left\{\omega_n=(\sigma_n,H_n)\right\}$ is dense in $\Omega$, there exists a subsequence $(\omega_{n'}=(\sigma_{n'},H_{n'}))$ which converges to $\omega=(\sigma,H)$ locally uniformly. So, for any $k\in \N^*$, there exists $n'_k\geq k$ such that $\omega$ belongs to $B_{U,k}(\omega_{n'_k})$.  In particular, for any $x\in \R^d$, 
$$
\sup_{n\in \N}  m_\mu(x,S,\omega_n)\indc_{B_{U,k}(\omega_n)}\geq m_\mu(x,S,\omega_{n'_k}).
$$
Since $\{ m_\mu(\cdot, S,\omega_{n'_k})\}_{k\in\N}$ converges locally uniformly to $m_\mu(\cdot,\omega)$ by the local uniform convergence of $(\omega_{n'_k})$ to $\omega$, we obtain the desired inequality upon sending $k\to +\infty$.

We now prove the last statement of the lemma. Fix $\omega=(\sigma,H)\in \Omega$ and $x\in \R^d$. For any $k\in \N$, let $n_k\in \N$ be such that $\omega\in B_{U,k}(\omega_{n_k})$ and 
\begin{equation*}
\sup_{n\in \N}  m_\mu(x,S,\omega_n)\indc_{B_{U,k}(\omega_n)}  \leq m_\mu(x,S,\omega_{n_k})+ \frac1k.
\end{equation*}
By uniform continuity of the elements of $\Omega$, there exists $\omega'=(\sigma',H')\in \overline \Omega$ and a subsequence $\{\omega_{n_k}\}_{k\in\N}$ which converges locally uniformly to $\omega'$. It follows that  $\{m_\mu(\cdot,S,\omega_{n_k})\}_{k\in\N}$ converges locally uniformly to 
$m_\mu(\cdot, S,\omega')$ and therefore
$$
m_\mu^U(x,S,\omega)=\lim_{k\to\infty}  m_\mu(x,S,\omega_{n_k})= m_\mu(x,S, \omega').
$$
Using $\omega\in B_{U,k}(\omega_{n_k})$, we deduce that $(\sigma',H')=(\sigma,H)$ in $\R^d\times U$.
\end{proof}

We next show that $m_\mu^U$ is a good approximation of $m_\mu$ in the sublevel sets of $m_\mu$ which are contained in $U$. 

\begin{lemma}
\label{l.mvsmmu1}
With $\overline a_\mu$ defined as in Proposition \ref{p.localizationPC}, fix $t\geq \overline a_\mu$ and assume that 
\begin{equation*}
\left\{m_\mu(\cdot,S)\leq t\right\}\subseteq U.
\end{equation*}
Then  
\begin{equation*}
m_\mu(\cdot,S)\leq m_\mu^U(\cdot,S)\leq m_\mu(\cdot,S)+l_\mu \quad \mbox{in} \  
\left\{m_\mu(\cdot,S)\leq t-\bar a_\mu\right\},
\end{equation*}
and, for any $s\in [0,t-\bar a_\mu]$, 
\begin{equation*}
\{m_\mu^U(\cdot,S)\leq s\} \subset  \{m_\mu(\cdot,S)\leq s\} \subset \{m_\mu^U(\cdot,S)\leq s\}+ \overline B_3.  
\end{equation*}
\end{lemma}
\begin{proof} 
Let $\omega\in \Omega$ and $x\in \{m_\mu(\cdot,S, \omega)\leq t-\bar a_\mu\}$. According to Lemma \ref{l.basicpptmU} there exists $\omega'=(\sigma',H')\in  \overline \Omega$ such that $m_\mu^U(x,S,\omega)= m_\mu(x,S,\omega')$ and $(\sigma',H')=(\sigma,H)$ in $\R^d\times U$. 

\smallskip

By assumption, we have $(\sigma',H')=(\sigma,H)$ in $U \supseteq \{m_\mu(\cdot,S, \omega)\leq t\}$ and thus
Proposition~\ref{p.localizationPC} yields
\begin{equation*}
\left| m_\mu(\cdot,S,\omega)-m_\mu(\cdot,S,\omega')\right|\leq l_\mu \quad \mbox{in} \  \left\{m_\mu(\cdot,S,\omega)\leq t-\bar a_\mu\right\}.
\end{equation*}
Applying this to $x$, we obtain $|m_\mu(x,S,\omega)-m_\mu^U(x,S,\omega)|\leq l_\mu$.

\smallskip

Now fix $s\leq t-\bar a_\mu$. We already know from Lemma~\ref{l.basicpptmU} that $m_\mu(\cdot,S,\omega)\leq m_\mu^U(\cdot,S,\omega)$, so that 
$$
\{m_\mu^U(\cdot,S,\omega)\leq s\}\subset \{m_\mu(\cdot,S,\omega)\leq s\}.  
$$
Conversely, if $x$ belongs to $\{m_\mu(\cdot,S,\omega)\leq s\}$, then, by \eqref{e.movefron}, one can find $y\in \R^d$ such that $|y-x|\leq 3$ and $m_\mu(y,S,\omega)\leq s-l_\mu$. Let $\omega'=(\sigma',H')\in \overline \Omega$ be associated to $y$ as in Lemma \ref{l.basicpptmU}: $m_\mu^U(y,S,\omega)=m_\mu(y,S,\omega')$ and $(\sigma',H')=(\sigma,H)$ in $\R^d\times U$. Then, by Proposition \ref{p.localizationPC}, $m_\mu(y,S,\omega')\leq m_\mu(y,S,\omega)+l_\mu\leq s$. So $m_\mu^U(y,S,\omega)\leq s$, which  proves the second part of the lemma. 
\end{proof}

We next obtain a result like the previous lemma, with the important difference that the hypothesis requires only that the $t$--sublevel set of $m_\mu^U$ is contained in $U$, rather than the $t$--sublevel set of $m_\mu$. We define a new constant 
\begin{equation*}
a_\mu:=2(L_\mu+l_\mu)+\bar a_\mu(2+L_\mu l_\mu^{-1}),
\end{equation*}
where $\bar a_\mu$ is defined in Proposition \ref{p.localizationPC}. Note that, as $\bar a_\mu = C(\mu l_\mu)^{-1}$, we have that $a_\mu= C\mu^{-1} l_\mu^{-2}$. 

\begin{lemma}
\label{l.localization2} 
Fix $t\geq a_\mu$ and assume that 
\begin{equation*}
\left\{m_\mu^U(\cdot,S)\leq t\right\}\subseteq U
\end{equation*}
Then  
\begin{equation*}
m_\mu(\cdot,S)\leq m_\mu^U(\cdot,S)\leq m_\mu(\cdot,S)+l_\mu \quad \mbox{in} \  
\left\{m_\mu^U(\cdot,S)\leq t-\bar a_\mu\right\},
\end{equation*}
and, for any $s\in [0,t-\bar a_\mu]$, 
\begin{equation*}
\{m_\mu^U(\cdot,S)\leq s\} \subset  \{m_\mu(\cdot,S)\leq s\} \subset \{m_\mu^U(\cdot,S)\leq s\}+ \overline B_3.  
\end{equation*}
\end{lemma}
\begin{proof}
Let $s$ be the largest real number such that $\{m_\mu(\cdot,S)\leq s\}\subset \{m_\mu^U(\cdot,S)\leq t\}$. Then, since 
$$
\{m_\mu(\cdot,S)\leq s\}\subset \{m_\mu^U(\cdot,S)\leq t\}\subset U,
$$
Lemma~\ref{l.mvsmmu1} implies that 
\begin{equation}\label{e.khjerbvlzbres}
|m_\mu(\cdot,S)-m_\mu^U(\cdot,S)|\leq l_\mu \quad \mbox{in} \  
\left\{m_\mu(\cdot,S)\leq s-\bar a_\mu\right\}.
\end{equation}
and, for any $\tau \in [0,s-\bar a_\mu]$, 
$$
\{m_\mu^U(\cdot,S)\leq s\} \subset  \{m_\mu(\cdot,S)\leq s\} \subset \{m_\mu^U(\cdot,S)\leq s\}+ \overline B_3. 
$$

It remains to show that $\{m_\mu^U(\cdot,S)\leq t- a_\mu\}\subset \{m_\mu(\cdot,S)\leq s\}$. Observe that, by the definition of $s$, there exists $x$ such that 
$m_\mu(x,S)= s$ and $m_\mu^U(x,S) =t$. 
By \eqref{e.movefron}, there exists $y\in \R^d$ such that $|x-y|\leq l_\mu^{-1}\bar a_\mu+2$ and $m_\mu(y,S)= s-\bar a_\mu$. Then, by the Lipschitz estimates,  
\begin{multline*}
|s-t| 
\leq \bar a_\mu +|s-\bar a_\mu-t|
\leq \bar a_\mu+ |m_\mu(y,S)-m_\mu^U(y,S)|+|m_\mu^U(y,S)-m_\mu^U(x,S)|\\
 \leq \bar a_\mu+l_\mu +L_\mu ( l_\mu^{-1}\bar a_\mu +2).
\end{multline*}
Let now $\tau$ be the largest real number such that 
\begin{equation*}
\{m_\mu^U(\cdot,S)\leq \tau\}\subset \{m_\mu(\cdot,S)\leq s-\bar a_\mu\}.
\end{equation*}
Then there exists $x'$ such that 
$m_\mu^U(x')= \tau$ and $m_\mu(x')=s-\bar a_\mu$. So, by \eqref{e.khjerbvlzbres}, $|\tau-s+\bar a_\mu|\leq l_\mu$. This gives the result since
\begin{equation*}
 a_\mu = ( \bar a_\mu+l_\mu +L_\mu (l_\mu^{-1}\bar a_\mu+2)) + (\bar a_\mu +l_\mu). \qedhere 
\end{equation*}
\end{proof}

We conclude this subsection by slightly modifying the statement of the previous lemma, putting it in a form better suited for our purposes in the following subsection.  

\begin{proposition}
\label{p.localization}
Suppose that $S\subseteq U\subseteq \Rd$ and $t\geq 1$  satisfy 
\begin{equation*}
\left\{ x \in \R^d \,:\,  m_\mu^U (x,S) \leq t \right\} \subseteq U_{R_0}, 
\end{equation*}
where $R_0:= 5 + l_\mu^{-1}a_\mu$. (Note that $R_0= C\mu^{-1} l_\mu^{-3}$.) Then 
\begin{equation}
\label{e.localization}
0\leq m_\mu^U(\cdot,S)  - m_\mu(\cdot,S) \leq l_\mu  \quad \mbox{in} \ \left\{ m^U_\mu(\cdot,S) \leq t  \right\}
\end{equation}
and, for any $s\in [0,t]$, 
\begin{equation} \label{e.local-levelsets}
\left\{ m_\mu^{U}(\cdot,S)  \leq s \right\} \subseteq \left\{ m_\mu(\cdot,S)  \leq s \right\} 
\subseteq\left\{ m_\mu^{U}(\cdot,S)  \leq s \right\} + \overline B_3.
\end{equation}
\end{proposition}
\begin{proof}
In view of Lemma~\ref{l.localization2}, one just needs to check that, under our assumptions, $\{m_\mu^U(\cdot)\leq t+a_\mu\}\subset U$. Here is the proof: as $\{m_\mu^U(\cdot)\leq t\}\subset U$, Lemma~\ref{l.localization2} implies that 
$$
\{m_\mu\leq t-a_\mu\}\subset \{m_\mu^U\leq t-a_\mu\}+\overline B_3.
$$
Thus, from \eqref{e.movefron}, 
\begin{equation*}
\{m_\mu^U\leq t\}\subset \{m_\mu\leq t\}\subset \{m_\mu\leq t-a_\mu\} +\overline B_{l_\mu^{-1}a_\mu+2}
\subset \{m_\mu^U\leq t-a_\mu\}+\overline B_{R_0} \subset U.\qedhere
\end{equation*}
\end{proof}

\subsection{Construction of the martingale}
In this subsection, we perform a construction similar to the one in~\cite[Section 3.3]{AC}. The eventual goal is to define a filtration $\{ \mathcal G_{t}\}_{t\geq 0}$ on $\Omega$ so that the martingale $\E \left[ m_\mu(x,S) \,\vert\, \G_{t} \right]$ has bounded increments in $t$, almost surely with respect to~$\P$, permitting the application of Azuma's inequality. As in~\cite{AC}, the filtration $\G_t$ is a ``perturbation" of the smallest $\sigma$--algebra which makes the $t$--sublevel set of $m_\mu$ measurable, but nevertheless is sufficiently localized in its dependence on the coefficients that we may make use of independence. 

\smallskip

We begin by introducing a discretization of the set of compact subsets of $\Rd$ which contain $S$. We denote this set by
\begin{equation*}
\K:= \mbox{the set of compact $K\subseteq \Rd$ such that $S \subseteq K$.}
\end{equation*}
We endow $\K$ with the Hausdorff metric $\dist_H$, which is defined by
\begin{align*}
\dist_H\left(K,K'\right)  := & \max \left\{ \inf_{x\in K} \sup_{y\in K'} |x-y|, \, \inf_{y\in K'} \sup_{x\in K} |x-y| \right\} \\ 
 = & \inf\left\{ r>0 \,:\, K \subseteq K'+\overline B_r \ \mbox{and} \ K' \subseteq K+\overline B_r \right\}. 
\end{align*}
Since the metric space~$(\K,\dist_H)$ is locally compact, it follows that there exists a pairwise disjoint partition~$\left( \Gamma_i \right)_{i\in\N}$ of $\K$ into Borel subsets of~$\K$ such that~$\diam_H(\Gamma_i) \leq 1$ for each~$i\in\N$. For each $i\in\N$, we define
\begin{equation*}
K_i := \overline{\bigcup_{K\in \Gamma_i} K}+\overline B_1.
\end{equation*}
By definition, $K_i$ has the interior ball condition of radius $1$ (i.e. satisfies \eqref{e.IntBallPpt}) and 
\begin{equation*}
K \in \Gamma_i \quad \implies \quad K \subseteq K_i \subseteq K+\overline B_2. 
\end{equation*}
We define compact sets 
\begin{equation*}
K_i \subseteq K_i' \subseteq K_i'' \subseteq \tilde K_i 
\end{equation*}
for each $i\in\N$ by
\begin{equation*}
K_i' := K_i + \overline B_{R_0}, \quad K_i'':= K_i' + \overline B_{14} \quad \mbox{and} \quad \tilde K_i := K_i''+\overline B_1. 
\end{equation*}
We enlarge $\Gamma_i$ by setting
\begin{equation*}
\tilde \Gamma_i := \left\{ K\in \K \,:\, K \subseteq K_i \subseteq K + \overline B_4 \right\}. 
\end{equation*}
We next define the following subsets of $\Omega$, for each $t>0$ and $i\in\N$:
\begin{equation*}
F_i(t) :=  \mbox{the event that} \ \  \left\{ x \in \R^d \,:\, m_\mu^{K_i'} (x,S) \leq t  \right\}  \in \tilde \Gamma_i.
\end{equation*}
We next show that the events $\{ F_i(t) \}_{i\in\N}$ covers $\Omega$. 
\begin{lemma}
\label{l.martin2}
For every $t>0$, 
\begin{equation*}
\bigcup_{i\in\N} F_i(t) = \Omega. 
\end{equation*}
\end{lemma}
\begin{proof}
Fix $i\in\N$ such that $\left\{ x \in \Rd \,:\, m_\mu(x,S) \leq t \right\} \in \Gamma_i$. We claim that $F_i(t)$ holds. Our assumption implies that
\begin{equation*} \label{}
\left\{ x \in K_i' \,:\, m_\mu^{K_i'}(x,S) \leq t \right\} \subseteq K_i. 
\end{equation*}
Thus Proposition~\ref{p.localization} is applicable and we deduce from~\eqref{e.local-levelsets} for $s=t$ and our choice of~$i$ that $\left\{ m_\mu^{K_i'}(\cdot,S) \leq t \right\} \in \tilde \Gamma_i$. That is, $F_i(t)$ holds. 
\end{proof}

We next show that, if $F_i(t)$ holds, then the $t$--sublevel set of $m_\mu$ is close to~$K_i$.

\begin{lemma}
\label{l.martin1}
Suppose that $t>0$ and $F_i(t)$ holds. Then 
\begin{equation} \label{kjhbzkheb1}
\left\{ m_\mu(\cdot,S) \leq t \right\} \subseteq K_i+ \overline B_3 \quad \mbox{and} \quad K_i\subseteq \left\{ m_\mu(\cdot,S) \leq t \right\}  +\overline B_4,
\end{equation}
\begin{equation} \label{kjhbzkheb2}
\sup_{K_i}\left( m_\mu^{K_i'}(\cdot,S) - m_\mu(\cdot,S) \right)  \leq 4L_\mu+l_\mu
\end{equation}
and
\begin{equation} \label{kjhbzkheb3}
\sup_{\partial K_i}\left| m_\mu(\cdot,S) - t \right|  \leq 8L_\mu+2l_\mu.
\end{equation}
\end{lemma}
\begin{proof}
If $F_i(t)$ holds, then 
\begin{equation} \label{e.martin12}
\left\{ m_\mu^{K_i'}(\cdot,S)  \leq t \right\} \subseteq K_i \subseteq \left\{ m_\mu^{K_i'}(\cdot,S)  \leq t \right\} + \overline B_4.
\end{equation}
Note that the first inclusion of~\eqref{e.martin12} ensures the applicability of Proposition~\ref{p.localization}. The inclusions in \eqref{kjhbzkheb1} are therefore obtained from~\eqref{e.local-levelsets} and~\eqref{e.martin12}. The  inequality~\eqref{kjhbzkheb2} is obtained from~\eqref{e.localization},~\eqref{e.martin12} and the Lipschitz estimate. To prove~\eqref{kjhbzkheb3}, note that if $x\in \partial K_i$, then by the first inclusion in~\eqref{e.martin12}, we have $m_\mu^{K_i'}(x)\geq t$. By the second inclusion~\eqref{e.martin12}, there exists $y\in \R^d$ such that $m_\mu^{K_i'}(y)\leq t$ and $|x-y|\leq 4$. Thus, by the Lipschitz estimate, $m_\mu^{K_i'}(x)\leq t+4L_\mu$. Combining this with the second conclusion proves the claim. 
\end{proof}

We now create a partition of $\Omega$ by defining
\begin{equation*}
E_1(t):= F_1(t), \quad E_{i+1}(t) : = F_{i+1}(t) \setminus \left( E_1(t) \cup\cdots\cup E_i(t) \right), \quad i\in\N. 
\end{equation*}
By the previous lemma, $\{ E_i(t) \}_{i\in\N}$ is a pairwise disjoint partition  of~$\Omega$. 

\smallskip

An important property of $E_i(t)$ is that it is measurable with respect to the restriction of the coefficient fields to $K_i''$, which is the assertion of the following lemma. 

\begin{lemma}
\label{l.martin3}
For every $0< s \leq t$ and $i,j\in\N$, 
\begin{equation*}
F_i(s) \cap F_j(t) \neq \emptyset \quad \implies \quad K_i'\subset K_j''\qquad {\rm and}\qquad E_i(s) \in \mathcal F\left( K_j''\right).
\end{equation*}
In particular, $E_i(t) \in \mathcal F\left( K_i''\right)$.
\end{lemma}
\begin{proof}
According to Lemma~\ref{l.martin1}, if $F_i(s) \cap F_j(t) \neq \emptyset$ for some $0<s\leq t$, then 
\begin{equation*} \label{}
K_i\subseteq \left\{ m_\mu(\cdot,S) \leq s \right\}  +\overline B_4  \subseteq \left\{  m_\mu(\cdot,S) \leq t \right\}  +\overline B_4 \subseteq K_j + \overline B_7.
\end{equation*}
In particular, $F_i(s) \in \F\left( K_j' + \overline B_7 \right)$. 
We also obtain the expression
\begin{equation*} \label{}
E_i (s) = F_i(s) \setminus \bigcup_{n\in D(i)} F_n(s),
\end{equation*}
where
\begin{equation*}
D(i) : = \left\{ n \in\N\,:\, 1 \leq n \leq i-1, \ K_i\subseteq K_n + \overline B_7, \ \mbox{and} \ K_n \subseteq K_i + \overline B_7\right\}.
\end{equation*}
Notice that $n\in D(i)$ implies $K_n' \subseteq K_j'+\overline B_{14} = K_j''$ and thus $F_n(s) \in \F\left( K_j'' \right)$. It follows that $E_i(s) \in \F \left(  K_j'' \right)$, as desired.
\end{proof}

In view of Lemma~\ref{l.martin2}, we may define, for $t>0$, a random element $S_t$ of $\K$ by 
\begin{equation*}
S_t:= K_i \quad \mbox{if} \ E_i(t) \ \mbox{holds}, \ i\in\N.
\end{equation*}
Note that $S_t$ is an approximation of the set $\left\{ x\in \H \,: \, m_\mu(x,S) \leq t \right\}$, but with more local dependence on the coefficient fields, as witnessed by Lemma~\ref{l.martin3}. 

\smallskip

We now define the filtration $\{\G_t\}_{t\geq 0}$ by~$\G_0:= \{ \emptyset, \Omega\}$ and, for $t>0$, by
\begin{multline*}
\G_t:=    \mbox{$\sigma$--algebra on $\Omega$ generated by  events of the form $E_i(s) \cap F$,} \\
 \mbox{where $0<s\leq t$, $i\in\N$ and $F\in \mathcal \F(K_i'')$.}
\end{multline*}
The martingale we are interesting in is $\E \left[ m_\mu(x,S) \,\vert\, \G_t \right]$ for a fixed~$x\in\Rd$. The eventual goal, which is completed in the following subsection, is to show that this martingale possesses bounded increments and to deduce from Azuma's inequality bounds on its fluctuations for $t\gg1$ large. We conclude this subsection with three lemmas containing some estimates we need. 

\smallskip

\begin{lemma}
\label{l.martingale1}
For every $t>0$ and $x\in \Rd \setminus S_t$, 
\begin{equation*} \label{}
\left| m_\mu(x,S) - (t+m_\mu(x,S_t)) \right| \leq 8L_\mu+2l_\mu.
\end{equation*}

\end{lemma}
\begin{proof}
By the maximality of $m_\mu(\cdot,S_t)$ (the last statement of Theorem~\ref{t.LipschitzAppenA1}), we have that, for every $x\in \R^d\setminus S_t$,
\begin{equation}
m_\mu(x,S) - \sup_{y\in \partial S_t} m_\mu(y,S) \leq m_\mu(x,S_t).
\end{equation}
By the maximality of $m_\mu(\cdot,S)$, we have, for every $x\in \Rd \setminus S_t$,
\begin{equation*}
m_\mu(x,S_t) + \inf_{y\in \partial S_t} m_\mu(y,S) \leq m_\mu(x,S).
\end{equation*}
Finally, we note by  Lemma~\ref{l.martin1} that, for every $t >0$,
\begin{equation*}
\sup_{y\in \partial S_t}  \left| t -m_\mu(y,S) \right| \leq 8L_\mu+2l_\mu.
\end{equation*}
Combining the above yields the lemma. 
\end{proof}

\begin{lemma}
\label{l.martingale2}
For every $0< s\leq t$ and $x\in\Rd$,
\begin{equation*} \label{}
\left| m_\mu(x,S)  - \E \left[ m_\mu(x,S) \, \vert \, \G_t \right]  \right| \indc_{ \{ x \in S_s \}} \leq 8L_\mu+2l_\mu.
\end{equation*}
In particular, if $t\geq L_\mu \dist(x,S)$, then 
$$
\left| m_\mu(x,S)  - \E \left[ m_\mu(x,S) \, \vert \, \G_t \right]  \right|  \leq 8L_\mu+2l_\mu\qquad {\rm a.s.}
$$
\end{lemma}
\begin{proof}
Fix $x\in \Rd$ and $0<s\leq t$. Define a random variable $Z$ by
\begin{equation*}
Z:= \sum_{i\in \{ j\in \N\,:\, x\in K_j\} } m_\mu^{K_i' } (x,S) \indc_{E_i(s)}.
\end{equation*}
Notice that $Z$ is $\G_s$--measurable by the definition of the filtration, since $m_\mu^{K_i'} (x,S)$ is $\F\left(K_i'' \right)$--measurable. Since the event that $x\in S_s$ is the union of $E_i(s)$ over $i\in \{ j\in \N\,:\, x\in K_j\}$, we obtain from Lemma~\ref{l.martin1} that 
\begin{equation*}
\left| Z - m_\mu(x,S) \indc_{\{ x\in S_s \}} \right| \leq \sum_{i\in \{ j\in \N\,:\, x\in K_j\} } \left| m_\mu(x,S) - m_\mu^{K_i' } (x,S) \right| \indc_{E_i(s)}  \leq 4L_\mu+l_\mu.
\end{equation*}
Since $Z$ and $ \indc_{\{ x\in S_s\}}$ are $\G_t$--measurable, we have 
\begin{align*}
\lefteqn{ 
\left| m_\mu(x,S) - \E \left[ m_\mu(x,S) \,\vert\, \G_t \right] \right| \indc_{\{ x\in S_s\}} 
} \qquad & \\
& \leq \left| Z  - m_\mu(x,S) \right|  \indc_{\{ x\in S_s\}} +  \left| \E \left[ Z  - m_\mu(x,S) \indc_{\{ x\in S_s\}} \,\vert\, \G_t \right] \right| \\
& \leq 2 \left\|Z  - m_\mu(x,S)\indc_{\{ x\in S_s\}} \right\|_{L^\infty(\Omega,\P)} \\
& \leq 2(4L_\mu+l_\mu).
\end{align*}
This completes the proof of the first statement.

In order to check that, for $t\geq L_\mu \dist(x,S)$, $m_\mu(x)$ is ``almost" $\G_t-$measurable, let us note that, in $E_i(t)$, we have 
$$
\{ m_\mu^{K_i'} (\cdot,S)\leq t\} \subset K_i, 
$$
where, in view of \eqref{control2},  
$$
m_\mu^{K_i'}(x,S)\leq L_\mu \dist(x,S)\leq t.
$$
So, in $E_i(t)$, $x$ belongs to $K_i$,  $\P$--a.s., so that $x\in S(t)$, $\P$--a.s.
\end{proof}

\begin{lemma}
\label{l.martingale3}
For every $t>0$ and $x\in\Rd$, 
\begin{equation} \label{e.margingale3}
\left| \E \left[ m_\mu(x,S_t) \,\vert\, \G_t \right] - \sum_{i\in\N} \E \left[ m_\mu(x,K_i) \right] \indc_{E_i(t)} \right| \leq 2L_\mu(R_0+15).
\end{equation}
\end{lemma}
\begin{proof}
We first argue that, for each $x\in\Rd$, $i\in\N$, $t>0$,
\begin{equation}
\label{e.martin3wts}
\E \left[  m_\mu\left(x,\tilde K_i \right)  \indc_{E_i(t)}  \,\vert\, \G_{t} \right] = \E \left[ m_\mu\left(x,\tilde K_i \right)\right] \indc_{E_i(t)}.
\end{equation}
Since $E_i(t) \in \G_t$ by definition, to establish~\eqref{e.martin3wts} it suffices to show that, for every $A\in \G_t$, 
\begin{equation}
\label{e.martin3wts2}
\E \left[ m_\mu\left(x,\tilde K_i \right)  \indc_{A \cap E_i(t)}  \right] = \E  \left[ m_\mu\left(x,\tilde K_i \right)    \right] \P \left[ A \cap  E_i(t)   \right].
\end{equation}
We obtain~\eqref{e.martin3wts2} from the fact that
\begin{equation}
\label{e.martin3wts3}
A \in \G_t \quad \implies  \quad A \cap E_i(t) \in \F\left(K_i'' \right),
\end{equation}
which we will check below, the fact that $m_\mu(x,\tilde K_i)$ is $\F\left( \Rd \setminus \tilde K_i \right)$--measurable, and the independence of $\F\left( \Rd \setminus \tilde K_i \right)$ and $\F\left( K_i'' \right)$.

\smallskip

We now give the proof of~\eqref{e.martin3wts3}. We may assume that $A$ takes the form $A= F \cap E_j(s)$, with $j\in\N$, $F\in \F\left(K_j''  \right)$ and $s\in (0,t]$, since such events generate~$\G_t$. Then $A\cap E_i(t) = F \cap E_j(s)\cap E_i(t)$. Now, either $E_j(s)\cap E_i(t)  = \emptyset$ and there is nothing more to show, else or $K_j' \subseteq K_i''$ by Lemma~\ref{l.martin3}. In the latter case, we have $F \cap E_j(s) \in \F\left(K_j'  \right) \subseteq \F\left( K_i''\right)$. By Lemma~\ref{l.martin3} again, $E_i(t) \in \F(K_i'')$ and thus $F \cap E_j(s)\cap E_i(t) \in \F\left( K_i'' \right)$, as desired. This completes the proof of~\eqref{e.martin3wts3}.

\smallskip

We now complete the proof, deriving~\eqref{e.margingale3} from~\eqref{e.martin3wts} using the Lipschitz estimates: we have 
\begin{align*}
\E \left[ m_\mu(x,S_t) \,\vert\, \G_t \right] 
& = \sum_{i\in \N} \E \left[ m_\mu(x,K_i) \indc_{E_i(t)} \,\vert\, \G_t \right]  \\
& \leq \sum_{i\in \N} \E \left[ m_\mu(x,\tilde K_i) \indc_{E_i(t)} \,\vert\, \G_t \right] + L_\mu(R_0+15) \\
& = \sum_{i\in \N} \E \left[ m_\mu(x,\tilde K_i)\right]  \indc_{E_i(t)} + L_\mu(R_0+15) \\
& \leq  \sum_{i\in \N} \E \left[ m_\mu(x, K_i)\right]  \indc_{E_i(t)} + 2L_\mu(R_0+15).
 \end{align*}
The reverse inequality is proved in the same way. 
\end{proof}

\subsection{Proof of the fluctuations estimate}
Using the results of the previous subsection, we are now ready to derive Proposition~\ref{p.fluctuations} from Azuma's inequality. 

\begin{proof}[{{Proof of Proposition~\ref{p.fluctuations}}}]
We fix $x\in \Rd$ and consider the $\G_t$--adapted martingale $\{X_t\}_{t\geq0}$ defined by 
\begin{equation*}
X_t:= \E \left[ m_\mu(x,S) \,\vert\, \G_t \right] - \E \left[  m_\mu(x,S) \right].
\end{equation*}
Note that $X_0 = 0$ since $\G_0=\{ \emptyset, \Omega\}$.

\smallskip

The main step in the argument is to show, using Lemmas~\ref{l.martingale1}, \ref{l.martingale2} and~\ref{l.martingale3} that, for every $t,s>0$, 
\begin{equation}
\label{e.boundeddiff}
\left| X_t - X_s \right| \leq C\left( l_\mu^{-1} |s-t| +1+\mu^{-1} l_\mu^{-3}\right), 
\end{equation}
where $C$ is bounded for bounded $\mu$'s. 
We may assume that $s\leq t$. Then the event that $x\in S_s$ is $\G_s$--measurable and hence 
\begin{equation*}
\left\{ \begin{aligned}
& X_t = \E \left[ m_\mu(x,S) \,\vert\, \G_t \right]\indc_{\{x\in S_s\}} + \E \left[ m_\mu(x,S)\indc_{\{x\not\in S_s\}}  \,\vert\, \G_t \right] -\E \left[  m_\mu(x,S) \right], \\
& X_s = \E \left[ m_\mu(x,S) \,\vert\, \G_s \right]\indc_{\{x\in S_s\}} + \E \left[ m_\mu(x,S)\indc_{\{x\not\in S_s\}}  \,\vert\, \G_s \right] - \E \left[  m_\mu(x,S) \right].
\end{aligned}\right. 
\end{equation*}
Subtracting these and applying Lemma~\ref{l.martingale2}, we get 
\begin{equation*}
\left| X_s - X_t \right| \leq 4( 4L_\mu+l_\mu) + \left|\E \left[ m_\mu(x,S)\indc_{\{x\not\in S_s\}}  \,\vert\, \G_t \right]- \E \left[ m_\mu(x,S)\indc_{\{x\not\in S_s\}}  \,\vert\, \G_s \right] \right|.
\end{equation*}
Applying Lemma~\ref{l.martingale1} twice, we obtain
\begin{equation*}
\left| X_s - X_t \right| \leq 12( 4L_\mu+l_\mu) + \left| \E \left[ m_\mu(x,S_s)\indc_{\{x\not\in S_s\}} \,\vert\, \G_t\right] - \E \left[ m_\mu(x,S_s)\indc_{\{x\not\in S_s\}} \,\vert\, \G_s\right]  \right| .
\end{equation*}
Let us now estimate $\dist_H(S_s,S_t)$. We set $\Re_{\mu,t}:= \{m_\mu(\cdot,S)\leq t\}$. By Lemma~\ref{l.martin1} and the growth of $t\to \Re_{\mu,t}$ in \eqref{e.movefron}, we have 
\begin{align*}
\lefteqn{
\dist_H(S_s,S_t) = \sum_{i,j\in \N} \dist_H(K_i,K_j) \indc_{E_i(s)\cap E_j(t)}
} \qquad & \\
& \leq 
\sum_{i,j\in \N} \left( \dist_H(K_i,\Re_{\mu,s}) +  \dist_H(\Re_{\mu,s},\Re_{\mu,t})+  \dist_H(\Re_{\mu,t}, K_j) \right)  \indc_{E_i(s)\cap E_j(t)}  \\
& \leq  l_\mu^{-1} |s-t|+10. 
\end{align*}
From the Lipschitz estimates we obtain therefore that 
\begin{equation*}
\left| m_\mu( x,S_s) - m_\mu(x,S_t) \right|  \indc_{\{x\not\in S_s\}} \leq L_\mu\left(l_\mu^{-1} |s-t|+10\right).
\end{equation*}
Plugging this inequality into  the estimate of $\left| X_s - X_t \right|$ and using that the event~$\{x\not\in S_s\}$ is  $\G_s$--measurable,  we obtain that 
\begin{equation*}
\left| X_s - X_t \right| \leq  \left| \E \left[ m_\mu(x,S_t) \,\vert\, \G_t\right] - \E \left[ m_\mu(x,S_s) \,\vert\, \G_s\right]  \right|
+L_\mu l_\mu^{-1} |s-t|+C .
\end{equation*}
To complete the proof of~\eqref{e.boundeddiff}, we are left to show that
\begin{equation}
\label{e.boundeddiffwts}
\left| \E \left[ m_\mu(x,S_t) \,\vert\, \G_t\right] - \E \left[ m_\mu(x,S_s) \,\vert\, \G_s\right]  \right| \leq 
L_\mu  (l_\mu^{-1} |s-t|+40+4R_0).
\end{equation}
For this, we combine Lemma~\ref{l.martingale3}, the estimate for  $\dist_H(S_s,S_t) $ established above and with the Lipschitz estimate, to get
\begin{align*}
\lefteqn{
\left| \E \left[ m_\mu(x,S_t) \,\vert\, \G_t\right] - \E \left[ m_\mu(x,S_s) \,\vert\, \G_s\right]  \right|
}
\qquad & \\
& \leq 4L_\mu(R_0+15)+ \sum_{i,j\in \N} \left| \E \left[ m_\mu(x,K_j) \right] - \E \left[ m_\mu(x,K_i) \right]  \right| \indc_{E_i(s) \cap E_j(t)}  \\
& \leq 4L_\mu(R_0+15)+  L_\mu \sum_{i,j\in\N}  \dist_H(K_i,K_j) \indc_{E_i(s) \cap E_j(t)} \\
& \leq 4L_\mu(R_0+15)+  L_\mu \dist_H(S_s,S_t)  \\
& \leq 4L_\mu(R_0+15)+  L_\mu  (l_\mu^{-1} |s-t|+10), 
\end{align*}
as desired. This yields~\eqref{e.boundeddiff} because $R_0= C\mu^{-1} l_\mu^{-3}$.

\smallskip

We now complete the proof of the proposition. We fix $T\geq 1$ large and set $\tau= \mu^{-1}l_\mu^{-2}$ and $N=T/\tau$. Applying Azuma's inequality yields that for every $\lambda>0$, 
\begin{equation*}
\P \left[ \left| X_T - X_0 \right| > \lambda \right] \leq 2 \exp\left(-\frac{\lambda^2}{CN} \right).
\end{equation*}
Using $X_0=0$ and the choice of $N$, we obtain, in view of~\eqref{e.boundeddiff}, that for every $\lambda>0$ and $T\geq 1$,
\begin{equation*}
\P \left[ \left| X_T \right| > \lambda \right] \leq 2\exp\left( -\frac{\mu l_\mu^{2} \lambda^2}{C T} \right). 
\end{equation*}
If we choose $T= L_\mu \dist(x,S)$, we have by Lemma \ref{l.martingale2}: 
$$
\left|X_T-\left(m_\mu(x,S)  - \E \left[ m_\mu(x,S) \right]\right)\right|\leq C\qquad {\rm a.s.}
$$
Plugging this inequality into the former one, we obtain 
\begin{equation*}
\P \left[ \left| m_\mu(x,S)  - \E \left[ m_\mu(x,S)\right] \right| > \lambda \right] \leq C\exp\left( -\frac{\mu l_\mu^{2} \lambda^2}{C \dist(x,S)} \right). 
\end{equation*}
This is~\eqref{e.fluctuations}. 
\end{proof}

\section{Convergence of the mean distance to a plane}
\label{s.means}

In the previous section, we obtained good control of the stochastic fluctuations of the solutions of the planar metric problem at points far from the boundary plane. To complete the proof of Proposition~\ref{p.planarquantitative}, it remains to study the asymptotic behavior of the quantity $\E \left[ m_\mu(x,\H_e^-) \right]$ as $x\cdot e \to \infty$. The precise statement we need is presented in the following proposition, the proof of which is the focus of this section.

\begin{proposition}
\label{p.statbiasBIS}
For each $L\geq 1$, there exists $C(\data,L)\geq 1$ and, for each $\mu\in (0,L]$ and $e\in\partial B_1$, a positive real number $\overline m_\mu(e)$ such that, for every $x\in \H_e^+$,
\begin{equation*}
\left| \E \left[ m_\mu(x,\H_e^-) \right] - \overline m_\mu (e)(x\cdot e) \right|  \leq C_\mu  (x\cdot e)^{\frac12} \log^{\frac12} (1+(x\cdot e)),
\end{equation*}
where  $C_\mu := C \mu^{-\frac32} \left(1+\left|\log \mu\right|\right)^{-\frac12}$. 
Moreover, $(\mu,e) \mapsto \overline m_\mu(e)$ is continuous on~$(0,\infty) \times \partial B_1$.
\end{proposition}

The proof of Proposition~\ref{p.statbiasBIS} requires a new localization argument which is a generalization of the finite speed of propagation property for first order Hamilton-Jacobi equations. This is presented in the next subsection and the proof of Proposition~\ref{p.statbiasBIS} is given in Section~\ref{s.ssbias}.

\smallskip

Throughout, we fix $L\geq 1$, $\mu \in (0,L]$ and $e \in\partial B_1$ and set $\H^\pm :=  \H_e^\pm$. 

\subsection{Propagation of influence}
\label{s.propagation}
An important property of the planar metric problem is an ``approximate finite speed of propagation" property. To be more precise, what we show is that, while the diffusion term of course creates an infinite speed of propagation, the behavior of the boundary condition outside of a ball centered at $x\in \H^+$ of radius $\gg (x\cdot e)^{\frac 92}$ has essentially negligible influence on the value of $m_\mu(x,\H^-)$. The result is summarized in the following proposition. 

\begin{proposition}
\label{p.finitespeedMCM}
Let $m^1,m^2\in W_{\mathrm{loc}}^{1,\infty}(\H^+)$ be, respectively, a subsolution and supersolution of the equation
\begin{equation}\label{e.MCMgen}
-\tr\left( A(Dm,y) D^2m \right) + H\left( Dm,y \right) = \mu \quad \mbox{in} \ \H^+.
\end{equation}
Suppose also that there exist constants $K,R,{M}\geq 1$ such that, for every $i\in \{1,2\}$ and $x\in\H^+$,
\begin{equation} \label{e.milip}
\esssup_{y\in \H^+} \left| Dm^i(y) \right| \leq K,
\end{equation}
\begin{equation} \label{e.migrowth}
0\leq  m^i(x) \leq {M}+L  |x|
\end{equation}
and
\begin{equation} \label{e.miord}
m^1\leq  m^2  \quad \mbox{in} \ \partial \H^+\cap B_R.
\end{equation}
Then there exists $C(\data, K,L)\geq1$ such that, if $s\geq 1$ and $R\geq C\mu^{-5} (1+{M}+s)^{\frac92}$, then we have
\begin{equation} \label{e.FSconc}
m^1(se)\leq m^2(se)+1.
\end{equation}
\end{proposition}

We expect that the exponent $\frac 92$ appearing in the conclusion of Proposition~\ref{p.finitespeedMCM} is suboptimal. Indeed, our argument for obtaining it is somewhat crude and, for example, the proof is easier and we obtain a better estimates in the semilinear case in which $A(\xi,x)=A(x)$. This is of no consequence for the results in the paper, however, because Corollary~\ref{c.fluctuationsH} provides \emph{exponential} estimates on the fluctuations of $m_\mu(x,\H^-)$ which overwhelms any finite power such as $\frac92$. Thus the statement above turns out to be more than enough for what we need. 

\smallskip

The proof of Proposition~\ref{p.finitespeedMCM} is inspired by the proof of the finite speed of propagation for first-order equations, with the role of time being played by the unit direction~$e$. Indeed, the argument relies on a comparison between the planar metric problem and a time-dependent one which is captured in the following simple lemma. 

\begin{lemma}\label{l.testfunction1} 
Suppose $K\geq 1$ and $m\in W^{1,\infty}_{\mathrm{loc}}(\H^+)$ is a nonnegative subsolution of~\eqref{e.MCMgen} satisfying the Lipschitz bound
\begin{equation} \label{e.m0lip}
\esssup_{x\in \H^+} \left| Dm(x) \right| \leq K.
\end{equation}
Fix $\lambda,\nu>0$ and define $w\in W^{1,\infty}_{\mathrm{loc}}(\H^+ \times(0,\infty))$ by
\begin{equation*}
w(x,t):= -\frac1\lambda \log \left( \exp\left(-\lambda m(x)\right) + \exp\left(-\lambda \nu t \right) \right).
\end{equation*}
Then there exists $C=C( \data )$ such that, if $\nu \leq \mu - C\lambda K^2$, then $w$ is a subsolution of the time-dependent equation
\begin{equation}
\label{e.Wsubsolution}
\partial_t w - \tr\left( A(Dw,x) D^2w \right) + H\left(Dw,x \right) \leq \mu \quad \mbox{in} \ \H^+ \times (0,\infty). 
\end{equation}
\end{lemma}
\begin{proof}  
We give the proof assuming that $m$ is smooth; the general case is obtained by performing analogous computations on smooth test functions, in the usual way. For convenience, denote $Z(x,t) : = \exp\left(-\lambda m(x)\right) + \exp\left(-\lambda \nu t\right)$ so that we may write $w(x,t)=-\lambda^{-1} \log Z(x,t)$. Straightforward computations give
\begin{equation}
\left\{ \begin{aligned}
\partial_t Z(x,t) & = -\lambda \nu\exp\left( -\lambda \nu t \right), \\
D Z(x,t) & = -\lambda \exp\left( -\lambda m(x) \right) Dm(x),  \\
D^2Z(x,t) & = -\lambda \exp\left( -\lambda m(x) \right) \left( D^2m(x) - \lambda Dm(x) \otimes Dm(x) \right)
\end{aligned} \right.
\end{equation} 
and thus
\begin{align*}
\partial_t w(x,t) & = -\frac{1}{\lambda Z(x,t)} \partial_tZ(x,t) = \frac{\nu \exp\left( -\lambda \nu t\right) }{ Z(x,t)},\\
D w(x,t) & = -\frac{1}{\lambda Z(x,t)} DZ(x,t) = \frac{\exp\left( -\lambda m(x) \right)}{Z(x,t)} Dm(x), \\
D^2w(x,t) & = -\frac{1}{\lambda Z(x,t) } \left( D^2Z(x,t) - \frac{1}{Z(x,t)} DZ(x,t) \otimes DZ(x,t) \right) \\
& = \frac{\exp\left( -\lambda m(x) \right) }{Z(x,t)}\left( D^2m(x) -   \frac{\lambda\exp\left( -\lambda \nu t \right)}{Z(x,t)} Dm(x) \otimes Dm(x)  \right) .
\end{align*}
Assembling these together and using the positive homogeneity and nonnegativity of $H$, the $0$-homogeneity of $A(\cdot,x)$ and~\eqref{e.m0lip}, we find that 
\begin{align*}
\lefteqn{ \partial_t w(x,t) -\tr\left( A(Dw,x) D^2w(x,t) \right) + H\left( Dw(x,t),x \right) } \qquad & \\
& \leq \left( \frac{\exp\left(-\lambda \nu t \right)}{Z(x,t)} \right) \nu 
+ \left(\frac{\exp\left( -\lambda m(x) \right)}{Z(x,t)}\right)\left( \frac{\exp\left( -\lambda \nu t \right)}{Z(x,t)} \right) \left( C \lambda \left| Dm(x) \right|^2 \right) \\
& \qquad + \left(\frac{\exp\left( -\lambda m(x) \right)}{Z(x,t)} \right) \left( -\tr\left( A(Dm,x) D^2m(x) \right) + H\left( Dm(x),x \right) \right)  \\
& \leq \mu - ( \mu - \nu) \left( \frac{\exp\left( -\lambda \nu t\right)}{Z(x,t)} \right)  + C \lambda L^2  \left( \frac{\exp\left( -\lambda \nu t \right)}{Z(x,t)} \right).
\end{align*}
Therefore we get~\eqref{e.Wsubsolution} provided we select $\nu \leq \mu - C \lambda L_\mu^2$, as claimed.
\end{proof}

\begin{proof}[{Proof of Proposition~\ref{p.finitespeedMCM}}]  
The proof is broken into six steps. The main part is a comparison argument that is similar to several others appearing later in the paper. It comes in Steps~2-4. We then derive the conclusion of the proposition from the result of the comparison in the last two steps. Throughout, $C$ and $c$ denote positive constants depending only on $(\data,K,L)$ and may vary in each occurrence. 

\smallskip

\emph{Step 1.} We set up the comparison argument. We fix parameters $\lambda,\nu,\delta,\ep,T>0$ (selected below), with $\ep \leq \delta^{\frac13}$, and define 
\begin{equation*} \label{}
w_1(x,t):=-\frac1\lambda \log \left( \exp\left(-\lambda m^1(x)\right) + \exp\left(-\lambda \nu t \right) \right).
\end{equation*}
Provided that $\nu \leq \mu - C\lambda$, Lemma \ref{l.testfunction1} asserts that $w_1$ is a subsolution of
\begin{equation}\label{e.evolMCM}
\partial_t w_1 - \tr\left(A(Dw_1,y) D^2w_1 \right) + H\left( Dw_1,y \right)=0\quad \mbox{in} \  \H^+ \times (0,\infty). 
\end{equation}
Fix another parameter $\eta>1$ (to be chosen below in Step~2) and select a smooth, nonnegative, nondecreasing and convex function $g:\R \to [0,\infty)$ satisfying (here we assume $\ep$ is sufficiently small)
\begin{equation*} \label{}
g(t) = \ep t+ \frac14 \quad \mbox{for} \ t\in [0,\infty), \quad \sup_{t\in\R} g'(t) \leq \ep \quad \mbox{and} \quad  \sup_{t\in\R} g''(t) \leq \ep.
\end{equation*}
Define
\begin{equation*} \label{}
\psi(x,t):=g\big( \left( 1+|x|^2 \right)^{\frac12} -  \left(T-t \right)\big).
\end{equation*}
We note for later use that, with $(\cdot) := \left( 1+|x|^2 \right)^{\frac12} -  \left(T-t \right)$, 
\begin{equation*}
\left\{ 
\begin{aligned}
& \partial_t \psi(x,t) =  g'(\cdot ), \\
& D\psi(x,t) = g'(\cdot )\left( 1+|x|^2 \right)^{-\frac12}x, \\
& D^2\psi(x,t) = g'(\cdot )\left( 1+|x|^2 \right)^{-\frac12}  \left( I_d -  \left( 1+|x|^2 \right)^{-1} x\otimes x \right)+ g''(\cdot )\left( 1+|x|^2 \right)^{-1}x\otimes x,
\end{aligned} \right.
\end{equation*}
and therefore, for every $(x,t) \in\H^+\times [0,\infty)$,
\begin{equation}
\label{e.psiDbounds}
|D\psi(x,t)|\leq \ep \quad \mbox{and} \quad  |D^2\psi(x,t)|\leq C\ep.
\end{equation}
We next introduce the auxiliary function $\Psi:\overline{\H^+} \times\overline{\H^+} \times [0,T] \to \R$ defined by
\begin{equation}\label{e.defPsi}
\Psi (x,y,t) := w_1(x,t)-\eta m^2 (y,t)-\frac{|x-y|^4}{4\delta}-\psi(x,t).
\end{equation}
Observe that $\Psi$ attains its supremum on $\overline{\H^+}\times\overline{\H^+}\times[0,T]$ at some point $(x_0,y_0,t_0)\in\overline{\H^+}\times \overline{\H^+}\times [0,T]$. Indeed, $w_1$ is bounded on this set, $m^2$ is nonnegative and $\psi(x,t)\to+\infty$ as $|x|\to+\infty$ uniformly with respect to $t\in [0,T]$. By the Lipschitz assumption~\eqref{e.milip}, we have 
\begin{equation}\label{e.boundx-y}
|x_0 - y_0|\leq (K\delta)^{\frac13} \leq C\delta^{\frac13}.
\end{equation}
The main claim, which is proved in the following three steps, is that, under a suitable choice of $\eta$, we have
\begin{equation} \label{e.trifecta}
\mbox{either} \quad t_0=0 \quad \mbox{or} \quad  x_0\in \partial \H^+ \quad \mbox{or} \quad y_0\in \partial \H^+.
\end{equation}
To prove~\eqref{e.trifecta}, we proceed by contradiction and assume that $t_0>0$ and $(x_0,y_0) \in\H^+\times\H^+$, which means that $(x_0,y_0,t_0)$ is an interior local maximum of the function~$\Psi$. 

\smallskip

\emph{Step 2.} The application of viscosity solution theoretic comparison machinery. Fix another parameter $\gamma>0$ to be selected below. By the parabolic version of the maximum principle for semicontinuous functions~\cite[Theorem 8.3]{CIL}, there exist symmetric matrices $X,Y\in \R^{d\times d}$ such that
\begin{equation*} \label{}
\left\{ 
\begin{aligned}
& \left(X, \xi_0+D\psi(x_0,t_0) , \partial_t\psi(x_0,t_0)\right)\in \overline{\mathcal P}^{2,+}w_1(x_0,t_0), \\
& \left(Y, \eta^{-1}  \xi_0\right)\in \overline{\mathcal J}^{2,-}m^2(y_0),
\end{aligned}
\right.
\end{equation*}
and
\begin{equation} \label{e.matfunwgamma}
-\left(\frac1\gamma  +  |M| \right) I_{2d} \leq  \begin{pmatrix} 
X+D^2\psi(x_0,t_0) &0 \\  0 & -\eta Y 
\end{pmatrix} 
\leq M + \gamma M^2,
\end{equation}
where~$\xi_0:=\delta^{-1} |x_0-y_0|^2(x_0-y_0) $ and
\begin{equation*} \label{}
M:= \frac1\delta \begin{pmatrix} 
N & -N  \\ -N & N
\end{pmatrix}, \quad N:=  |x_0-y_0|^2  I_d + 2(x_0-y_0) \otimes (x_0-y_0).
\end{equation*}
If $x_0=y_0$, then we take $\gamma:=1$. Otherwise, if $x_0 \neq y_0$ (we will show below that this is indeed the case), we set $\gamma:= \delta |x_0-y_0|^{-2}$. In the latter case, we obtain 
\begin{equation} \label{e.matsveryfun}
\frac{C}{\delta} \left| x_0 - y_0 \right|^2 I_{2d}
 \leq 
 \begin{pmatrix} 
X+D^2\psi(x_0,t_0) &0 \\  0 & -\eta Y 
\end{pmatrix} \leq
\frac{C}{\delta} \left| x_0 - y_0 \right|^2 \begin{pmatrix} 
I_d & -I_d  \\ -I_d & I_d
\end{pmatrix}.
\end{equation}
In particular, we have 
\begin{equation*} \label{}
\left| X \right| + \left| Y \right|  \leq  C\delta^{-1} \left| x_0 - y_0 \right|^2+C\ep \leq C K^{\frac23} \delta^{-\frac13} \leq C\delta^{-\frac13}.
\end{equation*}

\smallskip

\emph{Step 3.} We prove a lower bound for $|x_0 - y_0|$. The claim is that
\begin{equation} \label{e.x0y0lower}
\left| x_0 - y_0 \right|  \geq c (\mu \delta)^{\frac12}.
\end{equation}
In particular, this implies that we are in the case~$x_0 \neq y_0$. First we note that~\eqref{e.psiDbounds} and~\eqref{e.matsveryfun} give a lower bound for $Y$ which, in both cases $x_0=y_0$ (when we have $M=0$) and $x_0\neq y_0$ can be written as
\begin{equation*} \label{}
Y \geq  - C \eta^{-1} \delta^{-1} |x_0-y_0|^2I_d \geq - C\delta^{-1} |x_0-y_0|^2I_d .
\end{equation*}
Next we observe that the equation for $m^2$ yields
\begin{equation*} \label{}
-\tr_*\left( A\left(\eta^{-1}  \xi_0 ),y_0 \right) Y \right) + H\left(\eta^{-1}  \xi_0, y_0 \right) \geq \mu.
\end{equation*}
Combining the previous two inequalities and using the positive homogeneity and coercivity of $H$ and~\eqref{e.psiDbounds} gives
\begin{equation*} \label{}
C\delta^{-1} |x_0-y_0|^2+ C\delta^{-p}\left| x_0-y_0\right|^{3p}  \geq \mu.
\end{equation*}
This yields~\eqref{e.x0y0lower}. 

\smallskip

\emph{Step 4.} We complete the proof of~\eqref{e.trifecta} by deriving a contradiction, using the results of the previous two steps and the equations respectively satisfied by $w_1$ and $m^2$. Using that $w_1$ is a subsolution of~\eqref{e.evolMCM} and $\partial_t\psi\geq 0$, we have
\begin{multline*} \label{}
- \tr^*\left( A(\xi_0+D\psi(x_0,t_0)),x_0) X \right) + H\left( \xi_0+D\psi(x_0,t_0) , x_0 \right) \\
\leq \partial_t \psi(x_0,t_0) - \tr\left( A(\xi_0+D\psi(x_0,t_0)),x_0) X \right)  + H\left( \xi_0 +D\psi(x_0,t_0), x_0 \right) \leq \mu. 
\end{multline*}
By~\eqref{e.psiDbounds} and~\eqref{e.x0y0lower}, we have
\begin{multline*}
\left| \tr^*(A\left(\xi_0+D\psi(x_0,t_0) , x_0 \right)X)   - \tr(A\left(\xi_0, x_0 \right)X) \right| \\
 \leq  C \left| \xi_0 \right|^{-1} \left| D\psi(x_0,t_0)\right|  |X|
 \leq C\delta \left| x_0-y_0 \right|^{-3} \ep \delta^{-\frac13}
 \leq C\mu^{-\frac32}\delta^{-\frac56} \ep
\end{multline*}
and, by~\eqref{e.psiDbounds} and~\eqref{e.boundx-y},
\begin{equation*} \label{}
\left| H\left( \xi_0+D\psi(x_0,t_0), x_0 \right) - H\left( \xi_0, x_0 \right)  \right| \leq C\left| \xi_0 \right|^{p-1} \left| D\psi(y_0,t_0) \right|\leq C  \ep.
\end{equation*}
Combining the previous three lines and using $\left| D^2\psi(x_0,t_0)\right| \leq C\ep$, we get 
\begin{equation} \label{e.ineq1}
- \tr\left( A(\xi_0,x_0) (X+D^2\psi(x_0,t_0)) \right) + H\left( \xi_0, x_0 \right) \\
\leq \mu + C\mu^{-\frac32}\delta^{-\frac56}\ep.
\end{equation}
Using the equation for $m^2$, the homogeneity of $A$ and $H$ and the fact that $\eta \geq 1$ and $H\geq 0$, we obtain
\begin{equation} \label{e.ineq2}
-\tr\left( A\left(\xi_0, y_0 \right) \eta Y \right) + H\left( \xi_0, y_0 \right) \geq \eta \mu.
\end{equation}

The goal is to show that~\eqref{e.ineq1} and~\eqref{e.ineq2} are incompatible. 
Using the matrix inequality~\eqref{e.matsveryfun}, we find
\begin{multline*}
\tr\left( A(\xi_0,x_0)(X+D^2\psi(x_0,t_0))  - A(\xi_0,y_0)\eta Y\right) \\
\leq C\delta^{-1} \left| x_0-y_0\right|^2 \left| \sigma(\eta,x_0) - \sigma(\eta,y_0) \right|^2  \leq C\delta^{-1} \left| x_0-y_0\right|^4   \leq C\delta^{\frac13}.
\end{multline*}
We note that 
\begin{equation*} \label{}
\left| H\left( \xi_0 , x_0 \right) - H\left( \xi_0 , y_0 \right) \right| \leq C \left| \xi_0 \right|^{p-1} \left| x_0 - y_0 \right| \leq C\delta^{\frac13}.
\end{equation*}
Taking the difference of~\eqref{e.ineq1} and~\eqref{e.ineq2} and using the previous two inequalities, we get
\begin{equation*} \label{}
(\eta-1) \mu \leq C\left( \mu^{-\frac32}\delta^{-\frac56}\ep+  \delta^{\frac13} \right).
\end{equation*}
If we choose $\eta :=1+ C\mu^{-1}\left( \mu^{-\frac32}\delta^{-\frac56}\ep+  \delta^{\frac13} \right)$, for a sufficiently large constant $C$, we obtain the desired contradiction. This completes the proof of~\eqref{e.trifecta}.

\smallskip

\emph{Step 5.} We deduce that, under certain restrictions on $R$ and $\delta$, that 
\begin{equation} \label{e.w1w2}
w_1(x,t)\leq \eta m^2 (x,t)+ \frac14 \qquad \forall (x,t)\in \H\times [0,T].
\end{equation}
We consider the three alternatives provided by~\eqref{e.trifecta}.
In the case $t_0=0$, we get
\begin{equation*} \label{}
\sup_{(x,t) \in \H^+\times [0,T]} \left( w_1 - \eta m^2  \right) \leq \max \Psi = w_1(x_0,0)-\eta m^2 (y_0,0)-\frac{|x_0-y_0|^4}{4\delta}-\psi(x_0,0)\leq 0
\end{equation*}
because $w_1(\cdot,0)< 0$ by construction and $0\leq m^2(\cdot,0)$ by assumption. Therefore we consider the case that~$x_0$ or~$y_0$ belongs to~$\partial \H^+$. We give only the argument in the case that $x_0\in\partial \H^+$, the other case being analogous. We divide this case into two sub cases: $x_0\in B_R$ or $x_0 \not\in B_R$. If $x_0\in B_R$, then the assumption~\eqref{e.miord} yields that 
\begin{align*}
\max \Psi & =   w_1(x_0,t_0)-\eta m^2 (y_0)-\frac{\left|x_0-y_0\right|^4}{4\delta}-\psi(x_0,t_0) \\
& \leq m^1(x_0)- \eta m^2(y_0) 
 \leq m^2(x_0)- m^2(y_0) \leq K \left|x_0-y_0 \right| \leq C \delta^{\frac13}\leq \frac14,
\end{align*}
provided that $\delta$ is sufficiently small. If, on the contrary, $x_0\in \partial \H^+\setminus B_R$, then
\begin{equation*} \label{}
\max \Psi \leq \nu T- \eta m^2(y_0) -\psi(x_0,t_0)\leq \nu T - g(R - T).
\end{equation*}
If $R \geq T $, then, since $g(s)= \ep s+1/4$ for $s\geq 0$,
\begin{equation*}
\max \Psi \leq \nu T -\ep(R-T). 
\end{equation*}
Therefore, if 
\begin{equation} \label{e.condcond1}
R\geq T( 1+ \ep^{-1} \nu),
\end{equation}
then we obtain $\max\Psi \leq \frac14$. This completes the proof of~\eqref{e.w1w2}. 

\smallskip

\emph{Step 6.} We complete the argument by deriving~\eqref{e.FSconc}. Fix $s\geq 1$ and choose $\lambda := c\mu/L_\mu^2$ and $\nu := \mu/2$. Then  
\begin{equation*}
w_1(x,t) \geq m^1(x) - \frac {CL^2}{\mu} \exp\left( -c\mu \left(\frac{\mu t}2  -m^1(x) \right ) \right).
\end{equation*}
By assumption~\eqref{e.migrowth}, we have $m^1(se)\leq {M}+K s$, and thus we can choose $t\geq Cs+2\mu^{-1} {M}$ so that the last term in right-hand side is at most $\frac14$.
On the other hand, if $1+s -  \left(T-t \right)\leq 0$, then, since $g\leq 1/4$ on $(-\infty, 0]$,  we have
\begin{align*}
\eta m^2(se)+\psi(se,t) & \leq  \eta m^2(se)+g\left( \left( 1+s^2 \right)^{\frac12} -  \left(T-t \right)  \right) \\
& \leq  m^2(se)+ (\eta -1)({M}+K s) + \frac14,
\end{align*}
where we used assumption~\eqref{e.migrowth} again in the last line. Note that we can choose~$t$ such that both conditions $t\geq C_\mu s+2\mu^{-1} {M}$ and $1+s -  \left(T-t \right)\leq 0$ hold, provided
\begin{equation}\label{e.condcond2}
T\geq 2\mu^{-1} {M}+C  s+ 1.
\end{equation}
In this case, we get, by the choice of $\eta$, 
$$
m^1(se)\leq m^2(se) +C\mu^{-1}\left( \mu^{-\frac32}\delta^{-\frac56}\ep+  \delta^{\frac13} \right)({M}+Ks) + \frac34.
$$
We conclude by selecting $\delta := \mu^{\frac97}\ep^{\frac67}$ and $\ep:= c \mu^5({M}+Ks)^{-\frac72}$. Then, if we choose $T := C\mu^{-1}({M}+s)$ and $R = C\mu^{-5} ({M}+s)^{\frac92}$, so that \eqref{e.condcond1} and \eqref{e.condcond2} hold, we obtain $m^1(se)\leq m^2(se)+1$ as claimed. This completes the proof.
\end{proof}

\subsection{Convergence of the means}
\label{s.ssbias}

The idea of the proof of Proposition~\ref{p.statbiasBIS} is to use the estimate of the stochastic fluctuations and the ``approximate" finite speed of propagation to compare~$m_\mu(\cdot,\H_e^-)$ to $m_\mu(\cdot,\H_e^-+se) + \E \left[ m_\mu(se,\H_e^-) \right]$. If these functions are close, then the map $t\to \E[m_\mu(te,\H_e^-)]$ is almost linear.

\smallskip

We begin by extending the fluctuation estimate to large balls, using union bounds. For $t>0$ and $R> 1$, we define
\begin{equation}\label{e.defN+N-}
N^+_R(t):= \sup_{x \in B_{R t} \cap (\partial \H_e^++te) } m_\mu(x,\H_e^-), \quad N^-_R(t):= \inf_{x \in B_{R t} \cap (\partial \H_e^++te) } m_\mu(x,\H_e^-).
\end{equation}

\begin{lemma}\label{l.N+N-}
There exists $C(\data,L)\geq 1$ such that, for every $R,t>1$,
\begin{multline}
\E \left[ \left| N^+_R(t) - \E \left[ m_\mu(te,\H_e^-) \right] \right| \right]  + \E\left[  \left| N^-_R(t) - \E \left[ m_\mu(te,\H_e^-) \right] \right| \right]  \\
\leq C( \mu l_\mu^2)^{-\frac12} t^{\frac12} \log^{\frac12} (1+Rt ). 
\end{multline} 
\end{lemma}
\begin{proof}
We prove the estimate for $N^+_R(t)$, the one for $N^-_R(t)$ being obtained in a similar way. By the Lipschitz estimate and a union bound, 
\begin{align*}
\lefteqn{ \P[ |N^+_R(t)-\E[m_\mu(te,\H_e^-)]\geq \lambda] } 
\qquad & 
\\
& \leq \P\left[ \sup \left\{  \left| m_\mu(x+te,\H_e^-)-\E[m_\mu(te,\H_e^- )] \right| \,:\, x\in \partial \H_e^+ \cap B_{R t}\cap c\lambda\Z^d \right\}  \geq \frac \lambda 2 \right] \\
& \leq \sum_{ x\in \partial \H_e^+ \cap B_{R t}\cap c\lambda\Z^d}  \P\left[ \left|m_\mu(x+te,\H_e^-)-\E[m_\mu(te, \H_e^-)]\right| \geq \frac \lambda2 \right]. 
\end{align*}
There are $C(Rt)^{d-1}\lambda^{1-d}$ terms in the sum in the previous line. 
Therefore, by stationarity and Corollary~\ref{c.fluctuationsH},
\begin{equation*} \label{}
\P\left[ \left| N^+_R(t)-\E\left[m_\mu(te,\H_e^-)\right] \right| \geq \lambda\right] \leq C(Rt)^{d-1}\lambda^{1-d} \exp\left( -\frac{c\mu l_\mu^2\lambda^2}{t} \right).
\end{equation*}
Integrating this inequality, we get, for every $\alpha\geq 1$,
\begin{align*}
\lefteqn{ 
\E\left[ \left|N^+_R(t)-\E\left[m_\mu(te,\H_e^-)\right]\right|\right] }
\qquad & \\
& = \int_0^\infty \P\left[\left|N^+_R(t)-\E\left[m_\mu(te, \H_e^-)\right] \right| \geq \lambda\right]\, d\lambda  \notag \\
& \leq  \alpha \left( t\log(1+Rt) \right)^{\frac12} + C(Rt)^{d-1} \int_{ \alpha \left( t\log(1+Rt) \right)^{\frac12}}^\infty \lambda^{1-d} \exp\left( -\frac{c\mu l_\mu^2\lambda^2}{t} \right)\, d\lambda. \notag
\end{align*}
Choosing $\alpha=( \mu l_\mu^2)^{-\frac12}$, we find, after a change of variable in the integral: 
\begin{align*}
\lefteqn{ 
\E\left[ \left|N^+_R(t)-\E\left[m_\mu(te,\H_e^-)\right]\right|\right] }
\qquad & \\
& \leq  \alpha \left( t\log(1+Rt) \right)^{\frac12} + CR^{d-1}t^{\frac{d}{2}} (\mu l_\mu)^{\frac{d}{2}-1}\int_{ \left( \log(1+Rt) \right)^{\frac12}}^\infty s^{1-d} \exp\left( -s^2 \right)\, ds  \notag \\
& \leq   ( \mu l_\mu^2)^{-\frac12} \left( t\log(1+Rt) \right)^{\frac12} + C  \notag 
\end{align*}
This yields the lemma.
\end{proof}

Using the previous lemma, the approximate finite speed of propagation and a very simple comparison argument, we next show that the quantity~$\E\left[ m_\mu(te) \right]$ is almost linear in~$t$.
 
\begin{lemma}\label{l.almostlinearBIS}
There exists $C(\data,L)\geq1$ such that, for every $s,t> 1$,
\begin{multline}
\left| \E\left[ m_\mu((t+s)e,\H_e^-) \right] - \E\left[ m_\mu(te,\H_e^-) \right] - \E\left[ m_\mu(se,\H_e^-) \right] \right| \\
 \leq C_\mu \left( s+t \right)^{\frac12} \log^{\frac12} \left( 1+s+t \right),
\end{multline}
where $C_\mu := C \mu^{-\frac32} (1+\left|\log\mu\right|)^{-\frac12}$. 
\end{lemma}

\begin{proof}
We may assume that $s\leq t$. Let $N^+(t)$ and $N^-(t)$ be defined by~\eqref{e.defN+N-}, where we take $R:=C\mu^{-5} (1+L_\mu t +s)^{\frac92}$, with~$C$ the large constant appearing in Proposition~\ref{p.finitespeedMCM}. We apply Proposition~\ref{p.finitespeedMCM} to $m^1(x):= m_\mu(x)$ and $m^2(x)= m_\mu(x, \H_e^-+te)+N^+(t)$ in the domain~$\H_e^++te$. To check the hypotheses of the proposition, we note that both $m^1$ and $m^2$ are solutions of 
\begin{equation}\label{e.dansHte}
-\tr\left( A(x) D^2m \right) + H\left( Dm,x \right) = \mu \quad \mbox{in} \ \H_e^++te.
\end{equation}
By the Lipschitz estimates, we have $\|Dm^1\|_\infty,\|Dm^2\|_\infty\leq L_\mu$. In particular, $0\leq m^i(x)\leq {M}+L_\mu |x-te|$ for every $x\in\H_e^-+te$ and $i\in\{1,2\}$, with ${M}:=L_\mu t$. Finally, we observe that $m^1\leq  m^2$ on $(\partial \H_e^-+te)\cap B_R(te)$ by the definition of~$N^+(t)$. Therefore the conclusion of Proposition~\ref{p.finitespeedMCM} yields
\begin{equation*} \label{}
m_\mu((t+s)e,\H_e^-)\leq m_\mu((t+s)e, \H_e^-+te)+N^+(t) +1.
\end{equation*}
Taking expectations and applying Lemma~\ref{l.N+N-}, we get
\begin{align*}
\lefteqn{
\E\left[m_\mu((t+s)e,\H_e^-)\right] 
}
\quad & \\
& \leq \E\left[m_\mu((t+s)e, \H_e^-+te)\right]+\E\left[m_\mu(te,\H_e^-)\right]+ C (\mu l_\mu^2)^{-\frac12}t^{\frac12} \log^{\frac12} (1+Rt)+1\\
& = \E\left[m_\mu(se,\H_e^-)\right]+ \E\left[m_\mu(te,\H_e^-)\right]+ C (\mu l_\mu^2)^{-\frac12} \log^{\frac12} (1+Rt).
\end{align*}
Since $R \leq C\mu^{-5}(1+t)^{\frac92}$ and $\ell_\mu\geq c\mu$,  we obtain
\begin{align*} \label{}
\lefteqn{ \E\left[m_\mu((t+s)e,\H_e^-)\right]  } 
\qquad & \\
& \leq \E\left[m_\mu(te,\H_e^-)\right]+ \E\left[m_\mu(se,\H_e^-)\right] + C(\mu l_\mu^2|\log(\mu)|)^{-\frac12} t^{\frac12} \log^{\frac12} (1+t) \\
& \leq \E\left[m_\mu(te,\H_e^-)\right]+ \E\left[m_\mu(se,\H_e^-)\right]  + C\mu^{-\frac32} \left|\log \mu\right|^{-\frac12} t^{\frac12} \log^{\frac12} (1+t).
\end{align*}
To obtain the reverse inequality, we take $m^1(x)= m_\mu(x, \H_e^-+te)+N^-(t)$ and $m^2(x)= m_\mu(x,\H_e^-)$ and follow similar reasoning.
This completes the proof. 
\end{proof}

We next use Lemma~\ref{l.almostlinearBIS} to show that the expectation of the solution of the planar metric problem is approximately affine far from the boundary plane, thereby obtaining the first statement of Proposition~\ref{p.statbiasBIS}. 

\begin{lemma}
\label{l.cauchymeans}
There exist $\overline m_\mu(e) >0$ and $C(\data,L)\geq 1$ such that, for every $t>1$,
\begin{equation} \label{e.cauchymeans}
\left| \frac1t \E \left[ m_\mu(te,\H_e^-)\right] - \overline m_\mu(e) \right| \leq C_\mu t^{-\frac12} \log^{\frac12}(1+t),
\end{equation}
where $C_\mu := C \mu^{-\frac32}(1+\left|\log\mu\right|)^{-\frac12}$. 
\end{lemma}
\begin{proof} 
As the dependence with respect to $\mu$ of the various constants $C$ come from Lemma~\ref{l.almostlinearBIS}, we omit this dependence throughout the proof. For simplicity, we denote $m_\mu(x) = m_\mu(x,\H_e^-)$. We break the argument into three steps.

\smallskip

\emph{Step 1.} The application of Lemma~\ref{l.almostlinearBIS}. It is convenient to denote
\begin{equation*} \label{}
G(t):= \frac1t \E \left[ m_\mu(te,\H_e^-)\right].
\end{equation*}
Then Lemma~\ref{l.almostlinearBIS} gives, for every $1 < t < s$,
\begin{align} \label{e.almostlinapp}
\left| G(s) - G(t) \right| 
& = \frac1s \left| \E \left[ m_\mu(se) \right]  - \E \left[ m_\mu(te) \right] - \left( \frac{s}{t}-1 \right)  \E \left[ m_\mu(te) \right]  \right| \\
& \leq \frac1s \left| \E \left[ m_\mu(se) \right]  - \E \left[ m_\mu(te) \right] - \E \left[ m_\mu((s-t)e) \right] \right| \notag \\
& \qquad + \left( 1-\frac ts \right) \left| \frac1t \E \left[ m_\mu(te) \right]  - \frac1{s-t}\E\left[ m_\mu((s-t)e) \right] \right| \notag \\
& \leq C s^{-\frac12} \log^{\frac12}(1+s) + \frac{s-t}{s} \left| G(t) - G(s-t) \right|. \notag
\end{align}

\smallskip

\emph{Step 2.} We claim there exists $C\geq 1$ such that, for every~$1<t< s$,
\begin{equation} \label{e.twostep}
\left| G(s) - G(t) \right| \leq Ct^{-\frac12} \log^{\frac12} (1+t).
\end{equation}
The argument is by induction. Clearly the statement holds for a fixed constant $C_1$ and $t\in (1,2]$ by the Lipschitz estimate. Suppose that, for a fixed $m\in\N$, the statement holds whenever $t \in (1, 2^m]$ and with a constant $C=C_m\geq1$. Then for every $t\in \left(1,2^{m+1}\right]$ and $s\in \left(t,\frac32 t\right]$, we have that $s-t \leq 2^m$ and $s^{-1} (s-t)^{\frac12} \leq 2^{-\frac12} t^{\frac12}$ and therefore 
\begin{equation*}
\frac{s-t}{s} \left| G(t) - G(s-t) \right| \leq C_ms^{-1} (s-t)^{\frac12}  \log^{\frac12}(1+(s-t)) \leq 2^{-\frac12} C_m t^{\frac12} \log^{\frac12} \left( 1+t\right).
\end{equation*}
Combining with~\eqref{e.almostlinapp}, we deduce that, for every $t\in (1,2^{m+1}]$ and $s\in \left(t,\frac32t\right]$,
\begin{equation*} \label{}
\left| G(s) - G(t) \right| \leq \left( C + 2^{-\frac12} C_m \right) t^{-\frac12} \log^{\frac12}(1+t).
\end{equation*}
To remove the restriction on $s$, fix $t\in \left(1,2^{m+1}\right]$ and  $s\in (t,\infty)$. Denote $s_k:= 2^{-k} s$ for $k\in\N$. Take $n \in\N$ to be the unique positive integer such that $s_n \in \left(\frac34 t,\frac32 t\right]$. According to~\eqref{e.almostlinapp}, we get
\begin{align*}
\left| G(s_n) - G(s) \right| 
& \leq  \sum_{k=1}^n \left| G(s_{k})-G(s_{k-1}) \right|  \\
& \leq  \sum_{k=1}^n Cs_k^{-\frac12} \log^{\frac12}(1+s_k) \\
& \leq Ct^{-\frac12} \log^{\frac12} \left( 1+t \right).
\end{align*}
Since either $s_n\leq 2^{m+1}$ and $t\in \left[s_n,\frac32 s_n\right)$ or else $s_n \in \left(t,\frac32 t\right]$, we have that
\begin{equation*} \label{}
\left| G(s_n) - G(t) \right| \leq \left( C + 2^{-\frac12} C_m \right) t^{-\frac12} \log^{\frac12}(1+t).
\end{equation*}
Combining the previous two inequalities yields, for every $t \in \left(1,2^{m+1} \right]$ and $s\in (t,\infty)$, 
\begin{equation*} \label{}
\left| G(s) - G(t) \right| \leq \left( C + 2^{-\frac12} C_m \right) t^{-\frac12} \log^{\frac12}(1+t).
\end{equation*}
Thus we have shown that~\eqref{e.twostep} holds for every $t \in \left(1,2^{m+1} \right]$ and $s\in (t,\infty)$ with a constant $C_{m+1} = \left(C + (1-c) C_m\right)$, for some $c>0$. It is clear that the sequence of constants $\{ C_m \}$ remains bounded as $m\to \infty$. By induction, we therefore obtain~\eqref{e.twostep} for every $1<t<s$ and a fixed constant $C\geq1$.

\smallskip

\emph{Step 3.} We conclude. It is immediate from Step~2 that the sequence $\{ G(2^m) \}_{m\in\N}$ is Cauchy and therefore has a limit, which we denote by $\overline m_\mu(e)$. Taking $s=2^m$ in~\eqref{e.twostep} and sending $m\to \infty$ yields, for every $t>1$, 
\begin{equation*} \label{}
\left| G(t) - \overline m_\mu(e) \right| \leq C t^{-\frac12} \log^{\frac12} \left( 1+ t\right).
\end{equation*}
This complete the proof of~\eqref{e.cauchymeans} and the lemma.
\end{proof}

The previous lemma gives the first statement of Proposition~\ref{p.statbiasBIS}. To obtain the second statement of the proposition, it remains to check that $(\mu,e) \mapsto \overline m_\mu(e)$ is continuous, which relies on another application of Proposition~\ref{p.finitespeedMCM}.

\begin{lemma}
\label{l.mbarcont}
For any $L\geq 1$, there exist $C(\data,L)\ge1$ and $0<c(\data,L)\leq 1$ such that, 
for each $0 < \nu\leq \mu \leq L$ and $e_1,e_2\in \partial B_1$,
\begin{equation*}
\left| \overline m_\nu(e_1)- \overline m_\mu (e_2) \right|  \leq C\mu^{\frac1p-1} |\mu-\nu|+ C_\mu  \left|e_1-e_2\right|^{\frac18}\left( 1+ \left| \log|e_1-e_2| \right| \right)^{\frac12}
\end{equation*}
where $C_\mu= C\mu^{-\frac{31}{16}}(1+ \left|\log \mu \right|)^{\frac12}$ and 
\begin{equation*} \label{}
\bar m_\mu(e)\geq \bar m_\nu(e)+c\nu^{\frac1p-1}(\mu-\nu) .
\end{equation*}
Finally, for every $\mu>0$,
\begin{equation} \label{e.mbargrowth}
\left( \frac\mu{c_0} \right)^{\frac1p} \leq  \bar  m_\mu (e) \leq \left( \frac\mu{C_0} \right)^{\frac1p}.
\end{equation}
\end{lemma}
\begin{proof}
We use Proposition~\ref{p.finitespeedMCM} to approximate $\H_e^-$ by a compact target which is continuous in $e$. We then use the fact that $K\mapsto m_\mu(x,K)$ is Lipschitz with respect to the Hausdorff metric.

\smallskip

We fix $e_1,e_2\in\partial B_1$. Given $R\geq 2$, there exists a smooth convex set $K^i_R$, $i\in\{1,2\}$, with diameter at most $CR$ and curvature bounded by~$1$ which is contained in $\H^{-}_{e_i}$ and is such that $B_R \cap \partial \H^-_{e_i} \subseteq K_R \cap \partial \H^-_{e_i}$. Furthermore, we may assume that $K_R^2$ is the image of $K_R^1$ under the rotation that sends $e_1$ to $e_2$. Since the diameter of $\mathcal K_R^i$ is at most $CR$, we have
\begin{equation} \label{e.KRrot}
\dist_H\left( K_R^1,K_R^2 \right) \leq CR\left| e_1-e_2 \right|. 
\end{equation}

\smallskip

\emph{Step 1.} We show that, for every $1< s< (C^{-1}\mu^5 R)^{\frac29}-1$ and $i\in\{1,2\}$,
\begin{equation}\label{e.mHmKR}
\left| m_\mu(se_i, \H_{e_i}^-) -  m_\mu(se_i, K^i_R)\right|  \leq 1,
\end{equation}
where $C$ is the constant in Proposition \ref{p.finitespeedMCM}. Let $m^1:= m_\mu(\cdot, K_R^i)$ and $m^2:=m_\mu(\cdot,\H_{e_i}^-)$. Since $K_R^i\subset \mathcal H^-_{e_i}$, we have $m^2\leq m^1$. An application of Proposition~\ref{p.finitespeedMCM} then gives~\eqref{e.mHmKR}.

\smallskip

\emph{Step 2.} According to \eqref{e.Lipschpc}, we have for any $x\notin  K_{R}^1 \cup K_{R}^2$,
\begin{equation*} \label{}
\left| m_\mu(x, K_{R}^1)  -  m_\mu(x, K_{R}^2)\right|\leq L_\mu R |e_2-e_1|. 
\end{equation*}
By the Lipschitz estimate for $m_\mu(\cdot, K_R^i)$, we obtain, provided that $|e_2-e_1| <c$, 
\begin{align*}
\lefteqn{ \left| m_\mu(se_1, \H_{e_1}^-) -  m_\mu(se_2, \H_{e_2}^-) \right| 
}
\qquad & \\
 & \leq \left| m_\mu(se_1, \H_{e_1}^-) -  m_\mu(se_1, K_R^1) \right| + \left| m_\mu(se_1, K_R^1) - m_\mu(se_2, K_R^1)  \right|  \\
 & \qquad +  \left| m_\mu(se_2, K_R^1) -  m_\mu(se_2, K_R^2) \right| + \left| m_\mu(se_2, K_R^2) - m_\mu(se_2, \H_{e_2}^-)  \right|  \\
 & \leq 2 + C \left( R + s \right) \left| e_1 - e_2 \right|.
\end{align*}
Dividing by $s$, taking expectations and using~\eqref{e.cauchymeans}, we obtain 
\begin{equation*} \label{}
\left| \overline m_\mu(e_1) - \overline m_\mu(e_2) \right|  \leq 2s^{-1} + C\left( Rs^{-1} +1  \right) \left| e_1 - e_2 \right| + C_\mu s^{-\frac12} \log^{\frac12} (1+s),
\end{equation*}
where $C_\mu= C\mu^{-\frac32}(1+|\log(\mu)|)^{\frac12}$.
Taking
\begin{equation*}
R:= C\mu^{-5}\left(1+ s^{\frac92}\right) 
\quad \mbox{and} \quad
 s:=\left(\mu^{-\frac72} |e_1-e_2|\right)^{-\frac14}
\end{equation*}
we obtain, for every $e_1,e_2\in \partial B_1$ with $|e_1-e_2|\leq c$,
$$
|\bar m_\mu(e_1)-\bar m_\mu(e_2)| \leq C\mu^{-\frac{31}{16}}(1+|\log(\mu)|)^{\frac12}|e_1-e_2|^{\frac18}\left| \log |e_1-e_2| \right|^{\frac12}.
$$

\smallskip

\emph{Step 4.} We next check the continuity of $\mu\to \bar m_\mu(e)$. Fix $0<\nu\leq\mu<\infty$. It is clear that $m_\nu\leq m_\mu$ and thus $\bar m_\nu\leq \bar m_\mu$. By positive homogeneity of $H$, $(\mu/\nu) m_\nu$ is a supersolution of the planar metric problem for $\mu$. By comparison, we get $m_\mu\leq (\mu/\nu)m_\nu$, and thus  $\bar m_\mu\leq (\mu/\nu)\bar m_\nu$. On the other hand, we have $\overline m_\nu \leq C_0^{-1} \nu^{\frac1p}$. Combining these yields
\begin{equation*} \label{}
\left| \overline m_\mu - \overline m_\nu \right| \leq C\nu^{\frac1p-1} \left| \nu - \mu \right|.
\end{equation*}
  
\smallskip

\emph{Step 5.} We show that $\mu \to \bar m_\mu(e)$ is strictly  increasing. For $\nu<\mu$, let $w:= (1+\ep) m_\nu$. Then, by homogeneity of $H$,  
\begin{align*}
\lefteqn{ -\tr (A(Dw,x)D^2e)+H(Dw,x) }
\qquad & \\
&  = (1+\ep) \left(-\tr (A(Dm_\nu,x)D^2m_\nu) +H(Dm_\nu,x)\right) \\
& \qquad + (1+\ep)((1+\ep)^{p-1}-1)H(Dm_\nu, x) \\
& \leq (1+\ep)\nu +C\ep L_\nu \\
& \leq  \mu
\end{align*}
provided that $\ep \leq c (\mu-\nu)\wedge 1$. We deduce that $m_\mu\geq m_\nu+ (c (\mu-\nu)\wedge 1) m_\nu$, which implies that $\bar m_\mu(e)\geq \left( 1+ (c (\mu-\nu)\wedge 1)\right) \bar m_\nu(e)$, as desired. 

\smallskip

\emph{Step 6.} We conclude by noticing that~\eqref{e.mbargrowth} follows from~\eqref{RoughGrowthmmue}.
\end{proof}

The proof of Proposition~\ref{p.statbiasBIS} is now complete, as the statement follows immediately from Lemmas~\ref{l.cauchymeans} and~\ref{l.mbarcont}. 

\smallskip

We now present the proof of Proposition~\ref{p.planarquantitative}.

\begin{proof}[{Proof of Proposition~\ref{p.planarquantitative}}]

Fix $L\geq 1$, $e \in \partial B_1$ and $\mu\in (0,L]$. An immediate consequence of Corollary~\ref{c.fluctuationsH} and Proposition~\ref{p.statbiasBIS} is the following estimate, which is valid for every $\lambda>0$ satisfying 
\begin{equation} \label{e.lambdaineq}
\lambda \geq C\mu^{-\frac32} \left( 1+x\cdot e\right)^{\frac12} \left( 1+ \left| \log \mu \right| + \log \left( 1+x\cdot e \right) \right)^{\frac12}
\end{equation}
and $x\in \H_e^+$: 
\begin{equation} \label{e.muannoy}
\P \left[ \left| m_\mu(x,\H_e^-) - \overline m_\mu (e)(x\cdot e) \right| > \lambda \right]  \leq C \exp\left( \frac{-c\mu^3\lambda^2}{1+x\cdot e} \right).
\end{equation}
On the other hand, we have the following deterministic bounds from~\eqref{RoughGrowthmmue} and~\eqref{e.mbargrowth}:
\begin{equation*} \label{}
\left| m_\mu \left( x,\H_e^- \right) - \overline m_\mu(e) (x\cdot e) \right| \leq 2\left( \frac\mu{c_0} \right)^\frac1p (x\cdot e) \leq C  (x\cdot e) \mu.
\end{equation*}
We deduce that, if
\begin{equation} \label{e.muineq}
\mu\leq c\left(1+ x\cdot e\right)^{-\frac15} \log^{\frac15} (1+x\cdot e)
\end{equation}
and
\begin{equation} \label{e.lambdaineq2}
\lambda \geq C \left( 1+x\cdot e \right)^{\frac45} \left( 1+ \log \left( 1+x\cdot e \right) \right)^{\frac15}, 
\end{equation}
then
\begin{equation*} \label{}
\left| m_\mu \left( x,\H_e^- \right) - \overline m_\mu(e) (x\cdot e) \right| \leq \lambda \quad \mbox{$\P$--a.s.}
\end{equation*}
Notice that~\eqref{e.lambdaineq2} implies that at least one of~\eqref{e.muineq} and~\eqref{e.lambdaineq} must hold. We deduce that, for every $\mu>0$ and every $\lambda$ satisfying~\eqref{e.lambdaineq2}, we have
\begin{equation} \label{e.muannoy2}
\P \left[ \left| m_\mu(x,\H_e^-) - \overline m_\mu (e)(x\cdot e) \right| > \lambda \right]  \leq C \exp\left( \frac{-c\lambda^2 \log^{\frac35} (1+x\cdot e)}{\left(1+x\cdot e\right)^{\frac85}} \right).
\end{equation}
By adjusting the constant $C$, we obtain, for every $\lambda>0$,
\begin{equation} \label{e.muannoy3}
\P \left[ \left| m_\mu(x,\H_e^-) - \overline m_\mu (e)(x\cdot e) \right| > \lambda \right]  \leq C\exp\left( \frac{-\lambda^2}{(1+x\cdot e)^{\frac85}} \right).
\end{equation}
Indeed, the right side of~\eqref{e.muannoy3} larger than that of~\eqref{e.muannoy2} and also larger than~$1$ if~$\lambda$ does not satisfy~\eqref{e.lambdaineq2}. We have proved~\eqref{e.Pm-barm}. The fact that $(\mu,e) \mapsto \overline m_\mu(e)$ is continuous was already proved in Lemma~\ref{l.mbarcont}. 
\end{proof}

We conclude this section with a crude deterministic estimate concerning the continuity of $m_\mu(x,\H^--se)$ with respect to its parameters $\mu$ and $e$, which is obtained from a variation of the proof of Lemma~\ref{l.mbarcont}, above. This is needed in Section~\ref{s.outline} in the proof of the main result to ``snap to a grid" before taking union bounds in the final step of the argument.

\begin{lemma}
\label{l.mconte} There exists~$C(\data,L)\geq1$ such that, for every $0 < \nu\leq \mu \leq L$, $s\geq 1$, $e_1,e_2\in \partial B_1$ and $x\in (\H_{e_1}^+-se_1)\cap (\H_{e_2}^+-se_2)$,
\begin{multline}
\label{e.mconte}
\left|m_\mu(x, \H_{e_1}^--se_1)-m_\nu(x, \H_{e_2}^--se_2)\right| 
\\
\leq 2 +C\left( |\mu-\nu|^{\frac12}|x|+ 
(1+|x|+s)^{\frac92} |e_2-e_1|^{\frac{1}{5p+1}}\right).
\end{multline}
\end{lemma}
\begin{proof}
We first argue that 
\begin{multline}\label{RoughEstimmu}
\left|m_\mu(x, \H_{e_1}^--se_1)-m_\nu(x, \H_{e_2}^--se_2)\right| 
\\
\leq C\nu^{\frac1p-1} |\mu-\nu||x|+ 
2 + C\mu^{-5} (1+|x|+s)^{\frac92} |e_2-e_1|.
\end{multline}
The continuity estimate in~$\mu$ was essentially already obtained above in Step~4 of the proof of Lemma \ref{l.mbarcont}. There we found that $m_\nu\leq m_\mu\leq (\mu/\nu)m_\nu$, and this yields
\begin{equation*} \label{}
\left|m_\mu(x, \H_{e_1}^--se_1)-m_\nu(x, \H_{e_1}^--se_1)\right| \leq C\nu^{-1} |\mu-\nu| m_\nu(x)\leq C\nu^{\frac1p-1} |\mu-\nu| |x|. 
\end{equation*}
We turn to the continuity with respect to $e$. We denote by $z$ the projection of $x$ onto $\partial (\H_{e_1}^+-se_1)$ and let $\tilde s>0$ be such that $x=z+\tilde s e_1$. We also set  $x_2= z+\tilde s e_2$. For later use, we record the following algebraic relations:
$$
\H_{e_1}^-+z= \H_{e_1}^--se_1, \qquad \tilde s=x\cdot e_1+s, \qquad z= x-(x\cdot e_1+s)e_1. 
$$
We split the estimate of $\left|m_\mu(x, \H_{e_1}^--se_1)-m_\mu(x, \H_{e_2}^--se_2)\right|$ into three terms: 
\begin{multline*}
\left|m_\mu(x, \H_{e_1}^--se_1)-m_\mu(x, \H_{e_2}^--se_2)\right| \\
\leq
\left|m_\mu(x, \H_{e_1}^-+z)-m_\mu(x_2, \H_{e_2}^-+z)\right| 
+\left|m_\mu(x_2, \H_{e_2}^-+z)-m_\mu(x, \H_{e_2}^-+z)\right| \\
+ \left|m_\mu(x, \H_{e_2}^-+z)-m_\mu(x, \H_{e_2}^--se_2)\right|.
\end{multline*}
The first term on the right side can be handled as in Step 2 of the proof of Lemma \ref{l.mbarcont}, where we have translated the picture by~$z$. For $\tilde R := C\mu^{-5}(1+ \tilde s^{9/2})$, we have
\begin{multline*}
\left|m_\mu(x, \H_{e_1}^-+z)-m_\mu(x_2, \H_{e_2}^-+z)\right|\\
=\left|m_\mu(z+\tilde s e_1, \H_{e_1}^-+z)-m_\mu(z+\tilde s e_2, \H_{e_2}^-+z)\right| \leq 2 + C \left( \tilde R + \tilde s \right) \left| e_1 - e_2 \right|.
\end{multline*}
The last two terms are controlled by the Lipschitz estimate: 
\begin{align*}
\left|m_\mu(x_2, \H_{e_2}^-+z)-m_\mu(x, \H_{e_2}^-+z)\right| &
\leq L_\mu |x_2-x|  \\
& \leq C \tilde s \, |e_1-e_2| \\
& \leq C (|x|+s)  |e_1-e_2|
\end{align*}
and
\begin{align*}
\left|m_\mu(x, \H_{e_2}^-+z)-m_\mu(x, \H_{e_2}^--se_2)\right|
& \leq L_\mu \dist_H(\H_{e_2}^-+z, \H_{e_2}^--se_2) \\
& \leq  C \left|z\cdot e_2+s\right| \\
& \leq C\left(|x|+s\right) |e_2-e_1|.
\end{align*}
Combining the above inequalities,  we get
\begin{equation*} \label{}
\left|m_\mu(x, \H_{e_1}^--se_1)-m_\mu(x, \H_{e_2}^--se_2)\right|\leq
2 + C\mu^{-5} (1+|x|+s)^{\frac92} |e_2-e_1|.
\end{equation*}
This proves \eqref{RoughEstimmu}. To complete the proof of the lemma, we use~\eqref{RoughGrowthmmue}, which states that
$$
m_\mu(x, \H_{e}^--se)\leq C \mu^{\frac1p} (x\cdot e+s) \leq C  \mu^{\frac1p}  \left( |x|+s\right).
$$
We obtain~\eqref{e.mconte} by interpolating the previous inequality with \eqref{RoughEstimmu} (i.e., we use the latter for large $\mu$ and the former for small $\mu$). 
\end{proof}

\section{The approximate corrector problem}
\label{s.correctors}

In this section, we prove Proposition~\ref{p.correctors}, which links the planar metric and approximate corrector problems. It essentially states that the quantity $\left| \delta v^\delta(0,\xi)+ \overline H(\xi)\right|$ is controlled by the convergence of the planar metric problem to its limit. 

\smallskip

Rather than simply cite comparison results, the quasilinear viscous term and the desire for quantitative results forces us to use of the full uniqueness machinery for viscosity solutions of second-order equations~\cite{CIL} ``inside" the usual perturbed test function argument. The technical details unfortunately obscure the relatively simple and straightforward ideas, and so we recommend to the reader Sections~5 and~6 of~\cite{ACS}, where similar ideas are encountered in a simpler situation. 

\smallskip 

\begin{proof}[{Proof of Proposition~\ref{p.correctors}}]
Fix $\xi \in \Rd \setminus \{ 0 \}$. Denote $\mu:= \overline H(\xi) > 0$ and set $e:= \xi / |\xi|$ and $\H^\pm := \H_e^\pm$. The basic idea of the proof is to show that $v^\delta(\cdot,\xi)$ must be close enough to $m_\mu(\cdot,\H^--se)$, for $s \approx \delta^{-1}$, in a ball of radius $\approx \delta^{-1}$ centered at the origin, that we can infer the limit~\eqref{e.dvdlim} from~\eqref{e.MPlim}.

\smallskip

We assume that for a fixed $\delta,\lambda\in (0,1]$ and $s:= C'/\lambda \delta$, with $C'\geq 1$ to be selected below, we have
\begin{equation}
\label{e.goodevent}
\sup_{x\in B_{s/2}} \big|  m_\mu\left(x,\H^{-} - se\right) - \overline m_\mu(e)(s+(x\cdot e)) \big| \leq \frac{\lambda}{\delta}.
\end{equation}
The goal is to prove that 
\begin{equation}
\label{e.notbadevent}
-\delta v^\delta(0,\xi) \leq \mu + C \lambda^{\frac 15}.
\end{equation}
(The proof that $-\delta v^\delta(0,\xi) \geq \mu - C \lambda^{\frac 15}$ is very similar, and so we omit it.)

\smallskip

The argument for~\eqref{e.notbadevent} is based on a comparison argument, and so as usual we ``double the variables," introducing the auxiliary function
\begin{equation*} \label{}
\Phi(x,y):= m_\mu\left(x,\H^--se \right) - \left( \xi \cdot y + v^\delta(y,\xi) \right) - \lambda \underbrace{\left( ( 1+|x|^2 )^{\frac12} -1  \right)}_{=:\phi(x)} - \frac1{4\ep} \left| x-y \right|^4,
\end{equation*}
where $\ep \in (0,1]$ is selected below.

\smallskip

\emph{Step 1.} We show there exists point $(x_0,y_0) \in B_{s/4} \times B_{s/4}$ such that 
\begin{equation} \label{e.lockmax}
\Phi(x_0,y_0) = \sup_{x,y\in B_{s/2}} \Phi(x,y).
\end{equation} 
We take~$(x_0,y_0) \in \overline B_{s/2} \times \overline B_{s/2}$ to be any point at which $\Phi$ attains its maximum over~$\overline B_{s/2} \times \overline B_{s/2}$. We claim that $x_0$ and $y_0$ must lie in $B_{s/4}$, if we make a suitable choice of~$C'$. Observe that the Lipschitz estimate for $m_\mu$ (see Theorem~\ref{t.LipschitzAppenA1}) implies that
\begin{equation} \label{e.lipbnd1}
\left| x_0 - y_0\right| \leq C\ep^{\frac13}.
\end{equation}
The inequality $\Phi(x_0,y_0) \geq \Phi(0,0)$ gives
\begin{equation*} \label{}
\lambda \phi(x_0) \leq \left| m_\mu(x_0,\H^--se) - m_\mu(0,\H^--se) - \xi\cdot x_0 \right| + \left| \xi \right| \left|x_0-y_0\right|  + C \delta^{-1}.
\end{equation*}
Using the previous line,~\eqref{e.goodevent},~\eqref{e.lipbnd1} and taking $C'\geq C(\data,L)$, we obtain
\begin{equation*} \label{}
\phi(x_0) \leq \frac{2\lambda}{\delta}  +  \frac{C}{\lambda \delta} + \frac{C\ep^{\frac13}}{\lambda}  \leq \frac{C}{\lambda\delta}.
\end{equation*}
Since $\phi$ grows linearly, we get $|x_0| \leq C/(\lambda\delta)$ and then~\eqref{e.lipbnd1} gives $|y_0| \leq C/(\lambda\delta)$. These constants $C$ do not depend on the choice of $C'$, therefore we may further enlarge $C'$, if necessary, to obtain that $|x_0|+|y_0|< \frac14s$. This completes the proof of the claim.

\smallskip

\emph{Step 2.} Fix $\gamma>0$. According to Step~1 and the maximum principle for semicontinuous functions~\cite[Theorem 3.2]{CIL}, there exist symmetric matrices $X,Y\in\R^{d\times d}$ satisfying 
\begin{equation} \label{e.matgamma}
-\left( \frac{1}{\gamma} + \left| M \right| \right) I_{2d} \leq
\begin{pmatrix} 
X &0 \\  0 & -Y  
\end{pmatrix} 
\leq 
M + \gamma M^2
\end{equation}
and 
\begin{equation*} \label{}
\left\{ 
\begin{aligned}
& \left(X +\lambda D^2\phi(x_0),  \eta + \lambda D\phi(x_0) \right)\in \overline{\mathcal{J}}^{2,+}\left( m_\mu(\cdot,\H^--se) \right)(x_0), \\
& \left(Y,  \eta -\xi  \right) \in \overline{\mathcal{J}}^{2,-}\left(  v^\delta (\cdot,\xi) \right)(y_0).
\end{aligned}
\right.
\end{equation*}
where $\eta:= \ep^{-1} \left( x_0-y_0\right)^2 \left| x_0-y_0\right|$ and $M\in \R^{d\times d}$ is defined by 
\begin{equation} \label{e.matineq10}
M:= \frac1\ep \begin{pmatrix} 
N & -N  \\ -N & N
\end{pmatrix}, \quad N:= \left| x_0 - y_0 \right|^2 I_d + 2\left( x_0-y_0\right)\otimes \left( x_0-y_0\right).
\end{equation}
We choose $\gamma = \ep \left| x_0 -y_0\right|^{-2}$ so that~\eqref{e.matgamma} becomes
\begin{equation*} \label{}
- \frac{C\left| x_0 -y_0\right|^{2}}{\ep}  I_{2d}  \leq \begin{pmatrix} 
X &0 \\  0 & -Y  
\end{pmatrix}  \leq \frac{C\left| x_0 -y_0\right|^{2}}{\ep}\begin{pmatrix} 
I_d & -I_d  \\ -I_d & I_d
\end{pmatrix}.
\end{equation*}
Using the equations for $m_\mu$ and $v^\delta$, we obtain
\begin{equation} \label{e.mmuX}
-\tr\left( A \left( \eta + \lambda D\phi(x_0),x_0\right) 
\left( X + \lambda D^2\phi(x_0) \right) \right) 
+ H\left( \eta + \lambda D\phi(x_0), x_0\right) \leq \mu
\end{equation}
and
\begin{equation} \label{e.dvdY}
-\tr\left( A \left(\eta,y_0\right) Y \right) 
+ H\left( \eta, y_0\right) \geq -\delta v^\delta(y_0,\xi).
\end{equation}
The rest of the proof is concerned with deriving~\eqref{e.notbadevent} from~\eqref{e.lipbnd1},~\eqref{e.matgamma},~\eqref{e.mmuX} and~\eqref{e.dvdY}.

\smallskip

\emph{Step 3.} We next obtain a lower bound for $|x_0-y_0|$ and $|\eta|$. The claim is that 
\begin{equation} \label{e.x0y0lowerbound}
\left| x_0 - y_0 \right| \geq c \mu^{\frac12} \ep^{\frac12}\qquad {\rm and}\qquad |\eta|\geq c\mu^{\frac32}\ep^{\frac12}.
\end{equation}
Using~\eqref{e.dvdY} and \eqref{e.RoughBoundVdelta}, we have 
\begin{equation*} \label{}
c\mu = c\overline H(\xi) \leq c_0 \left| \xi \right|^p \leq -\delta v^\delta(y_0,\xi)  \leq -\tr\left( A\left( \eta,y_0\right) Y \right) + H\left( \eta,y_0 \right) \leq C_0 \left( |M|+\left|\eta\right|^p \right).
\end{equation*}
As $|M| \leq C\ep^{-1}|x_0 - y_0|^2$ and $|\eta|=\ep^{-1} \left| x_0 - y_0 \right|^{3}$, we get~\eqref{e.x0y0lowerbound}.

\smallskip

\emph{Step 4.} We complete the proof of~\eqref{e.notbadevent}. Using~\eqref{e.lipbnd1} and the matrix inequality \eqref{e.matineq10}, we obtain
\begin{align*}
\tr\left( A(\eta,x_0)X - A(\eta,y_0)Y\right) &  = \tr\left( \begin{pmatrix} 
X &0 \\  0 & -Y  
\end{pmatrix}  \begin{pmatrix} 
\sigma(\eta,x_0)  \\  \sigma(\eta,y_0)   
\end{pmatrix}  \begin{pmatrix} 
\sigma(\eta,x_0)  \\  \sigma(\eta,y_0)   
\end{pmatrix}^T  \right)  \\
& \leq C\ep^{-\frac13} \left| \sigma(\eta,x_0) - \sigma(\eta,y_0) \right|^2 \\
& \leq C\ep^{-\frac13} \left| x_0-y_0\right|^2 
 \leq C\ep^{\frac13}.
\end{align*}
Thus 
\begin{align*} \label{}
\lefteqn{ \tr\left( A \left( \eta + \lambda D\phi(x_0),x_0\right) X \right) -\tr\left( A \left(\eta,y_0\right) Y \right) } 
\qquad & \\
& \leq \left|A \left( \eta + \lambda D\phi(x_0),x_0\right) -A \left( \eta ,x_0\right) \right| \left| X \right| + \tr\left(A \left( \eta  ,x_0\right) X -A \left( \eta ,y_0\right)Y \right) \\
& \leq \left| \eta \right|^{-1} \lambda \left| X \right| + C\ep^{\frac13} \\
& \leq C \left( \lambda\mu^{-\frac32} \ep^{-\frac56}  + \ep^{\frac13} \right).
\end{align*}
Similarly,
\begin{equation*} \label{}
\left| H\left( \eta + \lambda D\phi(x_0),x_0 \right) - H\left( \eta ,y_0 \right) \right|  \leq C \left| \eta \right|^{p} \left( \lambda + \left| x_0 - y_0 \right| \right) \leq C \left( \lambda + \ep^{\frac13} \right). 
\end{equation*}
Subtracting~\eqref{e.mmuX} from~\eqref{e.dvdY} and combining the above inequalities yields 
\begin{equation*}
-\delta v^\delta(y_0,\xi) \leq \mu+C \left( \lambda\mu^{-\frac32} \ep^{-\frac56}  + \ep^{\frac13} \right). 
\end{equation*}
We now select the value of $\ep$ by optimizing the last expression, which leads to the choice $\ep:= \lambda^{\frac67}\mu^{-\frac97}$. We obtain
\begin{equation} \label{e.notbady0}
-\delta v^\delta(y_0,\xi) \leq \mu + C \mu^{-\frac37} \lambda^{\frac27}.
\end{equation}
We need to obtain the same inequality at the origin instead of~$y_0$. Since $|x_0-y_0|\leq |x_0| + |y_0| \leq s/2$,~\eqref{e.lockmax} gives that $\Phi(x_0-y_0,0) \leq \Phi(x_0,y_0)$. This inequality implies
\begin{equation*} \label{}
v^\delta(y_0,\xi) - v^\delta(0,\xi) \leq m_\mu(x_0,\H^--se) - m_\mu(x_0-y_0,\H^--se) - \xi\cdot y_0 +C\left|x_0-y_0\right|.
\end{equation*}
Another application of~\eqref{e.goodevent} then yields
\begin{equation*} \label{}
v^\delta(y_0,\xi) - v^\delta(0,\xi) \leq \frac{C\lambda}\delta + C\ep^{\frac13} \leq  \frac{C\lambda}\delta + C \mu^{-\frac37} \lambda^{\frac27}.
\end{equation*}
Multiplying by $\delta$, using $\lambda,\delta \leq1$ and combining with~\eqref{e.notbady0} gives
\begin{equation} \label{e.notbady1}
-\delta v^\delta(0,\xi) \leq \mu + C \mu^{-\frac37} \lambda^{\frac27}.
\end{equation}
This estimate obviously degenerates for small $\mu$. Therefore, to obtain~\eqref{e.notbadevent}, it is necessary to interpolate the previous inequality with the deterministic bound from~\eqref{e.RoughBoundVdelta}, which is useful precisely for small $\mu$:
\begin{equation*} \label{}
-\delta v^\delta(0,\xi) \leq C_0 \left|\xi\right|^p \leq C \mu \leq \mu + C\mu.
\end{equation*}
We obtain
\begin{equation*} \label{}
-\delta v^\delta(0,\xi) \leq  \mu + C \mu\wedge \left( \mu^{-\frac37} \lambda^{\frac27} \right) \leq \mu + C\lambda^{\frac15}. 
\end{equation*}
This completes the proof of~\eqref{e.notbadevent}.
\end{proof}

\begin{proof}[Proof of Proposition \ref{p.Hbarprops}] From Propositions~\ref{p.planarquantitative} and~\ref{p.correctors}, we deduce that 
\begin{equation*} \label{}
\delta v^\delta(0,\xi) \to -\overline H(\xi) \quad \mbox{in probability (with respect to $\P$) as $\delta\to 0$.}
\end{equation*}
The growth bound \eqref{e.Hbarpgrowth1} is therefore immediate from \eqref{e.RoughBoundVdelta}, while \eqref{e.VdeltaHolder} and \eqref{e.Hbarpgrowth1} imply that, for each $R>0$ and for each $\xi, \eta \in B_R$, 
$$
\left| \overline H(\xi)- \overline H(\eta)\right|\leq  \left( L_R \left(|\xi|\wedge |\eta|\}\right)^{-\frac{2p}7} |\xi-\eta|^{\frac27} \right) \wedge\left( C_0\left(|\xi|+|\eta|\right)\right),
$$
from which \eqref{e.Hbarholder1} follows. 
\end{proof}

\section{Homogenization of the time-dependent problem}
\label{s.timedependent}

In this section, we prove Proposition~\ref{p.PTFM}, which roughly asserts that the limit~\eqref{e.dvdlim} for the approximate correctors controls the homogenization of the full time-dependent initial-value problem for~\eqref{e.pde}. 
The argument is a quantitative version of the so-called perturbed test function argument, a technical device introduced in~\cite{E2}. Like the result in the previous section, this is a purely deterministic PDE fact derived from a comparison argument using the ``full" uniqueness machinery for second-order equations. 

\smallskip

 It has been pointed out~\cite{CM} that, for the mean curvature equation and other singular quasilinear equations, the argument from~\cite{E2} does not apply in a straightforward way and some regularization of the solutions is needed. We see this in our approach by the fact that, at the end of the argument, we cannot send the parameters to zero in a sequence to obtain a qualitative result; rather we have to send several of the parameters to zero at the same time. 

\begin{proof}[{Proof of Proposition~\ref{p.PTFM}}]
Fix $0<\ep \leq \delta \leq \lambda \leq 1$. We will argue that 
\begin{equation} \label{e.errorassump}
\sup_{(x,t) \in B_R \times [0,T] } \left( u^\ep (x,t) - u(x,t) \right) \geq \lambda 
\end{equation}
implies
\begin{equation} \label{e.bigtdassum}
\inf_{(x,\xi) \in  B_{C/  \ep  } \times \overline B_L} \left( -\delta v^\delta (x,\xi) - \overline{H}(\xi) \right) \leq -\frac\lambda{C}+ C \left( \frac{\ep}{\delta\lambda} \right)^{\frac1{10}}.
\end{equation}
This will prove only half of the proposition, since we also must establish that
\begin{multline} \label{e.otherdirection}
\inf_{(x,t)\in B_R\times[0,T]} \left( u^\ep(x,t) - u(x,t) \right) \leq -\lambda \\ 
\implies \quad \sup_{(x,\xi) \in  B_{C / \ep  } \times \overline B_L} \left( -\delta v^\delta (x,\xi) - \overline{H}(\xi) \right) \geq \frac\lambda{C} -   C \left( \frac{\ep}{\delta\lambda} \right)^{\frac1{10}}. 
\end{multline}
However, the argument we give for~$\mbox{\eqref{e.errorassump}} \Rightarrow \mbox{\eqref{e.bigtdassum}}$ can be modified in a straightforward way to yield~\eqref{e.otherdirection}. Therefore we omit the proof of~\eqref{e.otherdirection}. For the rest of the proof, we assume that~\eqref{e.errorassump} holds. 

\smallskip

Throughout the argument,~$C$ and~$c$ denote positive constants which depend on~$(\data,L,R,T)$ and may vary in each occurrence. We will prove the result under the further assumption that $\lambda\leq c$. We obtain the same result for all $\lambda \in(0,1]$ by adjusting the constant $C$ in the first term on the right side of~\eqref{e.bigtdassum}.

\smallskip

\emph{Step 1.} We begin by doubling the variables to find an initial touching point. This is not the main comparison argument: the objective here is merely to get an initial direction~$\xi_0$ which is then used in the next step to introduce the appropriate approximate corrector.

\smallskip

We fix parameters $\alpha\in [\ep,\lambda]$ and $\gamma\in (0,\lambda]$ to be selected below and define a first auxiliary function $\Phi: \Rd \times\Rd\times [0,T]\to\R$ by
\begin{equation*} \label{}
\Phi(x,y,t):= u^\ep(x,t) - u(y,t) -  \frac1{2\alpha}\left| x-y\right|^2 - \gamma \phi(x) - \left(\frac{\lambda}{2T} \right)t,
\end{equation*}
where $\phi:\Rd\to\R$ is the function given by $\phi(x):= \left( 1+|x|^2 \right)^{\frac12} - 1$. In order to check that $\Phi$ attains its supremum on~$\Rd\times\Rd\times[0,T]$, we notice that, according to~\eqref{e.LipLs} and the initial condition in~\eqref{e.ptfmpdes}, for every $(x,y,t) \in\Rd \times\Rd \times[0,T]$, 
\begin{align*}
\Phi(x,y,t) & \leq Lt + L\left |x-y \right|  - \frac1{2\alpha} \left |x-y \right|^2  -\gamma \phi(x) 
\\
& \leq Lt + \alpha L^2 - \frac1{4\alpha} \left| x-y \right|^2 - \gamma \left( |x| -1 \right) 
.
\end{align*}
The assumption~\eqref{e.errorassump} implies the existence of~$(\hat x, \hat t) \in B_R \times [0,T]$ for which
\begin{equation*} \label{}
\Phi\left(\hat x, \hat x, \hat t\right) \geq \lambda - \gamma \phi\left( \hat x \right)- \frac12\lambda \geq \frac12\lambda - \gamma R 
\end{equation*}
Imposing the restriction that $\gamma \leq \frac14\lambda R^{-1}$, we get
\begin{equation*} \label{}
\Phi\left(\hat x, \hat x, \hat t\right) \geq \frac14\lambda. 
\end{equation*}
Note that any $(x,y,t)$ such that $\Phi(x,y,t) \geq \frac18\lambda$ must then satisfy
\begin{equation*}
\frac18\lambda \leq Lt + \alpha L^2 + \gamma - \frac{1}{4\alpha} \left| x-y \right|^2 - \gamma |x|. 
\end{equation*}
In particular, 
\begin{equation} \label{e.boundxy}
\gamma |x| + \frac1{4\alpha} \left| x-y \right|^2 \leq LT + \alpha L^2 + 1 \leq C \left( LT + L^2 \right).
\end{equation}
and 
\begin{equation*} \label{}
t \geq \frac \lambda {8L} - \alpha L - \frac{\gamma}{L}. 
\end{equation*}
Imposing the restrictions $\gamma \leq \frac1{32} \lambda$ and $\alpha \leq \frac1{32} L^{-2}\lambda$ yields 
\begin{equation} \label{e.tpositive}
t \geq \frac \lambda {16L} > 0.
\end{equation}
We deduce that the set where $\Phi$ is at least $\frac 14 \lambda$ is nonempty, bounded and contained in the time interval $\left[ \frac \lambda {16L}, T\right]$. Thus there exists~
$$(x_0,y_0,t_0) \in \Rd\times\Rd\times\left[\tfrac\lambda {16L},T\right]$$
for which 
\begin{equation}
\label{e.Phimax}
\Phi(x_0,y_0,t_0) = \sup_{(x,y,t)\in \Rd\times\Rd\times[0,T]} \Phi(x,y,t) \geq \Phi\left( \hat x, \hat x, \hat t \right) \geq \frac14 \lambda. 
\end{equation}
It is useful to denote
\begin{equation*} \label{}
\xi_0:= \frac{x_0-y_0}{\alpha}.
\end{equation*}
According to the Lipschitz assumption~\eqref{e.LipLs} for $u(\cdot,t_0)$, we have
\begin{equation} \label{e.xi0lip}
\left| \xi_0 \right| \leq L. 
\end{equation}

\smallskip 

\emph{Step 2.} We introduce the approximate corrector, which requires a new auxiliary function and a tripling of the variables, and find a new touching point. Fixing an additional parameter $\beta \in (0,\ep]$, to be selected below, we set
\begin{multline*} \label{}
\Psi(x,y,z,t):= u^\ep(x,t) - u(y,t) - \ep v^\delta \left( \frac z\ep,\xi_0\right) + \xi_0\cdot (x-z) -\frac1{2\alpha}\left| x-y\right|^2 \\ - \frac1{4\beta^3} \left| x-z\right|^4 
 - \frac1{2\alpha} \left( \left| x-x_0\right|^2+\left| y-y_0 \right|^2 \right) - \gamma \phi(x) - \left( \frac{\lambda}{2T}\right)t.
\end{multline*}
Since, by \eqref{e.RoughBoundVdelta},
\begin{align*} \label{}
\Psi(x,y,z,t) & \leq \Phi(x,y,t) + \frac{C_0L^p\ep}{\delta} + L \left| x-z \right|- \frac1{4\beta^3} \left| x-z\right|^4 \\
& \leq  \Phi(x,y,t) + \frac{C_0L^p\ep}{\delta} + CL^{\frac43}\beta^3  -  \frac{1}{8\beta^3} \left| x-z \right|^4
\end{align*}
and
\begin{equation*} \label{}
\Psi(x_0,y_0,x_0,t_0) = \Phi(x_0,y_0,t_0) -  \ep v^\delta\left( \frac{x_0}{\ep} , \xi_0 \right) \geq  \Phi(x_0,y_0,t_0) \geq \frac14\lambda,
\end{equation*}
we obtain the existence of $(x_\beta,y_\beta,z_\beta,t_\beta)\in \Rd\times\Rd\times\Rd\times[0,T]$ such that 
\begin{equation}
\label{e.Psimax}
\Psi(x_\beta,y_\beta,z_\beta,t_\beta) = \sup_{(x,y,z,t)\in \Rd\times\Rd\times\Rd\times[0,T]} \Psi(x,y,z,t) \geq \frac14\lambda. 
\end{equation}
If we impose the restrictions $C_0L^p \ep \leq \frac1{16} \lambda \delta$ and $CL^{\frac43}\beta^3 \leq \lambda$, then we obtain in particular that 
\begin{equation*}
\Phi(x_\beta ,y_\beta,t_\beta ) \geq \frac12\lambda - \frac{C_0L^p\ep}{\delta}  - CL^{\frac43} \beta^3 \geq \frac18\lambda.
\end{equation*}
This implies by Step~1 that 
\begin{equation}
\label{e.trappings}
t_\beta \geq \frac{\lambda}{16L} \quad \mbox{and} \quad \left| x_\beta \right| \leq \frac{CL}{\gamma} \left( L+T \right).
\end{equation}
In particular, $t_\beta>0$.

\smallskip

We conclude this step with some estimates which are needed in the sequel. First we observe that the Lipschitz estimate for $v^\delta(\cdot,\xi_0)$ gives
\begin{equation*} \label{}
\left|  \frac{1}{\beta^3} \left| x_\beta-z_\beta \right|^2 \left(x_\beta-z_\beta \right) - \xi_0\right| \leq C.
\end{equation*}
Hence
\begin{equation} \label{e.xtzt}
\left| x_\beta-z_\beta \right| \leq C\beta
\end{equation}
and in particular, by~\eqref{e.trappings},  
\begin{equation}
\label{e.zthetatrapp}
\left| z_\beta \right| \leq \frac{CL}{\gamma}\left( L+T \right) + C\beta \leq \frac{CL}{\gamma}\left( L+T \right).
\end{equation}

We next claim that
\begin{equation}
\label{e.thetasclose}
\left| x_\beta - x_0 \right|^2 + \left| y_\beta - y_0 \right|^2   \leq \frac{C\ep\alpha}{\delta}.
\end{equation}
To see this, we use the fact that
\begin{align*} \label{}
\Phi(x_\beta,y_\beta,t_\beta) \leq \Phi(x_0,y_0,t_0) & = \Psi(x_0,y_0,x_0,t_0) + \ep v^\delta\left( \frac{x_0}{\ep} , \xi_0 \right) \\
& \leq \Psi(x_\beta,y_\beta,z_\beta,t_\beta) + \ep v^\delta\left( \frac{x_0}{\ep} , \xi_0 \right).
\end{align*}
This can be expressed equivalently (recall $\beta \leq \ep \leq \ep\delta^{-1}$) as
\begin{align*}
\lefteqn{
\frac1{2\alpha}  \left(\left| x_\beta-x_0\right|^2+\left| y_\beta-y_0\right|^2 \right)
} 
\qquad & \\
& \leq \ep v^\delta\left( \frac{x_0}{\ep} , \xi_0 \right) - \ep v^\delta\left( \frac{z_0}{\ep} , \xi_0 \right) + \xi_0 \cdot \left( x_\beta-z_\beta\right) - \frac{1}{4\beta^3} \left| x_\beta - z_\beta \right|^4 \\
& \leq C\left( \frac{\ep}{\delta} +  \beta \right) \leq \frac{C\ep}{\delta},
\end{align*}
which yields~\eqref{e.thetasclose}. 

\smallskip

\emph{Step 3.} The application of the parabolic version of the maximum principle for semicontinuous functions~\cite[Theorem 8.3]{CIL}. We obtain, for every $\eta>0$, the existence of $d$-by-$d$ symmetric matrices $X$, $Y$ and $Z$ and $b\in\R$ such that 
\begin{equation}
\label{e.triplematfun}
-\left( \frac{1}{\eta} + \left| M \right| \right) I_{3d} \leq
\begin{pmatrix} 
X &0 & 0 \\  0 & - Y & 0 \\ 0 & 0 & -\ep^{-1}Z 
\end{pmatrix} 
\leq 
M + \eta M^2,
\end{equation}
and
\begin{equation*} \label{}
\left\{
\begin{aligned}
& \left(X,  \xi_\beta+ \xi',b \right) \in \overline{\mathcal{P}}^{2,+} u^\ep(x_\beta,t_\beta), \\
& \left(Y, \xi_0+\xi'',b - \frac{\lambda}{2T} \right) \in \overline{\mathcal{P}}^{2,+} u(y_\beta,t_\beta), \\
& \left(Z, -\xi_0+ \xi_\beta) \right) \in \overline{\mathcal{J}}^{2,-} v^\delta(\cdot,\xi_0)\left( \frac{z_\beta}{\ep}\right),
\end{aligned} 
\right.
\end{equation*}
where
\begin{equation*}
M := \frac{1}{\alpha} \begin{pmatrix} 
2I_d & -I_d & 0\\ -I_d & 2I_d & 0 \\ 0 & 0 & 0
\end{pmatrix}  
+ \frac{1}{\beta^3} \begin{pmatrix} 
N & 0 & -N\\ 0 & 0 & 0 \\ -N & 0 & N
\end{pmatrix} 
+ \gamma \begin{pmatrix} 
P & 0 & 0\\ 0 & 0 & 0 \\ 0 & 0 & 0
\end{pmatrix} 
\end{equation*}
for 
\begin{equation*} \label{}
N:= \left| x_\beta-z_\beta \right|^2 I_d + 2\left( x_\beta-z_\beta\right)\otimes\left( x_\beta-z_\beta\right), \quad
 P:=  D^2\phi(x_\beta),
\end{equation*}
and
\begin{equation*}
\left\{
\begin{aligned}
& \xi_\beta:=  \frac{1}{\beta^3}\left| x_\beta-z_\beta \right|^2(x_\beta-z_\beta), \\
& \xi':= \frac1\alpha \left( 2(x_\beta-x_0)-(y_\beta-y_0)\right)   + \gamma D\phi(x_\beta), \\
& \xi'' :=  \frac1\alpha \left( (x_\beta-x_0)-2(y_\beta-y_0)\right).
\end{aligned}
\right.
\end{equation*}
Observe that (recall that $\beta \leq \ep,\alpha$ and $\gamma\leq1$)
\begin{align}
\label{e.Mbound}
\left| M \right| 
 \leq C \left( \frac{1}\alpha + \frac{|x_\beta-z_\beta|^2}{\beta^3} + \gamma |P| \right) \leq \frac{C}\beta
\end{align}
and, according to~\eqref{e.thetasclose},
\begin{equation}
\label{e.xiprbnds}
\left| \xi' \right| \leq  C\left(\frac{\ep}{\alpha\delta}\right)^{\frac12}+ \gamma \quad \mbox{and} \quad \left| \xi'' \right| \leq  C\left(\frac{\ep}{\alpha\delta}\right)^{\frac12}.
\end{equation}
Taking $\eta := \left| M \right|^{-1}$ in the matrix inequality, in view of the fact that $M\geq 0$, gives
\begin{equation}
\label{e.triplematfun2}
-C\left| M \right| I_{3d} \leq
\begin{pmatrix} 
X &0 & 0 \\  0 & - Y & 0 \\ 0 & 0 & -\ep^{-1}Z 
\end{pmatrix} 
\leq 
CM.
\end{equation}
As $t_\beta>0$, the equations satisfied by~$u^\ep$,~$u$ and~$v^\delta$ yield
\begin{equation}
\label{e.uepplug}
-\ep\tr\left( A \left(\xi_\beta+ \xi',\frac{x_\beta}{\ep}  \right) X \right) + H\left( \xi_\beta+ \xi',\frac{x_\beta}{\ep} \right)
 \leq -b,
\end{equation}
\begin{equation}
\label{e.uplug}
\overline{H}\left( \xi_0+\xi'' \right) \geq -b+ \frac{\lambda}{2T}, 
\end{equation}
and
\begin{equation}
\label{e.dvdplug}
-\tr\left( A \left(  \xi_\beta ,\frac{z_\beta}{\ep} \right) Z\right) + H\left(\xi_\beta,\frac{z_\beta}{\ep} \right) \geq -\delta v^\delta\left( \frac{z_\beta}{\ep},\xi_0\right).
\end{equation}

The rest of the proof is devoted to combining the previous five inequalities to obtain the desired conclusion. Due to the singularity of the diffusion term, it is natural to split the argument into two cases, depending on the size of $\left|\xi_\beta\right|$. For this purpose, we fix an additional parameter $s\in [\ep,1]$, to be selected below. 

\smallskip

\emph{Step 4.} We consider the case that $\left| \xi_\beta \right| \leq s$. Then~\eqref{e.triplematfun2} and revisiting the estimate for $|M|$ yields 
\begin{equation*}
\left| Z \right| \leq C \ep |M| \leq C\ep \left( \frac1\alpha  + \frac{s^{\frac23}}\beta \right)
\end{equation*}
and the left side of~\eqref{e.dvdplug} can then be estimated brutally:
\begin{equation*} \label{}
-\delta v^\delta\left( \frac{z_\beta}{\ep},\xi_0 \right) \leq -\tr\left( A \left(  \xi_\beta ,\frac{z_\beta}{\ep} \right) Z\right) + H\left(\xi_\beta,\frac{z_\beta}{\ep} \right) \leq C\ep \left( \frac1\alpha  + \frac{s^{\frac23}}\beta \right) + C s^p.
\end{equation*}
Since $s^p\leq s$ and $\overline{H}\geq 0$, we deduce that 
\begin{equation}\label{e.case1}
-\delta v^\delta\left( \frac{z_\beta}{\ep},\xi_0 \right) \leq  \overline{H}(\xi_0) - \frac\lambda {2T} + C \left( \frac\ep\alpha + \frac{ s^{\frac23}\ep}\beta \right) + C s. 
\end{equation}

\smallskip

\emph{Step 5.} We consider the case that $\left| \xi_\beta \right| \geq s$. Then we have
\begin{align*}
\left| A \left(\xi_\beta+ \xi',\frac{x_\beta}{\ep}\right)X - A \left(\xi_\beta,\frac{x_\beta}{\ep}\right) X \right| 
& \leq C \left| \xi_\beta \right|^{-1} \left| \xi' \right| \left| X \right| \\
& \leq \frac C s \left( \left(\frac{\ep}{\alpha\delta}\right)^{\frac12}+\gamma  \right) \left| M \right|.\notag
\end{align*}
Using~\eqref{e.Mbound}, we get
\begin{equation}
\label{e.case2A1}
\ep\left| A \left(\xi_\beta+ \xi',\frac{x_\beta}{\ep}\right)X - A \left(\xi_\beta,\frac{x_\beta}{\ep}\right) X \right| \\
\leq \frac {C\ep}{s\beta} \left(   \left(\frac{\ep }{\alpha\delta}\right)^{\frac12}+\gamma  \right).
\end{equation}
Multiplying the second inequality of~\eqref{e.triplematfun2} by the matrix
\begin{equation*}
 \begin{pmatrix} 
\sigma\left(\xi_\beta,\frac{x_\beta}{\ep}\right)  \\  0 \\ \sigma\left(\xi_\beta,\frac{z_\beta}{\ep}\right)   
\end{pmatrix}  \begin{pmatrix} 
\sigma\left(\xi_\beta,\frac{x_\beta}{\ep}\right)  \\  0 \\ \sigma\left(\xi_\beta,\frac{z_\beta}{\ep}\right)
\end{pmatrix}^T, 
\end{equation*} 
taking the trace of the resulting expression and using~\eqref{e.xtzt} yields 
\begin{align}
\label{e.case2A2}
\lefteqn{
-\tr\left( \ep A \left(  \xi_\beta ,\frac{x_\beta}{\ep} \right) X -A \left(  \xi_\beta ,\frac{z_\beta}{\ep} \right) Z\right)
}
\qquad \qquad & \\
& \leq \frac{C\ep\left| N \right|}{\beta^3}  \left| \sigma\left(\xi_\beta,\frac{x_\beta}{\ep}\right)- \sigma\left(\xi_\beta,\frac{z_\beta}{\ep}\right)\right|^2 +   \frac{C\ep}\alpha \notag \\
& \leq \frac{C}{\ep\beta} \left| z_\beta - x_\beta \right|^2+   \frac{C\ep}\alpha  \notag \\
& \leq C\left( \frac\beta\ep + \frac\ep\alpha   \right).\notag
\end{align}
For the $H$ terms, we similarly have
\begin{align}
\label{e.case2H}
\left| H\left( \xi_\beta+ \xi',\frac{x_\beta}{\ep} \right) - H\left(\xi_\beta,\frac{z_\beta}{\ep} \right) \right| 
& \leq C \left( \left| \xi' \right| +\ep^{-1} \left| x_\beta - z_\beta \right|  \right) \\
& \leq C \left( \frac{\beta}{\ep} +   \left(\frac{\ep }{\alpha\delta}\right)^{\frac12} + \gamma \right).  \notag
\end{align}
By the H\"older  continuity of $\overline{H}$ given~\eqref{e.Hbarholder1} combined with~\eqref{e.xiprbnds} gives
\begin{equation} \label{e.putholder}
\left| \overline{H}(\xi_0+\xi'') - \overline{H}(\xi_0) \right| \leq CK \left| \xi''\right|^a \leq CK \left( \frac{\ep }{\alpha \delta}\right)^{\frac 17}.
\end{equation}
Combining~\eqref{e.uepplug},~\eqref{e.uplug},~\eqref{e.dvdplug},~\eqref{e.case2A1},~\eqref{e.case2A2},~\eqref{e.case2H} and~\eqref{e.putholder} yields (recall that $\beta\leq \ep$ and $s\leq 1$)
\begin{equation*}
-\delta v^\delta\left(\frac{z_\beta}{\ep},\xi_0 \right) \leq  \overline{H}(\xi_0) -\frac{\lambda}{2T}+ C\mathcal E,
\end{equation*}
where we define
\begin{equation*} \label{}
\mathcal E:=  \frac{\beta}{\ep} + \frac\ep\alpha + \frac {\ep}{s\beta} \left(   \left(\frac{\ep }{\alpha\delta}\right)^{\frac12}+\gamma  \right) + K \left( \frac{\ep }{\alpha \delta}\right)^{\frac 17}.
\end{equation*}

\smallskip

\emph{Step 6.} The conclusion. Combining the results of Steps~4 and~5, we have shown that 
\begin{equation}
\label{e.bigaldn}
-\delta v^\delta\left(\frac{z_\beta}{\ep},\xi_0 \right) \leq   \overline{H}(\xi_0) - \frac{\lambda}{2T}+ C \mathcal E',
\end{equation}
for
\begin{equation*}
\mathcal E':= \frac{s^{\frac23}\ep } {\beta} + s + \frac{\beta}{\ep} + \frac\ep\alpha + \frac {\ep}{s\beta} \left(   \left(\frac{\ep }{\alpha\delta}\right)^{\frac12}+\gamma  \right) + K \left( \frac{\ep }{\alpha \delta}\right)^{\frac 17}.
\end{equation*}
Now we need to optimize the parameters in order to minimize the error~$\mathcal E'$. It may not be immediately obvious that the parameters can be chosen in such a way to make~$\mathcal E'$ small, ensuring even a qualitative proof. It is thus reassuring to notice that, if we take~$\beta$ to be a power of~$\ep$ slightly larger than one, $s$ to be a very small positive power of~$\ep$ and $\gamma$ to be equal to the other term sharing the parentheses  with it, then we can send~$\ep\to 0$ to get~$\mathcal E' \to 0$. 

\smallskip

To obtain a quantitative bound, an analysis leads to the following choices of the parameters (recall that we consider $\ep \leq \delta \leq \lambda$ to be given):
\begin{equation*} \label{}
\alpha: = \frac \lambda{32L},  \quad \beta:= \ep \left( s^{\frac13} +  s^{-\frac12} \left( \frac{\ep}{\alpha \delta} \right)^{\frac14} \right), 
\quad  s:= \left( \frac{\ep }{\alpha\delta } \right)^{\frac 3{10}}, \quad  \gamma:= \frac1{2 } \left( \frac{\ep }{\alpha\delta} \right)^{\frac12} \wedge \frac{\lambda}{32+4R}.
\end{equation*}
It is straightforward to check that each of the constraints we imposed on the parameters is fulfilled provided that $\lambda \leq c$ for a sufficiently small positive constant~$c$.  We obtain
\begin{equation*} \label{}
\mathcal E' = C \left( \frac{\ep}{\delta\lambda} \right)^{\frac1{10}}.
\end{equation*}
Moreover, due to~\eqref{e.xi0lip} and~\eqref{e.zthetatrapp}, which in view of the choice of~$\gamma$ yields
\begin{equation*} \label{}
\left| z_\beta \right| \leq \frac{C}{\gamma} \leq C \lambda^{\frac 12} \delta^{\frac12} \ep^{-\frac 12}   + C\lambda^{-1}  \leq C\ep^{-1},
\end{equation*}
provided that $\lambda \leq c$ for sufficiently small~$c$. Thus~\eqref{e.bigaldn} gives
\begin{equation*} \label{}
\inf_{(x,\xi) \in \overline B_{C/\ep}\times \overline B_L} \left( -\delta v^\delta (x,\xi) - \overline{H}(\xi) \right) \leq -\frac\lambda {2T} + C \left( \frac{\ep}{\delta\lambda} \right)^{\frac1{10}},
\end{equation*}
as desired. This completes the proof of~\eqref{e.bigtdassum}.
\end{proof}

\appendix

\section{Well-posedness and Lipschitz estimates}
\label{s.metricproblem}

The section is devoted to the well-posedness and global Lipschitz regularity of the metric problem and of the approximate corrector equation. 

\subsection{The role of the (LS) condition}
We start our discussion by explaining the technical role of the (LS) condition \eqref{e.LScoercivity}. Beside \eqref{e.LScoercivity}, we suppose that the pair $(\sigma, H)\in\Omega$ and thus satisfies the standing assumptions \eqref{e.Asigma}, \eqref{e.sigma}, \eqref{e.H}, \eqref{e.Hposhom} and \eqref{e.HLip}.
Note that, by $C^1$ regularity and the homogeneity assumptions on $\sigma$, we have:
\begin{equation}\label{e.LipschA}
|D_\xi A(\xi,y)|\leq C_0 \left|\xi\right|^{-1} \qquad \forall \xi\in\Rd\setminus\{ 0 \}.
\end{equation}

\smallskip

The (LS) condition~\eqref{e.LScoercivity} is devised to provide Lipschitz bounds on solutions. Its usefulness is captured in the following technical lemma.

\begin{lemma}\label{l.LSTechniPt} Let $\kappa $ and $\rho$  be as in the (LS) condition. Fix $\gamma>0$, $M\in \R$, $\xi,\eta\in \Rd\setminus\{0\}$, $x, y\in \Rd$ and symmetic matrices $X,Y\in \R^{d\times d}$. Assume that  
\begin{equation*}
-\tr\left(A(\xi,x)X\right) +H(\xi, x)\leq M\leq -\tr\left(A(\eta,y)Y\right) +H(\eta, y)
\end{equation*}
and
\begin{equation*} 
\begin{pmatrix} 
X &0 \\  0 & -Y  
\end{pmatrix} 
\leq (1+\kappa) \gamma \begin{pmatrix} 
I_d & -I_d \\ -I_d & I_d  
\end{pmatrix}.
\end{equation*}
Assume also that $\gamma (x-y)= \xi$ and $(1+\gamma) |\xi-\eta|\leq \kappa$. Then 
\begin{equation*}
\left|\xi\right| \leq \rho^{-1}\left(\frac{2M^2}{\kappa}+8 C_0 |\sigma|^2\kappa^2\right) \vee \left(\frac{1+M}{c_0}\right)^{\frac1p}.
\end{equation*}
\end{lemma}
\begin{proof} 
We may suppose that~$|\xi|\geq \Gamma:= \left(c_0^{-1}(1+M)\right)^{\frac1p}$, otherwise there is nothing to show.  
Let $\theta$ be as in the (LS) condition. According to Lemma \ref{l.matrixtrick} below, there exists a $C^1$ map $\lambda \mapsto Z_\lambda$ from $[0,1]$ to the set of $d$-by-$d$ symmetric matrices such that 
$Z_0=X$, $Z_1\leq Y$ and
\begin{equation*}
\frac{d}{d\lambda} Z_\lambda= \frac{1}{(1+\kappa)\gamma} Z_\lambda^2.
\end{equation*}
Define
$$
f(\lambda)= -\tr\left(A(\xi_\lambda,x_\lambda)Z_\lambda\right) +H(\xi_\lambda, x_\lambda). 
$$
where $\xi_\lambda := (1-\lambda) \xi+\lambda \eta$ and $x_\lambda:=(1-\lambda)x+\lambda y$. Then $f(0)\leq M$ while 
$$
f(1) = -\tr\left(A(\eta,y)Z_1\right) +H(\eta, y) \geq -\tr\left(A(\eta,y)Y\right) +H(\eta, y)\geq M. 
$$
We deduce the existence of a largest value $\lambda\in [0,1]$ for which $f(\lambda)=M$. Note that $f'(\lambda)\geq 0$. As $|\xi_\lambda|\geq 1$ (since $|\xi|\geq 2$), this means that 
\begin{multline}\label{e.fprime}
-\tr\left(A\frac{Z_\lambda^2}{(1+\kappa)\gamma}\right) -(\tr\left(A_\xi Z_\lambda \right)\cdot(\eta-\xi))\\ -(\tr\left(A_x Z_\lambda \right)\cdot(y-x))
+(H_\xi\cdot(\eta-\xi)+ (H_x\cdot(y-x))\geq 0. 
\end{multline}
where $A$, $H$ and their derivatives are evaluated at $(\xi_\lambda,x_\lambda)$. Without loss of generality we assume that $|\xi|,|\eta|\geq 1$.  
Following now the computation from~\cite[Theorem 1.1]{LS0}, we obtain 
\begin{multline*}
(\tr\left(A_x Z_\lambda \right).(y-x)) \\
\leq |\sigma_x| \tr(AZ_\lambda^2 )^{\frac12} |x-y|  \leq \theta \tr\left(A\frac{Z_\lambda^2}{(1+\kappa)\gamma} \right)+ \theta^{-1}(1+\kappa) \gamma |\sigma_x|^2 |x-y|^2 
\end{multline*}
and, in the same way,  
\begin{multline*}
(\tr\left(A_\xi Z_\lambda \right).(\eta-\xi)) \leq \theta \tr\left(A\frac{Z_\lambda^2}{(1+\kappa)\gamma} \right)+ \theta^{-1}(1+\kappa)\gamma|\sigma_\xi|^2 |\xi-\eta|^2 \\ 
\leq 
\theta \tr\left(A\frac{Z_\lambda^2}{(1+\kappa)\gamma} \right)+ \theta^{-1}(1+\kappa)  \gamma  C_0  |\xi-\eta|^2
\end{multline*}
because $\sigma_\xi(\xi_\lambda,x)$ is bounded by $ C_0 $ since $|\xi_\lambda|\geq 1$ thanks to assumption \eqref{e.LipschA}. Multiply \eqref{e.fprime} by $(1+\kappa)\gamma$ to obtain, since $\gamma(x-y)=\xi$ and $(1+\gamma)|\xi-\eta|\leq \kappa$, 
$$
\theta(1-2\theta)\tr\left(A Z_\lambda^2\right)\leq (1+\kappa)\left(\kappa \theta |H_\xi|+ \theta |H_x||\xi|+(1+\kappa)|\sigma_x|^2 |\xi|^2 + C_0 (1+\kappa)\kappa^2  \right).
$$
Next as $f(\lambda)=M$, we also have 
$$
H-M\leq \tr\left(A Z_\lambda\right)\leq |\sigma| \left(\tr\left(A Z_\lambda^2\right)\right)^{\frac12} . 
$$
By \eqref{e.Hposhom} and the assumption $|\xi|\geq \Gamma$, we have $H-M\geq 0$. So
$$
(1+\kappa)^{-1} H^2- \kappa^{-1} M^2 \leq  |\sigma|^2 \left(\tr\left(A Z_\lambda^2\right)\right) 
$$
and therefore
\begin{multline*}
\theta(1-2\theta)(1+\kappa)^{-1} H^2 \leq \theta(1-2\theta)\kappa^{-1}M^2  \\+ (1+\kappa)|\sigma|^2\left(\kappa \theta |H_\xi|+ \theta |H_x||\xi|+(1+\kappa)|\sigma_x|^2 |\xi|^2 + C_0 (1+\kappa)\kappa^2  \right).
\end{multline*}
Using the bound $\theta\leq 1$, we get
\begin{multline*}
\theta(1-2\theta) H^2 \leq (1+\kappa)\kappa^{-1}M^2+4(1+\kappa) C_0 |\sigma|^2\kappa^2   \\
+ (1+\kappa)^2|\sigma|^2\left(\kappa \theta |H_\xi|+ \theta |H_x||\xi|+(1+\kappa)|\sigma_x|^2 |\xi|^2  \right).
\end{multline*}
Applying the (LS) condition, we obtain $|\xi|\leq \rho^{-1} (2\kappa^{-1}M^2+8 C_0 |\sigma|^2\kappa^2)$. 
\end{proof}

The following lemma, which was used in the argument above, is borrowed from~\cite{BC}. 

\begin{lemma}\label{l.matrixtrick}
Let $\gamma>0$ and $X,Y\in\R^{d\times d}$ be symmetric matrices satisfying
\begin{equation*} 
\begin{pmatrix} 
X &0 \\  0 & -Y  
\end{pmatrix} 
\leq \gamma \begin{pmatrix} 
I_d & -I_d \\ -I_d & I_d  
\end{pmatrix}.
\end{equation*}
Then there exists a $C^1$ map $\lambda \to Z_\lambda$ from $[0,1]$ to the set of $d$-by-$d$ symmetric matrices, such that, for every $\lambda\in [0,1]$,  
\begin{equation*}
X= Z_0\leq Z_\lambda\leq Y \qquad {\rm and }\qquad  \frac{d}{d\lambda} Z_\lambda= \gamma^{-1} Z^2_\lambda.
\end{equation*}
\end{lemma}
\begin{proof} 
We express the matrix inequality in the following form: for every $\xi,\eta \in\Rd$,
\begin{equation*}
X\xi\cdot \xi- Y\eta\cdot \eta\leq \gamma|\xi-\eta|^2.
\end{equation*}
Thus, for every $\eta\in\Rd$,
\begin{equation*}
\sup_{\xi\in\Rd} \left( X\xi\cdot \xi -\gamma|\xi-\eta|^2 \right)  \leq Y\eta \cdot \eta.
\end{equation*}
For every $\lambda\in [0,1)$, we have $X\leq \gamma I_d < \lambda\gamma I_d$ and thus, for every $\eta \in\Rd$,
\begin{equation*}
\sup_{\xi\in\Rd} \left( X\xi\cdot \xi-\lambda \gamma\left|\xi-\eta\right|^2\right) = X\left(I_d-\frac{X}{\lambda \gamma}\right)^{-1}\eta\cdot \eta.
\end{equation*}
Define
$$
Z_\lambda:= X \left( I_d-\frac{X}{\lambda \gamma}\right)^{-1}. 
$$
Observe that $Z_\lambda$ is an increasing family of symmetric matrices bounded above by~$Y$. It follows that $Z_0$ and $Z_1$ are well-defined and $Z_0=X$. Moreover,
\begin{equation*}
\frac{d}{d\lambda} Z_\lambda= \frac{X^2}{\gamma}   \left( I_d-\frac{\lambda X}{\gamma}\right)^{-2}= \gamma^{-1} Z_\lambda^2.\qedhere
\end{equation*}
\end{proof}

\subsection{Comparison for the metric problem}

We fix a nonempty, closed subset $S\subseteq\R^d$ and consider the metric problem
\begin{equation}\label{e.metricAppA1}
\left\{
\begin{aligned}
& -\tr\left(A(Dm,x)D^2m\right) +H(Dm, x)=\mu & \mbox{in} & \ \Rd\setminus S, \\
& m=0 &\mbox{on} & \ \partial S.
\end{aligned}
\right.
\end{equation}
We assume that the coefficients $(\sigma,H)$ of the equation belong to the closure $\overline \Omega$ of $\Omega$ for the local uniform convergence (as usual, $A=\frac12\sigma\sigma^T$). In view of the definition of $\Omega$, we note that $\sigma \in C^{0,1}(\partial B_1\times \R^d)$ and $H \in C^{0,1}( \Rd \times\Rd)$ satisfy the following
\begin{equation*}
\left| \sigma(e,x) \right| + \left| D_x \sigma(e,x) \right| + \left| D_\xi \sigma(e,x) \right| \leq C_0\qquad \mbox{a.e. in} \ \partial B_1\times \R^d,
\end{equation*}
\begin{equation*}
H(t\xi,x) = t^pH(\xi,x)  \quad \mbox{and} \quad
c_0\left|\xi\right|^p \leq H(\xi,x) \leq C_0\left|\xi\right|^p,
\end{equation*}
and
\begin{equation}\label{e.LipHAppen}
\left| D_x H( \xi,x) \right| + \left| \xi \right| \left| D_\xi H(\xi,x) \right| \leq C_0 \left| \xi \right|^p\qquad \mbox{a.e. in} \ \left(\R^d\setminus\{0\}\right)\times \R^d.
\end{equation}

In order to show the well-posedness of~\eqref{e.metricAppA1}, we first study the following constrained problem, for a parameter $K\geq 1$:
\begin{equation}\label{e.metricLiLiAppA1} 
\left\{
\begin{aligned}
& \max\left\{ -\tr\left(A(Dm,x)D^2m\right) +H(Dm, x)-\mu,\, |Dm|-K\right\} =0 & \mbox{in}  & \ \R^d\setminus S, \\
& m=0 & \mbox{on} & \ \partial S.
\end{aligned}
\right.
\end{equation}

We begin with a comparison principle for~\eqref{e.metricLiLiAppA1}. 

\begin{proposition}
\label{p.MPconstraintComp}
Fix $K\geq1$ and $\mu>0$. Suppose $m^1\in \USC\big(\overline{\Rd\setminus S}\big)$ and $m^2\in \LSC\big(\overline{\Rd\setminus S}\big)$ are respectively a subsolution and nonnegative supersolution of~\eqref{e.metricLiLiAppA1}. Then $m^1 \leq m^2$ in $\Rd\setminus S$.
\end{proposition}
\begin{proof}
Throughout the proof, the constants $C$ and $c$ may depend on $\sigma$, $H$, $\mu$, $K$ and vary in each occurrence. 

\smallskip

For $\delta>0$ small, let $\psi:\R\to \R$ be a smooth map which is increasing, concave, bounded above with $\psi(0)=0$ and $0\leq \psi'\leq 1-2\delta$ and $-( C_0  K^2)^{-1}\mu\delta \leq \phi''\leq 0$. We set 
$w(x)= \psi(m^1(x))$. Let us check that $w$ is a subsolution of \eqref{e.metricLiLiAppA1}, using the homogeneity of $H$ and $A$ with respect to the gradient variable, we have (in the viscosity sense)
\begin{align*}
\lefteqn{ -\tr\left(A(Dw(x),x)D^2w(x)\right) +H(Dw(x), x) } 
\qquad & \\
& \leq 
\psi'(m^1(x)) \left( -\tr\left(A(Dm^1(x),x)D^2m^1(x)\right) +H(Dm^1(x), x)\right) \\
& \qquad  - \psi''(m^1(x)) \tr\left(A(Dm^1(x),x)Dm^1(x)\otimes Dm^1(x)\right) \\
& \leq \mu(1-2\delta) + \psi''(m^1(x))  C_0  |Dm^1(x)|^2 \\
&  \leq \mu(1-\delta).
\end{align*}
Moreover $|Dw|-(1-2\delta)K= (1-2\delta)(|Dm^1|-K)\leq 0$. Thus $w$ satisfies 
\begin{equation}\label{e.metricLiLibisAppA1}
\max\left\{ -\tr\left(A(Dw,x)D^2w\right) +H(Dm, x)-(1-\delta)\mu, \, |Dm|-(1-2\delta)K\right\} \leq 0.
\end{equation}
We claim that $w\leq m^2$. We argue by contradiction and suppose instead that $M:=\sup_{\R^d} (w-m^2)>0$. Set
\begin{equation*}
M_{\alpha, \beta} := \sup_{x,y\in \R^d\setminus S} \left( w(x)-m^2(y) -\frac{|x-y|^4}{4\alpha} -\frac{\beta}{2}|y|^2 \right).
\end{equation*}
The supremum on the right side is evidently attained at some point $$(x_\beta,y_\beta)\in \big( \overline{\R^d\setminus S}\big) \times \big( \overline{\R^d\setminus S}\big).$$ By construction, $w(x)-m^2(y)$ is bounded above and thus 
$$
M_\alpha:= \lim_{\beta\to0} M_{\alpha,\beta}=  \sup_{x,y\in\R^d\setminus S} \left( w(x)-m^2(y) -\frac{|x-y|^4}{4\alpha}  \right)
$$
and 
\begin{equation}\label{e.limbetay}
\lim_{\beta \to 0} \beta |y_\beta |^2= 0.
\end{equation}
As $w$ is $K$-Lipschitz continuous (as it is a subsolution of \eqref{e.metricLiLibisAppA1}), we  have
\begin{equation}\label{e.estix-yAppenA1}
|x_\beta-y_\beta|\leq (K\alpha)^{\frac13} \leq C \alpha^{\frac13}.
\end{equation}
If $\alpha$ is sufficiently small, we have $x_\beta\not\in\partial S$ and $y_\beta\not\in \partial S$; indeed, if $x_\beta\in \partial S$, then 
$$
M_{\alpha,\beta}= w(x_\beta)-m^2(y_\beta) -\frac{|x_\beta-y_\beta|^4}{4\alpha} -\frac{\beta}{2}|y_\beta|^2\leq 0,
$$ 
 while, if $y_\beta\in \partial S$, then, by Lipschitz continuity of $w$,  
$$
M_{\alpha,\beta}
\leq w(y_\beta)+K|x_\beta-y_\beta| \leq K^{\frac43}\alpha^{\frac13} \; <\; M
$$
for $\alpha$ small. Both cases are impossible since $M_{\alpha,\beta}\geq M>0$. We deduce that the maximum point belongs to the interior: $(x_\beta,y_\beta) \in \left(\Rd \setminus S\right)\times \left(\Rd \setminus S\right)$. 

\smallskip

Fix $\eta>0$ to be selected below and apply the maximum principle for semicontinuous functions~\cite[Theorem 3.2]{CIL} to obtain symmetric matrices $X_\beta,Y_\beta\in\R^{d\times d}$ such that 
\begin{equation*}
(X_\beta,\xi_\beta)\in \overline{\mathcal J}^{2,+}w(x_\beta),   \qquad (Y_\beta,\xi_\beta -\beta y_\beta)\in  \overline{\mathcal  J}^{2,-} m^2(y_\beta)
\end{equation*}
and
\begin{equation*} 
-\left(\frac{1}{\eta}+ |M|\right)I_{2d} \leq  \begin{pmatrix} 
X_\beta &0 \\  0 & -(Y_\beta+\beta I_d)  
\end{pmatrix} 
\leq  M+ \eta M^2
\end{equation*}
where $\xi_\beta:= \alpha^{-1} \left| x_\beta - y_\beta \right|^2(x_\beta-y_\beta)$ and 
\begin{equation*}
M:=\frac{1}{\alpha} \begin{pmatrix} 
N & -N \\ -N & N  
\end{pmatrix}, \quad N:=\left|x_\beta-y_\beta\right|^2 I_d+2 (x_\beta-y_\beta)\otimes (x_\beta-y_\beta).
\end{equation*}
We select $\eta:=  |M|^{-1}$. As $N\leq C\left| x_\beta-y_\beta \right|^2I_d$,  we obtain with this choice: 
\begin{equation}\label{e.estiXYbetaMetric}
-C|M|  I_{2d} \leq 
  \begin{pmatrix} 
X_\beta &0 \\  0 & -(Y_\beta+\beta I_d)  
\end{pmatrix} 
\leq  C|M|  \begin{pmatrix} 
I_d & -I_d \\  -I_d & I_d  
\end{pmatrix} .
\end{equation} 
As $w$ is a subsolution to \eqref{e.metricLiLibisAppA1} and $m^2$ is a supersolution to \eqref{e.metricLiLiAppA1}, we have 
$$
\max\left\{ -\tr^*\left(A(\xi_\beta,x_\beta)X_\beta\right) +H(\xi_\beta, x_\beta)- (1-\delta)\mu\; ,\; |\xi_\beta|-(1-2\delta)K\right\}\leq 0
$$
and 
$$
\max\left\{  -\tr_*\left(A(\xi_\beta-\beta y_\beta,y_\beta)Y_\beta\right)
+H(\xi_\beta-\beta y_\beta, y_\beta)- \mu\; , \; 
|\xi_\beta-\beta y_\beta|-K\right\}\geq 0. 
$$
Note that, as $|\xi_\beta|-(1-2\delta)K\leq 0$, we cannot have $|\xi_\beta-\beta y_\beta|\geq K$ for $\beta$ small enough. Then the above inequalities can be rewritten as 
\begin{equation}\label{e.relrel1}
-\tr^*\left(A(\xi_\beta,x_\beta)X_\beta\right) +H(\xi_\beta, x_\beta) \leq (1-\delta)\mu
\end{equation}
and 
\begin{equation}\label{e.relrel2}
-\tr_*\left(A(\xi_\beta-\beta y_\beta,y_\beta)(Y_\beta-\beta I_d)\right)+H(\xi_\beta-\beta y_\beta, y_\beta)\geq  \mu- C_0 \beta. 
\end{equation}
We may also take $\beta$ sufficiently small that $(1-\frac13\delta)\mu<  \mu -  C_0 \beta$. 

Next we provide a lower bound on $\left| \xi_\beta\right|$. By \eqref{e.estiXYbetaMetric},we have  $Y_\beta-\beta I_d\geq - C\gamma I_d$, and thus by \eqref{e.relrel2}, we obtain
$$
C\gamma   + C(|\xi_\beta|^p+ |\beta y_\beta|^p) \geq \mu\left(1-\frac13\delta\right).
$$
By \eqref{e.limbetay} we can choose $\beta$ so small that $C|\beta y_\beta|^p< \frac13\mu\delta$, so that, by the definition of $\gamma$, 
$$
C \alpha^{-\frac13}|\xi_\beta|^{\frac23} + C |\xi_\beta|^p \geq \mu\left(1-\frac23\delta\right).
$$
This proves the bound below for $|\xi_\beta|$ for $\delta$ small enough:  
\begin{equation}\label{e.boundbelowxi}
|\xi_\beta|\geq C^{-1}\alpha^{\frac12} .
\end{equation} 
In particular, assuming again that $\beta$ is so small that $\beta |y_\beta|\leq C^{-1}\alpha^{\frac12} $, we can remove the ``$*$'' in \eqref{e.relrel1} and \eqref{e.relrel2}. 

Note that, by \eqref{e.estix-yAppenA1}, we have the estimates on $\gamma$ and $\xi_\beta$: 
\begin{equation}\label{e.Lischestito}
|\xi_\beta|\leq K \qquad {\rm and}\qquad \gamma\leq K^{\frac23}\alpha^{-\frac13}.
\end{equation} 

Our aim now is to compute the difference between \eqref{e.relrel1} and \eqref{e.relrel2}: since $H$ is locally Lipschitz continuous, uniformly with respect to $x$, there is a constant $C=C_K$ such that
 $$
 H(\xi_\beta-\beta y_\beta, y_\beta)-H(\xi_\beta, x_\beta)\leq C(\beta|y_\beta|+|x_\beta-y_\beta|)\; \leq \; 
C(\beta|y_\beta|+\alpha^{\frac13}) .
 $$
 Using the regularity of $A$ with respect to $\xi$ in \eqref{e.LipschA} and the bound below for $\xi_\beta$ in \eqref{e.boundbelowxi}, we have 
 \begin{multline*}
 \left| \tr\left(A(\xi_\beta-\beta y_\beta,y_\beta)(Y_\beta-\beta I_d)\right)- 
 \tr\left(A(\xi_\beta,y_\beta)(Y_\beta-\beta I_d)\right)\right| \\
 \leq  C  |\xi_\beta|^{-1} \beta |y_\beta| |Y_\beta-\beta I_d| 
 \leq C \alpha^{-\frac12}  \beta |y_\beta| |Y_\beta-\beta I_d| 
 \end{multline*}
On another hand, since $A=\frac12\sigma\sigma^T$ with $\sigma$ satisfying \eqref{e.LipHAppen} and using \eqref{e.estix-yAppenA1}, \eqref{e.estiXYbetaMetric} and \eqref{e.Lischestito},
we also have: 
$$
 \tr\left(A(\xi_\beta,y_\beta)(Y_\beta-\beta I_d)\right) -\tr\left(A(\xi_\beta,x_\beta)X_\beta\right) \leq C \gamma   |x_\beta-y_\beta|^2
 \leq C \alpha^{\frac13}.
 $$
 The difference between \eqref{e.relrel1} and \eqref{e.relrel2} then  yields to 
$$
- C \alpha^{-\frac12}  \beta |y_\beta| |Y_\beta-\beta I_d|  
 -C(\beta|y_\beta|+\alpha^{\frac13}) 
 \leq \mu(1-\delta) -\mu + C \beta.
$$
 We now let $\beta\to 0$ (with $\alpha$ remaining fixed): taking into account the fact that $\beta y_\beta\to 0$ and the bound on the matrix $Y_\beta$, we obtain 
$$
- C\alpha^{\frac13}  \leq -\delta\mu .
 $$ 
 Choosing $\alpha$ small then yields to a contradiction.  

\smallskip 

Thus far we have proved that, for any $\delta>0$ and any  smooth map $\psi:\R\to \R$ which is increasing, concave, bounded above map with $\psi(0)=0$ and $0\leq \psi'\leq 1-2\delta$ and $-( C_0  K^2)^{-1}\mu\delta \leq \phi''\leq 0$, we have $\psi(m^1)\leq m^2$. We can now let $\delta\to 0$ and $\psi$ tend to the identity to obtain the desired inequality: $ m^1\leq m^2$. 
\end{proof}

We next deduce a comparison principle for the metric problem~\eqref{e.metricAppA1} from Proposition~\ref{p.MPconstraintComp}. 

\begin{proposition}\label{p.compaAppenA1} 
Fix $\mu>0$. Suppose that $m^1 \in W^{1,\infty}_{\mathrm{loc}} \big(\overline{\Rd\setminus S}\big)$ is a subsolution of \eqref{e.metricAppA1} satisfying
\begin{equation*}
\esssup_{x\in\Rd\setminus S} \left| Dm^1(x) \right| <\infty\end{equation*}
and $m^2$ is a nonnegative lower semicontinuous supersolution of \eqref{e.metricAppA1}. Then $m^1\leq m^2$ in $\R^d\setminus S$. 
\end{proposition}
\begin{proof}
With $K:= \esssup_{x\in\Rd\setminus S} \left| Dm^1(x) \right|$, we see that $m^1$ is a subsolution of~\eqref{e.metricLiLiAppA1}. Clearly $m^2$ is a nonnegative supersolution of~\eqref{e.metricLiLiAppA1}. Therefore we obtain the result from Proposition~\ref{p.MPconstraintComp}.
\end{proof}

In the proof of Lemma \ref{l.seige}, we needed the following comparison principle in bounded domains. The proof is similar to that of Proposition \ref{p.compaAppenA1}, without even a need of  a penalization at infinity: for this reason we omit it.

\begin{proposition}\label{p.compaAppenA2} Let $U$ be a bounded open subset of $\R^d$. Assume that $m^1$ is a globally Lipschitz continuous subsolution and $m^2$ is a nonnegative lower semicontinuous supersolution of
\begin{equation*}
-\tr\left(A(Dm,x)D^2m\right) +H(Dm, x)=\mu \quad \mbox{in} \ U
\end{equation*}
such that $m^1\leq m^2$ on $\partial U$. Then $m^1\leq m^2$ in $U$.
\end{proposition}

\subsection{Well-posedness of the metric problem}

We now show that the metric problem has a unique globally Lipschitz solution, and that the Lipschitz regularity is independent of the set $S$. 
As in the previous subsection, we assume that the coefficients $(\sigma, H)$ of the equation belong to the closure $\overline \Omega$ of $\Omega$ with respect to the topology of local uniform convergence. 

\begin{theorem}
\label{t.LipschitzAppenA1}
Fix $L\geq 1$, $\mu\in (0,L]$ and a nonempty closed subset $S\subseteq\Rd$ satisfying~\eqref{e.IntBallPpt2}. Then there exists a unique Lipschitz continuous solution $m_\mu(\cdot,S)\in C^{0,1}_{\mathrm{loc}}\left(\overline {\Rd\setminus S}\right)$ of the metric problem~\eqref{e.metricAppA1}. Moreover, there exists a constant $C(\data,L)$ such that, for all $x,y\in \overline{\Rd\setminus S}$,
\begin{equation*}
\left| m_\mu(x,S)-m_\mu(y,S)\right| \leq C|x-y|.
\end{equation*}
Finally, $m_\mu(\cdot,S)$ is the maximal Lipschitz continuous subsolution of  \eqref{e.metricAppA1}. 
\end{theorem}
\begin{proof} 
For $K\geq 2$, we consider the problem \eqref{e.metricLiLiAppA1}. Clearly the zero function is a subsolution. As this problem possesses a comparison principle (Proposition~\ref{p.MPconstraintComp}), Perron's method (see \cite{CIL}) therefore provides the existence and the uniqueness of a solution $m^K$ to~\eqref{e.metricLiLiAppA1}, which is identified as the maximal subsolution of~\eqref{e.metricLiLiAppA1}. It is clear that $m^K$ is Lipschitz continuous with Lipschitz constant $K$.

\smallskip

Our aim is to show that, if $K$ is large enough, then the constraint on the gradient is never in force and thus $m^K$ is a solution to \eqref{e.metricAppA1}. We do so by proving that $m^K$ is actually Lipschitz continuous with a Lipschitz constant $K'\in \left[1, \frac12K\right]$ to be chosen below.  

\smallskip

We note that, by standard stability argument for equation \eqref{e.metricLiLiAppA1}, we may assume for convenience that~$(\sigma, H)\in\Omega$. Otherwise, we approximate $(\sigma, H)$ by coefficients in $\Omega$, prove the Lipschitz estimates for the approximate equation, obtaining an upper bound on the Lipschitz constant depending only on~$(\data,L)$, and then pass to the limit in the approximation. 

\smallskip

We now replace $m^K$ by a bounded function. For $\delta>0$ small, let $\psi:\R\to \R$ be a smooth map which is increasing, concave, bounded above with $\psi(0)=0$ and $0\leq \psi'\leq 1-2\delta$ and $-( C_0  K^2)^{-1}\mu\delta \leq \phi''\leq 0$. We set $w(x)= \psi(m^K)$. As in the proof of Proposition \ref{p.compaAppenA1}, one easily checks that $w$  satisfies the subsolution inequality \eqref{e.metricLiLibisAppA1}. The main part of the proof consists in checking that for suitable choices of the parameters we have $w(x)\leq m^K(y)+K'|x-y|$ for  $x,y\in \Rd\setminus S$. 

\smallskip

We argue by contradiction. Assume that
\begin{equation*}
M:= \sup_{x,y\in\Rd\setminus S} \left( w(x)-m^K(y)-K'|x-y|\right) > 0.
\end{equation*}
For $\beta>0$ small, set
\begin{equation*}
M_\beta := \sup_{x,y \in \R^d\setminus S} \left( w(x)-m^K(y)-K'|x-y|-\frac{\beta}{2}|y|^2 \right).
\end{equation*}
The supremum is attained at some $(x_\beta, y_\beta)\in\big(\overline{\R^d\setminus S}\big)\times \big(\overline{\R^d\setminus S}\big)$. Note that
\begin{multline*}
0<M\leq M_\beta \leq m^K(x_\beta)-m^K(y_\beta)-K'|x_\beta-y_\beta| \\
\leq (K-K')\left|x_\beta-y_\beta\right| \leq \frac12K \left|x_\beta-y_\beta\right|,
\end{multline*}
which gives the lower bound
\begin{equation}\label{e.BdBelowx-yAppn1}
\left|x_\beta-y_\beta\right|\geq \frac{2M}{K}>0.
\end{equation}
In particular, $x_\beta \neq y_\beta$.
As in the proof of Proposition~\ref{p.compaAppenA1}, keeping the other parameters fixed, we have that $\lim_{\beta \to 0} \beta y_\beta = 0$.

\smallskip

We first assume that  $x_\beta\notin \partial S$ and $y_\beta\notin \partial S$. Fix $\eta>0$.
By the maximum principle for semicontinuous functions~\cite[Theorem 3.2]{CIL}, there exist symmetric matrices $X_\beta,Y_\beta\in \R^{d\times d}$ such that 
\begin{equation*}
\left\{
\begin{aligned}
& (X_\beta,\xi_\beta)\in \overline{\mathcal J}^{2,+}w(x_\beta), \\
& (Y_\beta,\xi_\beta -\beta y_\beta)\in  \overline{\mathcal  J}^{2,-} m^K(y_\beta)
\end{aligned}
\right.
\end{equation*}
and 
\begin{equation*}
-\left(\frac{1}{\eta}+ |M|\right)I_{2d} \leq  \begin{pmatrix} 
X_\beta &0 \\  0 & -(Y_\beta+\beta I_d)  
\end{pmatrix} 
\leq  M+ \eta M^2
\end{equation*}
where we denote $\gamma: =K'|x_\beta-y_\beta|^{-1}$, $\xi_\beta:=K'(x_\beta-y_\beta)/|x_\beta-y_\beta|$ and 
\begin{equation*}
M:=\gamma \begin{pmatrix} 
N & -N \\ -N & N  
\end{pmatrix}, \quad N:=I_d- \frac{x_\beta-y_\beta}{|x_\beta-y_\beta}\otimes \frac{x_\beta-y_\beta}{|x_\beta-y_\beta}.
\end{equation*}
Note that $M^2= 2\gamma^2 M$. Select $\eta:=  (2\gamma)^{-1} \kappa$, where $\kappa$ is the parameter in the~(LS) condition~\eqref{e.LScoercivity}. As $N\leq I_d$,  we obtain
\begin{equation*}
-C\gamma  \left(1+\frac1\kappa\right)I_{2d} \leq 
  \begin{pmatrix} 
X_\beta &0 \\  0 & -(Y_\beta+\beta I_d)  
\end{pmatrix}  \leq  (1+\kappa)\gamma  \begin{pmatrix} 
I_d & -I_d \\  -I_d & I_d  
\end{pmatrix}.
\end{equation*}
Using the equations satisfied by $w$ and $m^K$, we have 
\begin{equation*}
\max\left\{ -\tr^*\left(A(\xi_\beta,x_\beta)X_\beta\right) +H(\xi_\beta, x_\beta)-(1-\delta) \mu, \,|\xi_\beta|-(1-\delta)K\right\} \leq 0
\end{equation*}
and 
\begin{equation*}
\max\left\{-\tr_*\left(A(\xi_\beta-\beta y_\beta,y_\beta)Y_\beta\right) +H(\xi_\beta-\beta y_\beta, y_\beta)-\mu,\,  |\xi_\beta|-K\right\} 
\geq 0. 
\end{equation*}
As $0<|\xi_\beta|=K'\leq \frac12K$ and $\beta y_\beta$ is small, the above inequalities yield, for sufficiently small~$\beta$,
\begin{equation}\label{e.relrel21}
-\tr\left(A(\xi_\beta,x_\beta)X_\beta\right) +H(\xi_\beta, x_\beta) \leq (1-\delta) \mu
\end{equation}
and 
\begin{equation}\label{e.relrel22}
-\tr\left(A(\xi_\beta-\beta y_\beta,y_\beta)(Y_\beta-\beta I_d)\right) +H(\xi_\beta-\beta y_\beta, y_\beta) \geq \mu- C_0  \beta. 
\end{equation}
According to~\eqref{e.BdBelowx-yAppn1}, $\gamma\leq K'K/2M$. For $\beta$ so small that $(1+\gamma)\beta |y_\beta|\leq \kappa$, Lemma~\ref{l.LSTechniPt} gives the existence of $R(\data,L)$ such that $|\xi_\beta|\leq R$. As $|\xi_\beta|=K'$, we obtain a contradiction if $K'>R$. 

\smallskip

We have proved that, if $R<K'<\frac12K$, then the supremum in the definition of $M_\beta$, if positive, can only be achieved by $(x_\beta,y_\beta)$ if $x_\beta\in \partial S$ or $y_\beta\in\partial S$. If $x_\beta\in \partial S$, then 
$$
M_\beta=w(x_\beta)-m^K(y_\beta) -K'|x_\beta-y_\beta|-\frac{\beta}{2}|y_\beta|^2 \leq0,
$$ 
a contradiction.  Consider the case that $y_\beta\in\partial S$. As $S$ satisfies the interior ball condition of radius $1$, there exists $e\in \partial B_1$ such that the unit ball centered at $z_\beta:= y_\beta-e$ is contained in $S$. Then, for $r := c_0^{-1}(L+ C_0  (d-1))$, the map $w^2(x):= r(|x-z_\beta|-1)_+$ is a supersolution to \eqref{e.metricAppA1}, and thus to \eqref{e.metricLiLiAppA1}. 
By the comparison principle, $w(x_\beta)\leq m^K(x_\beta)\leq r(|x_\beta-z_\beta|-1)_+$. We deduce that 
\begin{multline*}
0<M_{\beta}= w(x_\beta)-m^K(y_\beta) -K'|x_\beta-y_\beta|-\frac{\beta}{2}|y_\beta|^2  \\
\leq r(|x_\beta-z_\beta|-1)_+-K'|x_\beta-y_\beta|\leq (r-K') |x_\beta-y_\beta|\leq 0,
\end{multline*}
provided that we also have  $K'\geq r$, which is another contradiction. 

\smallskip

In conclusion, we have shown that, if $C(\data,L):=R\vee r < K' <\frac12K$, then we have $w(x):= \psi(m^K(x))\leq m^K(y)+K'|x-y|$. As $\delta$ and $\psi$ are arbitrary, we can let $\psi$ tend to the identity to obtain that $m^K$ is Lipschitz with constant~$K' < \frac12K$. Thus the gradient constraint is never in force and $m^K$ is a Lipschitz continuous solution to~\eqref{e.metricAppA1}. By the comparison principle, it is the unique and the maximal one. 
\end{proof}

We next summarize several properties of the metric problem associated to a nonempty closed set $S\subset \R^d$ satisfying the uniform ball condition
\begin{equation} \label{e.IntBallPpt2}
S = \bigcup_{\overline B_1(x) \subseteq S} \overline B_1(x).
\end{equation}
As above, we consider coefficients $(\sigma, H)$ belonging to the $\overline \Omega$ of $\Omega$ with respect to the topology of local uniform convergence. 

\begin{lemma}\label{l.EstiMetric} 
Fix $L\geq 1$, $\mu\in (0,L]$ and closed subsets $S, S'\subseteq\Rd$ satisfying the uniform ball condition~\eqref{e.IntBallPpt2}. Then there exist $0< l_\mu\leq L_\mu$,  $C(\data,L)\geq1$ and $c(\data,L)\in (0,1]$ satisfying
\begin{equation}\label{lmuLmu}
 c\mu \leq l_\mu \leq L_\mu \leq C
\end{equation}
such that the following holds:
\begin{itemize}
\item[(i)] Lipschitz estimates: for every $x\in\Rd\setminus \Rd\setminus S$ and $y\in \Rd\setminus S'$,
\begin{equation}\label{e.Lipschpc}
\left|m_\mu(x,S)-m_\mu(y,S')\right|\leq L_\mu\left(|x-y|+ \dist_H(S,S')\right).
\end{equation}

\item[(ii)] Estimate on the level-sets:  for every $0\leq s\leq t$, 
\begin{equation}\label{e.movefron}
\dist_H\left( \{m_\mu(\cdot,S)\leq s\} , \{m_\mu(\cdot,S)\leq t\} \right) \leq \frac{1}{l_\mu} |s-t|+2.
\end{equation}
\item[(iii)] Growth property: for every $x\in\Rd$, 
\begin{equation}\label{control2}
l_\mu \dist(x,S)-2 \leq m_\mu(x,S ) \leq L_\mu \dist(x,S).
\end{equation}
When the target is the half space, i.e., ~$S=\H^-_e$ for some $e\in\partial B_1$, then we have, for every $x\in \H^+_e$, 
\begin{equation}\label{RoughGrowthmmue}
\left(\frac{\mu}{C_0} \right)^{\frac1p} x\cdot e \leq m_\mu(x,\H^-_e)\leq  \left(\frac{\mu}{c_0} \right)^{\frac1p}  x\cdot e.
\end{equation}
\end{itemize}
\end{lemma}

\begin{proof} The argument is standard and mostly borrowed from the proof of Proposition 2.2 in \cite{AC}. We only explain the differences. 

(i) Theorem \ref{t.LipschitzAppenA1} states that the $m_\mu(\cdot,S)$ is uniformly Lipschitz continuous. Then one can show, exactly as  for assertion (ii) in \cite[Proposition 2.2]{AC} that 
$$
|m_\mu(x,S)-m_\mu(x,S')|\leq L_\mu \dist_H(S,S').
$$

(ii) For the estimate on the level-sets, we follow again the proof of Proposition 2.2 of \cite{AC}: let us set $K:=\{m_\mu(\cdot, S)\leq s\}$ and let $\xi:\R^d\to \R$ be a standard mollification kernel. We denote, for $\ep>0$, $\xi_\ep(y) := \ep^{-d} \xi(y/\ep)$ and set $w_\ep := l_\mu \dist(\cdot, K)* \xi_\ep +s-Cl_\mu\ep$. Since  $\dist(\cdot, K)$ is Lipschitz continuous, we have $|Dw_\ep|\leq l_\mu$, $|D^2w_\ep|\leq l_\mu/\ep$ and $|w_\ep-\dist(\cdot, K)|\leq Cl_\mu\ep$. Then, for $l_\mu=c\mu$ with $c>0$ small enough, the function $w_\ep$ is a Lipschitz continuous subsolution of 
$$
-\tr\left(A(Dm,x)D^2m\right) +H(Dm, x)=\mu \qquad \mbox{in} \  \R^d\setminus \{m_\mu(\cdot,S)\leq s\}
$$
with $w_\ep\leq s\leq m_\mu(\cdot,S)$ on $\partial\{m_\mu(\cdot,S)\leq s\}$. By comparison (Proposition \ref{p.compaAppenA1}) we have therefore 
$$
l_\mu \dist(y, K) +s-2Cl_\mu\ep\leq w_\ep(y) \leq m_\mu(y,S) \qquad \mbox{in} \  \R^d\setminus \{m_\mu(\cdot,S)\leq s\}. 
$$
So, for $\ep$ small enough,
$$
 \{m_\mu(\cdot,S)\leq s\}\subset  \{m_\mu(\cdot,S)\leq t\}\subset  \{m_\mu(\cdot,S)\leq s\}+B_{l_\mu^{-1}(t-s)+2}.
 $$
 
The last two points are straightforward:  for (iii), the first inequality in  \eqref{control2} is a consequence of (ii) for $s=0$ while the second one holds by Lipschitz estimates. When the target is a plane, the map $x\to  (\mu/C_0)^{\frac1p} x\cdot e$ is a subsolution  to \eqref{e.planarmetric} by the homogeneity of $H$, while the map $x\to (\mu/C_0)^{\frac1p}  x\cdot e$ is a supersolution: this yields \eqref{RoughGrowthmmue} by comparison. 
 \end{proof}

\subsection{The approximate corrector problem}

In this subsection, we briefly discuss the approximate corrector problem: for $\delta > 0$ and $\xi\in\Rd$, we consider
\begin{multline} \label{e.approxcorrApp}
\delta v^\delta (\cdot,\xi) - \tr\left( A\left(\xi+Dv^\delta(\cdot,\xi),x\right)D^2v^\delta (\cdot,\xi) \right) + H\left(\xi+Dv^\delta(\cdot,\xi),x\right) \\
 = 0 \quad \mbox{in} \ \Rd.
\end{multline}
Here we assume that the coefficients $(\sigma, H)$ belong to $\Omega$. We summarize the facts we need in the following proposition, most of which is essentially known. The well-posedness of mean curvature-type equations has been well-understood since the pioneering papers \cite{CGG,ES,GGIS}. The only additional difficulty here comes the singularity of the matrix $A(\cdot+\xi,x)$, which is at~$-\xi$ instead of at~$0$. This can be treated as in the proof of Proposition~\ref{p.MPconstraintComp}. The $L^\infty$ bound~\eqref{e.RoughBoundVdelta} on $v^\delta$ is a straightforward consequence of assumption \eqref{e.Hposhom}. 
The Lipschitz estimate (in space) under coercivity condition \eqref{e.LScoercivity} for equation \eqref{e.approxcorrApp} has been established in \cite{LS0}, and subsequently discussed in \cite{CM}. The proof of \cite{LS0} relies on approximation by smooth solutions. The (suitably adapted) proof of Theorem~\ref{t.LipschitzAppenA1} above provides an alternative argument. Therefore, we only check the regularity of $\delta v^\delta$ with respect to $\xi$.

\begin{proposition} For any $\delta>0$ and $\xi\in \R^d$, equation \eqref{e.approxcorrApp} has a unique solution~$v^\delta(\cdot,\xi)\in W^{1,\infty}(\Rd)$ which satisfies 
\begin{equation}\label{e.RoughBoundVdelta}
-\frac{C_0|\xi|^p}{\delta}\leq v^\delta(x,\xi)    \leq -\frac{c_0|\xi|^p}{\delta}. 
\end{equation}
Moreover, for every $R\geq1$, there exists $C(\data,R)\geq 1$ such that, for every $x,y\in \R^d$ and $\xi,\eta\in B_R \setminus \{0\}$,
\begin{equation}\label{e.VdeltaHolder}
\left| v^\delta(x,\xi)-v^\delta(y,\eta)\right|\leq C\left(|\xi|+ |\eta|\right)|x-y|+ \frac{C}\delta\left(|\xi|\wedge |\eta| \right)^{-\frac{2p}7} \left|\xi-\eta\right|^{\frac27}.
\end{equation}
\end{proposition}
\begin{proof} 
As mentioned before the statement of the proposition, we only verify the regularity of $\delta v^\delta$ with respect to $\xi$. Throughout, we denote by $C$ and $c$ constants which depends on $(\data,R)$ and may vary in each occurrence. Given $\xi, \eta\in B_R\setminus \{0\}$, we denote for simplicity $v_1:= v^\delta(\cdot,\xi)$ and $v_2:= v^\delta(\cdot, \eta)$. Fix small parameters $\alpha,\beta>0$ to be chosen below, and consider 
$$
M_\beta:=  \sup_{(x,y) \in \R^d\times \R^d} \left(  v_1(x)-v_2(y) + \eta\cdot(x-y)-\frac{|x-y|^4}{4\alpha} -\beta |x|^2 \right)
$$
and denote by $(x_\beta,y_\beta)$ a maximum point of the problem. 
We shall send $\beta\to 0$ while keeping $\alpha>0$ fixed. Note that
\begin{equation*} \label{}
\lim_{\beta\to 0} M_\beta = \sup_{(x,y) \in \R^d\times \R^d} \left( v_1(x)-v_2(y) + \eta\cdot(x-y)-\frac{1}{4\alpha}|x-y|^4 \right)  \geq  \sup_{\R^d} (v_1-v_2)
\end{equation*}
and
\begin{equation*} \label{}
\limsup_{\beta\to 0}\beta |x_\beta|^2= 0. 
\end{equation*}
In view of the Lipschitz estimate on $v_2$, we have 
\begin{equation}\label{xbeta-ybeta}
|x_\beta-y_\beta|\leq C \alpha^{\frac13}.
\end{equation}
We set $\xi_\beta:= \alpha^{-1}|x_\beta-y_\beta|^2(x_\beta-y_\beta)$ and  $\hat \xi_\beta=\xi_\beta/|\xi_\beta|$ if $x_\beta\neq y_\beta$  ($\hat \xi_\beta=0$ otherwise). 

\smallskip

Fix $\gamma>0$ to be selected. By the maximum principle for semicontinuous functions~\cite[Theorem 3.2]{CIL}, there exist symmetric matrices $X_\beta,Y_\beta \in \R^{d\times d}$ such that 
\begin{equation*} \label{}
\left\{ 
\begin{aligned}
& (X_\beta,\xi_\beta-\eta+2\beta x_\beta)\in \overline{\mathcal J}^{2,+}v_1(x_\beta), \\
& (Y_\beta,\xi_\beta-\eta )\in  \overline{\mathcal  J}^{2,-} v_2(y_\beta)
\end{aligned}
\right.
\end{equation*}
with 
\begin{equation*}
-\left(\frac1\gamma +|M|\right) I_{2d} \leq 
  \begin{pmatrix} 
X_\beta+2\beta I_d & 0 \\  0 & -Y_\beta 
\end{pmatrix}  \leq  M+ \gamma M^2,
\end{equation*}
where
\begin{equation*}
M:=\frac1\alpha \begin{pmatrix} 
N & -N \\  -N & N  
\end{pmatrix}, \quad N:= |x_\beta-y_\beta|^2 I_d+2\left( x_\beta-y_\beta \right) \otimes \left(x_\beta-y_\beta  \right).
\end{equation*}
We choose $\gamma:=  \alpha|x_\beta-y_\beta|^{-2}$ if $x_\beta\neq y_\beta$ and $\gamma =1$ otherwise. As $N\leq C|x_\beta-y_\beta |^2I_d$,  we obtain:
\begin{equation}\label{e.estiXYbeta}
-C\frac{|x_\beta-y_\beta|^2}{\alpha} I_{2d} \leq 
 \begin{pmatrix} 
X_\beta+2\beta I_d &0 \\  0 & -Y_\beta 
\end{pmatrix}
 \leq  \frac{C|x_\beta-y_\beta|^2}{\alpha} \begin{pmatrix} 
I_d & -I_d \\  -I_d & I_d  
\end{pmatrix} 
\end{equation} 
where, in view of \eqref{xbeta-ybeta},  
\begin{equation*} \label{}
\frac{C|x_\beta-y_\beta|^2}{\alpha}\leq C\alpha^{-\frac13}.
\end{equation*}

\smallskip

Using the equations for $v_1$ and $v_2$, we have 
\begin{equation}\label{e.vdelta1App}
\delta v_1(x_\beta) -\tr^*\left(A\left(\xi_\beta+\xi-\eta+\beta x_\beta,x_\beta\right)X_\beta\right)\\ +H\left(\xi_\beta+\xi-\eta+\beta x_\beta, x_\beta\right)\leq 0
\end{equation}
and 
\begin{equation}\label{e.vdelta2App}
\delta v_2(y_\beta) -\tr_*\left(A\left(\xi_\beta,y_\beta\right)Y_\beta\right) +H\left(\xi_\beta, y_\beta\right)\geq 0. 
\end{equation}
Plugging \eqref{e.RoughBoundVdelta} into \eqref{e.vdelta2App} we find, since by \eqref{e.estiXYbeta} we have $|Y_\beta|\leq  C|x_\beta-y_\beta|^2/\alpha$: 
$$
C \frac{|x_\beta-y_\beta|^2}{\alpha} + C_0\frac{|x_\beta-y_\beta|^3}{\alpha} \geq -\delta v^\ep(y_\beta)\geq c_0|\eta|^p.
$$
This provides lower bound for $|x_\beta-y_\beta|$ and for $|\xi_\beta|$: 
$$
|x_\beta-y_\beta|\geq C^{-1} |\eta|^{\frac{p}2} \alpha^{\frac12}\qquad {\rm and}\qquad |\xi_\beta|\geq C^{-1} |\eta|^{\frac{3p}2} \alpha^{\frac12}.
$$
From now on we assume that $|\xi-\eta|\leq (3C)^{-1} |\eta|^{3p/2} \alpha^{\frac12}$ and that $\beta$ is so small that 
$\beta |x_\eta |\leq (3C)^{-1} |\eta|^{\frac{3p}2} \alpha^{\frac12}$. With this lower bound, we can simplify~\eqref{e.vdelta1App}: using
$$
\left|\xi_\beta+\xi-\eta+\beta x_\beta\right|\geq  (3C)^{-1} |\eta|^{\frac{3p}2} \alpha^{\frac12}, 
$$
we obtain, using the regularity of $A$ and \eqref{e.estiXYbeta}, 
\begin{multline*}
\left| A\left(\xi_\beta+\xi-\eta+\beta x_\beta,x_\beta\right)X_\beta
-
A\left(\xi_\beta,x_\beta\right)(X_\beta+2 \beta I_d))\right|  \\
\leq  C_0  \beta + C  |\eta|^{-\frac{3p}2} \alpha^{-\frac12} \left| \xi-\eta+\beta x_\beta\right||X_\beta| \\
\leq C|\eta|^{-\frac{3p}2} \alpha^{-\frac56} \left( \beta(1+|x_\beta|) + |\xi-\eta| \right).
\end{multline*}
On the other hand, by Lipschitz estimates,   
$$
H\left(\xi_\beta, y_\beta\right)
-
H\left(\xi_\beta+\xi-\eta+\beta x_\beta, x_\beta\right) \\
\leq 
C (|\xi-\eta|+\beta |x_\beta|+ |x_\beta-y_\beta|)
$$
while, by the structural assumptions on $A=\frac12\sigma\sigma^T$ and \eqref{e.estiXYbeta},
$$
\tr\left(A\left(\xi_\beta,y_\beta\right)Y_\beta\right) 
-
\tr\left(A\left(\xi_\beta,x_\beta\right)(X_\beta+2 \beta I_d)\right) 
\leq  C\alpha^{\frac13}. 
$$
Computing the difference between \eqref{e.vdelta1App} and \eqref{e.vdelta2App} and collecting the above inequalities we obtain (neglecting the lower order terms)
$$
\delta v_1(x_\beta) -\delta v_2(y_\beta)  \leq 
C|\eta|^{-\frac{3p}2} \alpha^{-\frac 56} \left( \beta(1+|x_\beta|) + |\xi-\eta| \right)+
 C\alpha^{\frac13}. 
$$
Therefore 
\begin{multline*}
M_\beta \leq  \delta v_1(x_\beta) -\delta v_2(y_\beta) +\eta\cdot(x_\beta-y_\beta)
\\ \leq 
C|\eta|^{-\frac{3p}2} \alpha^{-\frac 56} \left( \beta(1+|x_\beta|) + |\xi-\eta| \right)
+ C\alpha^{\frac13}. 
\end{multline*}
We now send $\beta\to 0$ to obtain 
\begin{equation*} \label{}
\sup_{\R^d} (v_1-v_2) \leq C|\eta|^{-\frac{3p}2} \alpha^{-\frac 56}  |\xi-\eta| + C\alpha^{\frac13}. 
\end{equation*}
Taking $\alpha= |\eta|^{-\frac{9p}7} |\xi-\eta|^{\frac67}$, we obtain
\begin{equation*} \label{}
\sup_{\R^d} (v_1-v_2) \leq C |\eta|^{-\frac{2p}7} |\xi-\eta|^{\frac27}.  \qedhere
\end{equation*}
 \end{proof}

\subsection{The time-dependent initial-value problem}

We briefly discuss the well-posedness of the problem
\begin{equation}\label{e.pdeA}
\left\{ 
\begin{aligned}
& \partial_t u - \tr\left( A\left(Du,x \right)D^2u \right) + H\left(Du,x \right) = 0 &  \mbox{in} & \ \Rd\times (0,T], \\
& u(\cdot,0) = g(\cdot) & \mbox{on} & \ \Rd,
\end{aligned}
\right.
\end{equation}
and its rescaled version
\begin{equation}\label{e.pdeAresc}
\left\{ 
\begin{aligned}
& \partial_t u^\ep -\ep \tr\left( A\left(Du^\ep,\frac x\ep \right)D^2u^\ep \right) + H\left(Du^\ep,\frac x\ep \right) = 0 &  \mbox{in} & \ \Rd\times (0,T], \\
& u^\ep(\cdot,0) = g(\cdot) & \mbox{on} & \ \Rd,
\end{aligned}
\right.
\end{equation}
Here we fix coefficients $(\sigma, H)$ belonging to $\Omega$. 

\begin{theorem} 
\label{e.WPcauchy}
For any~$g\in \BUC(\Rd)$, the initial-value problem~\eqref{e.pdeA} has a unique solution~$u \in \BUC(\Rd\times [0,T])$. If, in addition $g\in C^{1,1}(\R^d)$, then $u\in C^{0,1}(\Rd)$ and there exists a constant $L\left(\data, \|g\|_{C^{1,1}}\right)$ such that, for every $x,y\in \Rd$ and $s,t\in[0,T]$, 
\begin{equation*}
\left| u(x,t) - u(y,s) \right| \leq L \left( |x-y| + |t-s| \right).
\end{equation*}
\end{theorem}

We do not give the proof of Theorem~\ref{e.WPcauchy} here. Well-posedness and Lipschitz estimates for \eqref{e.pdeA} are proved in \cite{GGIS} and \cite{LS0}, respectively, and can also be obtained by modifying the arguments given here for the metric problem. 

\smallskip

Similar results for~\eqref{e.pdeAresc} can be obtained by applying Theorem~\ref{e.WPcauchy} and using the scaling $u^\ep(x,t)= \ep u(\frac x \ep, \frac t \ep)$. In particular, for a $C^{1,1}$ initial condition, the solution $u^\ep$ satisfies the same Lipschitz estimate: for every $x,y\in \R^d$ and $s,t\in[0,T]$,  
\begin{equation*}
\left| u^\ep(x,t) - u^\ep (y,s) \right| \leq L \left( |x-y| + |t-s| \right).
\end{equation*}

\smallskip

For the reader's convenience, let us briefly recall how Lipschitz bounds are established in this framework. Because of the $C^{1,1}$ regularity of $g$, the maps $(x,t)\to g(x)+Ct$ and $(x,t)\to g(x)-Ct$ are respectively super- and subsolutions of \eqref{e.pdeA} (where $C$ depends on $\data$ and $\|g\|_{C^{1,1}}$ only). This gives the estimate $|u(x,t)-g(x)|\leq Ct$. As the coefficients of the equation are  independent of time, the comparison principle then implies 
\begin{equation*}
\sup_{x\in\Rd} \left| u(x,t) - u(x,s) \right| \leq \sup_{x\in\Rd} \left| u(x,t-s) - g(x) \right| \leq C |t-s|,
\end{equation*}
and we have the Lipschitz regularity in time and a bound on $|\partial_t u|$. To get the Lipschitz regularity in space, one now argues as for stationary equations (see, e.g., the proof of Theorem \ref{t.LipschitzAppenA1}), using the bound on $|\partial_t u|$ and the (LS) condition~\eqref{e.LScoercivity}.

\bibliographystyle{plain}
\bibliography{meancurv}
\end{document}